\documentclass[a4paper,twoside]{article}

\usepackage{amsthm}
\usepackage{amssymb}
\usepackage{amsmath}
\usepackage{amscd}
\usepackage{makeidx}
\usepackage{fancyhdr}

\newcommand{\A}{\mathbb A}
\newcommand{\F}{\mathbb F}
\newcommand{\N}{\mathbb N}
\newcommand{\Q}{\mathbb Q}
\newcommand{\Z}{\mathbb Z}
\newcommand{\R}{\mathbb R}
\newcommand{\C}{\mathbb C}
\newcommand{\G}{\mathbb G}

\newcommand{\calA}{{\cal A}}
\newcommand{\calB}{{\cal B}}
\newcommand{\calC}{{\cal C}}
\newcommand{\calD}{{\cal D}}
\newcommand{\calE}{{\cal E}}

\newcommand{\calH}{{\cal H}}
\newcommand{\calI}{{\cal I}}
\newcommand{\calJ}{{\cal J}}
\newcommand{\calK}{{\cal K}}
\newcommand{\calL}{{\cal L}}
\newcommand{\calM}{{\cal M}}
\newcommand{\calN}{{\cal N}}
\newcommand{\calO}{{\cal O}}
\newcommand{\calR}{{\cal R}}
\newcommand{\calS}{{\cal S}}
\newcommand{\calU}{{\cal U}}
\newcommand{\calV}{{\cal V}}
\newcommand{\calX}{{\cal X}}
\newcommand{\calY}{{\cal Y}}

\def\CC{{\bf C}}
\newcommand{\kk}{{\mbox{\boldmath $k$}}}

\def\wtA{\widetilde{A}}
\def\wtB{\widetilde{B}}
\def\wtS{\widetilde{S}}
\def\wtf{\widetilde{f}}
\def\wtXi{\widetilde{\Xi}}

\def\tif{\tilde{f}}
\def\tixi{\tilde{\xi}}

\newcommand{\wt}[1]{\widetilde{#1}}
\newcommand{\wh}[1]{\widehat{#1}}

\newcommand{\whO}{\widehat{\calO}}

\def\uw{{\underline{w}}}

\def\al{{\alpha}}
\def\be{{\beta}}
\def\de{{\delta}}

\def\De{{\Delta}}
\def\ep{{\epsilon}}
\newcommand{\vep}{{\varepsilon}}
\def\ga{{\gamma}}
\def\Ga{{\Gamma}}
\newcommand{\ka}{{\kappa}}
\def\la{{\lambda}}
\def\La{{\Lambda}}

\def\si{{\sigma}}
\def\Si{{\Sigma}}
\def\th{{\theta}}
\def\vth{{\vartheta}}
\def\Th{{\Theta}}
\def\ze{{\zeta}}
\def\om{{\omega}}
\def\Om{{\Omega}}
\newcommand{\oom}{\underline{\omega}}

\newcommand{\gotA}{\mathfrak A}

\newcommand{\gotC}{\mathfrak C}

\newcommand{\gote}{\mathfrak e}
\newcommand{\gotH}{\mathfrak H}
\newcommand{\gotm}{\mathfrak m}
\newcommand{\gotn}{\mathfrak n}
\newcommand{\gotoo}{\mathfrak o}
\newcommand{\gotp}{\mathfrak p}
\newcommand{\gotq}{\mathfrak q}
\newcommand{\gotR}{\mathfrak R}
\newcommand{\gots}{\mathfrak s}
\newcommand{\gotS}{\mathfrak S}
\newcommand{\gotg}{\mathfrak g}
\newcommand{\gotl}{\mathfrak l}

\def\Aut{{\rm Aut}}

\def\CM{{\rm CM}}

\def\Dop{D^{\mathrm{op}}}
\def\End{{\rm End}}
\def\Ext{\mathrm{Ext}}
\def\Gal{{\rm Gal}}
\def\GL{{\rm GL}}
\def\GO{{\rm GO}}
\def\GSp{{\rm GSp}}
\def\gr{\mathrm{gr}}
\def\Hom{{\rm Hom}}
\newcommand{\Imm}{\mathrm{Im}}
\def\KS{{\rm KS}}
\def\Lie{{\rm Lie}}

\def\M{{\rm M}}
\newcommand{\orn}{{\rm or}}
\def\SL{{\rm SL}}
\def\SO{{\rm SO}}
\def\Sp{{\rm Sp}}
\newcommand{\Spec}{\mathrm{Spec}}
\newcommand{\Spf}{\mathrm{Spf}}
\def\Sym{\mathrm{Sym}}
\newcommand{\tra}{{\rm tr}}

\newcommand{\Gauss}{Gau{\ss}}

\newcommand{\bs}{\backslash}

\newcommand{\derham}[1]{\calH^{#1}_{\mathrm{dR}}}
\newcommand{\Derham}[1]{H^{#1}_{\mathrm{dR}}}

\newcommand{\dual}[1]{{#1}^\vee}

\newcommand{\mcd}[2]{\mathrm{gcd}(#1,#2)}

\newcommand{\GaZ}{{\Ga_{0}}}
\newcommand{\GaU}{{\Ga_{1}}}
\newcommand{\GaE}{{\Ga_{\vep}}}

\newcommand{\ideles}{{K_{\A}^{\times}}}

\newcommand{\inv}[1]{{#1}^{-1}}

\newcommand{\invol}[1]{{#1}^\dagger}

\newcommand{\iisom}{\buildrel\sim\over\longrightarrow}
\newcommand{\isom}{\buildrel\sim\over\rightarrow}

\newcommand{\jet}{{\rm jet}}

\newcommand{\map}{\rightarrow}
\newcommand{\mmap}{\longrightarrow}
\newcommand{\hmap}{\hookrightarrow}

\newcommand{\nr}[1]{{#1}^{\rm nr}}

\newcommand{\PH}{{\rm Pr}_{\infty}}
\newcommand{\PFr}{{\rm Pr}_{p}}

\def\smallmat#1#2#3#4{
    \left(\sopra{#1}{#3}\sopra{#2}{#4}\right)}
 
 \def\liminj#1{\mathrel{\mathop{\kern0pt\lim}
               \limits_{\kern3pt\longrightarrow\,#1}}}
                
 \def\limproj#1{\mathrel{\mathop{\kern0pt\lim}
                \limits_{\kern3pt\leftarrow\,#1}}}

\newcommand{\scal}[2]{\left\langle#1,#2\right\rangle}
\newcommand{\scalq}[2]{\left[#1,#2\right]}
\newcommand{\scaldR}[2]{{\scal{#1}{#2}}_{\rm dR}}

\newcommand{\sopra}[2]{\genfrac{}{}{0pt}{1}{#1}{#2}}

\newcommand{\vass}[1]{\left|#1\right|}
\newcommand{\vvass}[1]{\left|\!\vass{#1}\!\right|}

\newcommand{\vvec}[2]{\left(\sopra{#1}{#2}\right)}

\newcommand{\srh}{\sqrt{p_o}}

\def\timestamp{%
\scratch=\time
\divide\scratch by 60
\edef\hours{\the\scratch}
\multiply\scratch by 60
\minutes=\time
\advance\minutes by -\scratch
\the\day/\the\month $\,$---$\,$\hours:\null
\ifnum\minutes< 10 0\fi
\the\minutes}

\def\rmarginsert[#1]{\hglue 0pt\vadjust
{\null\vskip -\baselineskip\rightline{\abstractfont\rlap{\hfil\  #1}}}}

\def\today{\number\day\space\ifcase\month\or
January\or February\or March\or April\or May\or June\or
July\or August\or September\or October\or November\or December\fi,
\number\year}

\def\qedmark{\hbox{\vrule height 6pt width 5pt}}
\def\qedskip{\vrule height 4pt width 0pt depth 1pc}
\def\qed{\penalty 1000 $\,$\penalty 1000{\qedmark\qedskip}\par}

\def\pr#1{{#1^\prime}}

\newtheorem{dfn}{Definition}[section]
\newtheorem{thm}[dfn]{Theorem}
\newtheorem{rem}[dfn]{Remark}
\newtheorem{rems}[dfn]{Remarks}
\newtheorem{pro}[dfn]{Proposition}
\newtheorem{lem}[dfn]{Lemma}

\newtheorem{notat}[dfn]{Notation}

\newcommand{\pf}{\par\noindent{\bf Proof. }}

\begin{document}

\title{\LARGE Power Series Expansions of Modular Forms\\
       and Their Interpolation 
       Properties}

\author{\LARGE Andrea Mori\\
         \large Dipartimento di Matematica\\
         \large Universit\`a di Torino}

\date{}

\maketitle

\pagestyle{fancy}

\fancypagestyle{plain}{%
               \fancyhf{}
               \renewcommand{\headrulewidth}{0pt}}
\fancyhead{}
\fancyfoot{\normalsize\thepage}
\fancyhead[RE,LO]{}
\fancyhead[RO]{Power series expansions of modular forms}
\fancyhead[LE]{A. Mori}
\fancyfoot[RO,LE]{}
\fancyfoot[RE,LO]{}
\renewcommand{\headrulewidth}{0pt}


\begin{abstract}
We define a power series expansion of a holomorphic modular 
form $f$ in the $p$-adic neighborhood of a CM point $x$ of type $K$ 
for a split good prime $p$. The modularity group can be either a 
classical conguence group or a group of norm 1 elements in an order 
of an indefinite quaternion algebra. The expansion 
coefficients are shown to be closely 
related to the classical Maass operators and give $p$-adic 
information  on the ring of definition of $f$. By letting the CM 
point $x$ vary in its Galois orbit, the $r$-th coefficients define a 
$p$-adic $K^{\times}$-modular form in the sense of Hida. By coupling 
this form with the $p$-adic avatars of algebraic Hecke characters 
belonging to a suitable family and using a Rankin-Selberg type formula 
due to Harris and Kudla along with some explicit computations of Watson and of Prasanna,
 we obtain in the even weight case a 
$p$-adic interpolation for the square roots of 
a family of twisted special values of the automorphic 
$L$-function associated with the base change of $f$ to $K$.

\medskip\noindent
2000 Mathematics Subject Classification 11F67    
\end{abstract}


\section*{Introduction}
The idea that the power series expansion of a modular form at a CM 
point with respect to a well-chosen local parameter should have an 
arithmetic significance goes back to the author's  thesis, \cite{Mori94}. 
The goal of the thesis was to prove an expansion 
principle, namely a characterization of the ring of algebraic $p$-adic  
integers of definition of an elliptic modular form in terms of the 
coefficients of the expansion. Such a result would be analogous to the 
classical $q$-expansion principle  based on the Fourier expansion
(e.÷g. \cite{Katz73}), with the advantage of being  
generalizable in principle to groups of modularity without 
parabolic elements where Fourier series are not available. The 
simplest such situation is that of a Shimura curve attached to an 
indefinite non-split quaternion algebra $D$ over $\Q$ (quaternionic 
modular forms).

The basic idea in \cite{Mori94} was to consider a prime $p$ of good 
reduction for the modular curve that is split 
in the quadratic field of complex multiplications $K$ and use the 
Serre-Tate deformation parameter to construct a local parameter at the 
CM point $x$ corresponding to a fixed embedding of $K$ in the split quaternion algebra. 
The coefficients of the resulting power series are
related to the values obtained evaluating the 
$C^{\infty}$-modular forms $\de_{k}^{(r)}f$ at a lift $\tau$ 
of $x$ in the complex upper half-plane, where $k$ 
is the weight of $f$ and $\de_{k}^{(r)}$ is the $r$-th iterate, in the 
automorphic sense, of the basic Maass operator 
$\de_{k}=-\frac1{4\pi}\left(2i\frac{d}{dz}+\frac{k}{y}\right)$.
Our first goal in this paper is to prove a version of the expansion principle 
valid also for quaternionic modular forms
without making use of the local complex geometry and completely 
$p$-adic in nature. 
The realization of modular forms as global sections of a line bundle $\calL$
suitable for the Serre-Tate theory is subtler in the non-split case because
for Shimura curves the Kodaira-Spencer 
map $\KS\colon\Sym^2\oom\map\Om^1_{\calX}$ is not an isomorphism (for a trivial 
reason: the push-forward $\oom=\pi_{*}\Om^{1}_{\calA/\calX}$ for the 
universal family of ``false elliptic curves'' has rank 2). This motivates the introduction of 
$p$-ordinary test triples (definition  \ref{th:testpair}) that require moving to an 
auxiliary quadratic extension.
The abelian variety of dimension $\leq2$ 
corresponding to the CM point $x$ defined over the ring of $p$-adic algebraic 
integers $\calO_{(v)}$ is either a CM curve $E$ with $\End_{0}(E)=K$ or an 
abelian surface isogenous to a twofold product $E\times E$ of such a CM 
curve. To it we associate a complex period 
$\Om_{\infty}\in\C^{\times}$ and a $p$-adic period 
$\Om_p\in\calO_v^{{\rm nr},\times}$. If also the modular form is defined 
over $\calO_{(v)}$ and $\sum_{r=0}^{\infty}(b_{r}(x)/r!)T_{x}^{r}$ is 
its expansion obtained form the Serre-Tate theory, we  establish in 
theorem \ref{thm:equality} an equality
\begin{equation}
    c^{(r)}_{v}(x)=
    \de_{k}^{(r)}(f)(\tau)\Om_{\infty}^{-k-2r}=
    b_{r}(x)\Om_{p}^{-k-2r}
    \label{eq:intro1}
\end{equation}
of elements in $\calO_{(v)}$. The expansion 
principle, theorem \ref{thm:expanprinc}, asserts that if $f$ is a 
holomorphic modular form such that the numbers $c^{(r)}_{v}(x)$ 
defined by the complex side of the equality \eqref{eq:intro1} are in 
$\calO_{v}$ and the $p$-adic integers $\Om_p^{2r}c^{(r)}_{v}(x)$ satisfy 
the Kummer-Serre congruences, then $f$ is defined over the integral 
closure of $\calO_{(v)}$ in the compositum of all finite 
extensions of the quotient field of $\calO_{(v)}$ in which $v$ splits 
completely.

Suppose again that the holomorphic modular form $f$ is defined over a ring 
$\calO_{(v)}$ of $p$-adic integers. The numbers $c^{(r)}_{v}(x)$ are
related to the coefficients of a $p$-integral power series, i.e. to a 
$p$-adic measure on $\Z_{p}$, naturally attached to $f$. One may 
wonder about the interpolation properties of this measure. In the 
introduction of \cite{HaTi01} Harris and Tilouine suggest that in the 
case of an eigenform $f$ the author's techniques may be used in 
conjunction with the results of Waldspurger \cite{Waldsp85} to 
$p$-adically interpolate the square roots of the special values of the 
automorphic $L$-functions $L(\pi_{K}\otimes\xi,s)$, where $\pi_{K}$ is 
the base change to $K$ of the $\GL_{2}$-automorphic representation 
$\pi$ associated to $f$ (possibly up to Jacquet-Langlands correspondence) 
and $\xi$ belongs to a suitable family of Gr\"ossencharakters
for $K$. 

Our second goal for this paper is to partially fulfill this expectation 
when $f$ has even weight $2\ka$. 
A key observation (proposition \ref{th:meascrfx}) is that the set of 
values $c^{(r)}_{v}(x)$ for $x$ ranging in a full set of 
representatives of the copy of the generalized ideal class group 
$\ideles/K^{\times}\C^{\times}\wh{\calO}_{c}^\times$ embedded in the 
modular (or Shimura) curve extends to a Hida \cite{Hida86} $p$-adic
$\GL_{1}(K)$-modular form $\hat{c}_{r}$, which is essentially the 
$r$-th moment of a $p$-adic measure on $\Z_{p}$ with values in the 
unit ball of the $p$-adic Banach space of such $p$-adic forms. The 
scalar obtained by coupling the form $\hat{c}_{r}$ with the $p$-adic 
avatar of a Gr\"ossencharakter $\xi_r$ for $K$ trivial on $\wh\calO_{c}^\times$ 
and of suitable weight twisted by a power of the idelic norm is proportional to the integral
\begin{equation}
    J_{r}(f,\xi_r,\tau)=
    \int_{\ideles/K^\times\R^\times}\phi_r(td_\infty)\xi_r(t)\,dt
    \label{eq:intro2}
\end{equation}
where $\phi_r$ is the adelic lift of $\de_{2\ka}^{(r)}(f)$, 
$\tau\in\gotH$ represents $x$ and 
$d_{\infty}\in\SL_{2}(\R)$ is the standard parabolic matrix such 
that $d_{\infty}i=\tau$. When $\xi_r$ is of the form $\xi_r=\chi\xi^r$
and satisfies some technical conditions the value so obtained
is essentialy the $r$-th moment of a $p$-adic measure $\mu(f,x;\chi,\xi)$ on $\Z_p$.

On the other hand, the square of the integral \eqref{eq:intro2} 
is a special case of the generalized Fourier coefficients 
$L_{\underline{\xi}}(\Phi)$ studied by Harris and Kudla in 
\cite{HaKu91}. Building on results of Shimizu \cite{Shimi72} and 
refining the techniques of Waldspurger \cite{Waldsp85}, Harris and 
Kudla use the seesaw identity associated with the theta 
correspondence between the similitude groups $\GL_{2}$ and $\GO(D)$ 
and the splitting $D=K\oplus K^{\perp}$ to express the generalized 
Fourier coefficients $L_{\underline{\xi}}(\theta_{\varphi}(F))$ 
where $F\in\pi$ and $\varphi$ is a split primitive Schwartz-Bruhat 
function on $D_{\A}$ as a Rankin-Selberg Euler product. 
Thus, we can use the explicit version of Shimizu's theory worked out by Watson \cite{Wat03},
the local non-archimedean computations of Prasanna \cite{Pra06} together with
some local archimedean computations to obtain a formula relating the square of 
the $r$-th moment of $\mu(f,x;\chi,\xi)$ to the values $L(\pi_K\otimes\chi\xi^r,\frac12)$
whose local correcting terms are explicit outside the primes dividing the conductor
of the Gr\"ossencharakter and the primes dividing the non square-free part of the level
(theorem \ref{thm:maininterpolation}).

Some natural questions arise. 
First of all, one would like to compute the special 
values of the $p$-adic $L$-function attached to the measure $\mu(f,x;\chi,\xi)$. 
Secondly, one may ask if the methods can be extended to treat different or
more general families of
Gr\"ossencharakters, in particular if one can control the interpolation as 
the ramification at $p$ increases. Proposition 
\ref{th:oldforms} implies that, if anything, this cannot be achieved 
without moving the CM point. Thus, some kind of geometric construction 
in the modular curve may be in order, with a possible link to the 
question of the determination of the action of the Hecke operators on 
the Serre-Tate expansions. Another question is whether the 
reinterpretation of the integral \eqref{eq:intro2} as inner product in 
the space of $p$-adic $\GL_{1}(K)$-modular forms can be used to 
obtain an estimate of the number of non-vanishing special values 
$L(\pi_{K}\otimes\xi,\frac12)$.
We hope to be able to attack these problems in a future paper.

\paragraph{Acknowledgements.} The idea that the power series 
coefficients may be used to $p$-adically interpolate the special 
values $L(\pi_{K}\otimes\xi,\frac12)$ arose a long time ago in 
conversations with Michael Harris. I wish to thank Michael Harris for 
sharing his intuitions and for many useful suggestions.

Also, I wish to thank the anonymous referee of a previous version of 
the manuscript, whose suggestions helped greatly to remove some 
unnecessary hypotheses.

\paragraph{Notations and Conventions.}
The symbols $\Z$, $\Q$, $\R$, $\C$ and $\F_q$ denote, as usual, the integer, 
the rational, the real, the complex numbers and the field with $q$ elements respectively. 
We fix once for all an embedding $\imath\colon\overline{\Q}\map\C$ and 
by a number field we mean a finite subextension of the field 
$\overline{\Q}$ of algebraic numbers.
If $L$ is a number field, we denote $\calO_{L}$ its ring of integers and 
$\de_{L}$ its discriminant. If $L=\Q(\sqrt{d})$ is a quadratic field, for each 
positive integer $c$ we denote 
$\calO_{L,c}=\Z+c\calO_{L}=\Z[c\om_d]$ its order of conductor $c$, 
with $\om_d=\sqrt d$ if $d\equiv 2$, $3\bmod 4$ or $\om_d=(1+\sqrt d)/2$ if $d\equiv 1\bmod 4$.
If $[L:\Q]=n$ we denote $I_L=\{\si_1,\dots,\si_n\}$ the set of embeddings $\si_i:L\map\C$ and we assume $\si_1=\imath_{|L}$.

If $p$ is a rational prime we denote $\Z_{p}$ and $\Q_{p}$ the $p$-adic 
integers and the $p$-adic numbers respectively. By analogy, $\Q_\infty=\R$.
If $v|p$ is a  place of  the number field 
$L$ corresponding to the maximal ideal $\gotp_{v}\subset\calO_{L}$, we denote 
$\calO_{(v)}$, $L_{v}$, $\calO_{v}$, $k(v)$ the localization of $\calO_{L}$ at 
$\gotp_{v}$, the $v$-adic completion of $L$, the ring of 
$v$-adic integers in $L_{v}$ and the residue field respectively. 
The maximal ideal in $\calO_{v}$ is 
still denoted $\gotp_{v}$. Also, we denote $\nr{L}_v$ the maximal 
unramified extension of $L_{v}$ and $\nr{\calO}_v$ its ring of 
integers. 

We denote $\widehat{\Z}$ the profinite completion of $\Z$ and for 
each $\Z$-module $M$ we let $\widehat M=M\otimes\widehat\Z$. We 
denote $\A$ the ring of rational adeles and $\A_{f}$ the finite 
adeles, so that $\A=\R\times\A_{f}=\Q\R\widehat\Z$. For a number field $L$ we denote 
$\A_L=\A\otimes L$  and $L_\A^\times$ the corresponding ring of adeles and group of  \`{\i}deles respectively. If $\gotn\subseteq\calO_{K}$ is an ideal, we let 
$L^{\times}_{\gotn}=\{\la\in L^{\times}\mbox{ such that }\la\equiv1\bmod\gotn\}$
and denote $\calI_{\gotn}$ the group of fractional ideals of  $L$ prime with $\gotn$, 
$P_{\gotn}$ the subgroup of principal fractional ideals generated by the elements in
$L^{\times}_{\gotn}$ and $U_{\gotn}$ the subgroup of finite \`{\i}deles product of local units congruent to $1${} $\bmod\gotn$.

We fix an additive character $\psi$ of $\A/\Q$, by asking that $\psi_\infty(x)=e^{2\pi i x}$
and $\psi_p$ is trivial on $\Z_p$ with $\psi_p(x)=e^{2\pi i x}$ for $x\in\Z[\inv p]$ and finite $p$. On $\A$ we fix the Haar measure $dx=\prod_{p\leq\infty}dx_p$ where the local Haar measures $dx_p$ are normalized so that the $\psi_p$-Fourier transform is autodual. For a quaternion algebra $D$ with reduced norm $\nu$, we fix on $D_\A$ the Haar measure $dx=\prod_{p\leq\infty}dx_p$ where the local Haar measures $dx_p$ are normalized so that the Fourier transform with respect to the norm form is autodual. 
Let $(V,\scal{\,}{\,})$ be a quadratic space of dimension $d$ over $\Q$.
We denote $\calS_{\A}(V)=\bigotimes_{p\leq\infty}\calS_{p}$ the adelic 
Schwartz-Bruhat space, where for $p$ finite, $\calS_{p}$ is the space of 
Bruhat functions on $V\otimes\Q_p$ and $\calS_{\infty}$ is the space of Schwartz functions on 
$V\otimes\R$ which are finite under the natural action of a (fixed) maximal compact subgroup of the similitude group $\GO(V)$. The Weil representation $r_\psi$ is the representation of $\SL_2(\A)$ on 
$\calS_{\A}(V)$ which is explicitely described locally at $p\leq\infty$ by
\begin{subequations}\label{eq:Weilrep}
\begin{eqnarray}
r_\psi\left(\begin{array}{cc}1 & b \\0 & 1\end{array}\right)\varphi(x) & = & \psi_p\left(\frac12\scal{bx}{x}\right)\varphi(x),
\label{eq:Weilrep1} \\
r_\psi\left(\begin{array}{cc}a & 0 \\ 0 & \inv a\end{array}\right)\varphi(x) & = & \chi_V(a)\vass{a}_p^{d/2}\varphi(ax)
\label{eq:Weilrep2} \\
r_\psi\left(\begin{array}{cc}0 & 1 \\ -1 & 0\end{array}\right)\varphi(x) & = & \gamma_V\hat{\varphi}(x)
\label{eq:Weilrep3}
\end{eqnarray}
\end{subequations}
where $\gamma_V$ is an eighth root of 1 and $\chi_V$ is a quadratic character that are computed in our cases of interest in \cite{JaLa70} (see also the table in \cite[\S3.4]{Pra06}), while the Fourier transform $\hat{\varphi}(x)=\int_{V\otimes\Q_p}\varphi(y)\psi_p(\scal xy)\,dy$ is computed with respect to a $\scal{\,}{\,}$-self dual Haar measure on $V\otimes\Q_p$.

If $R$ is a ring and $M$ a $R$-module we denote $\dual M={\rm Hom}(M,R)$ 
the dual of $M$. The same notation applies to a sheaf of modules over a scheme. 
If $G$ is a subgroup of units in $R$ we say that non-zero elements 
$x$, $y\in M$ are $G$-equivalent and 
write $x\sim_{G}y$ if there exists $r\in G$ such that $rx=y$. 

The group $\SL_{2}(\R)$ acts on the complex upper half-plane $\gotH$ 
by linear fractional transformations, if $g=\smallmat abcd$ then 
$g\cdot z=\frac{az+b}{cz+d}$. The automorphy factor is defined to be 
$j(g,z)=cz+d$. The action extends to an action of the 
group $\GL^{+}_{2}(\R)$. 
If $\Ga<\SL_{2}(\R)$ is a Fuchsian group of the first kind we shall denote 
$M_{k}(\Ga)$ the space of modular forms of weight
$k\in\Z$ with respect to $\Ga$ i.e. the holomorphic functions $f$ 
on $\gotH$ such that
$$
\mbox{$f(\ga z)=f(z)j(\ga,z)^k$ for all $z\in\gotH$ and
$\ga\in\Ga$}
$$
and extend holomorphically to a neighborhood of each cusp (when cusps 
exist). The subspace of cuspforms, i.e. those modular forms that 
vanish at the cusps, will be denoted $S_{k}(\Ga)$.
The request that a holomorphic function on $\gotH$ extends 
holomorphically to a neighborhood of a cusp $s$ is equivalent to a certain 
growth condition as $z\to s$. Relaxing holomorphicity 
but mantaining the growth condition yields the much bigger spaces of 
$C^\infty$-modular and {}-cuspforms, which will be denoted 
$M_{k}^{\infty}(\Ga)$ and $S_{k}^{\infty}(\Ga)$
respectively. We will denote
$$
M_{k,\vep}(\De,N),\quad
S_{k,\vep}(\De,N),\quad
M_{k,\vep}^{\infty}(\De,N),\quad
S_{k,\vep}^{\infty}(\De,N)
$$
the above spaces of modular or cuspforms with respect to the groups 
$\Ga=\GaE(\De,N)$, $\vep\in\{0,1\}$,  defined in section \ref{se:curves}.
It is a well-known fact that $M_{k,\vep}(\De,N)$ is always 
finite-dimensional and trivial for $k<0$.


\section{Modular and Shimura curves}

\subsection{Quaternion algebras.}\label{se:quatalg}
Let $D$ be a quaternion algebra over $\Q$ with reduced norm $\nu$ 
and reduced trace $\tra$. For each place $\ell$ of $\Q$ let $D_\ell=D\otimes_Q\Q_\ell$.
Let $\Si_D$ be the set of places at which $D$ is \emph{ramified}, i.~e. 
$D_\ell$ is the unique, up to isomorphism, quaternion division 
algebra over $\Q_\ell$. If $\ell\notin\Si_D$ the algebra $D$ is 
\emph{split} at $\ell$, i.~e. $D_\ell\simeq\M_2(\Q_\ell)$. The set $\Si_D$ 
is finite and even and determines completely the isomorphism class 
of $D$. Moreover, every finite and even subset of places of $\Q$ is the set 
of ramified places of some quaternion algebra over $\Q$
(for these and the other basic results on quaternion algebras the standard 
reference is \cite{Vigner80}).
In particular, $M_2(\Q)$ is the only quaternion algebra up to 
isomorphism which is \emph{split}, i.~e. split at all 
places. The discriminant $\De=\De_D$ of $D$ is the product of 
the finite primes in $\Si_D$ if $\Si_D\neq\emptyset$, or $\De=1$ 
otherwise.
We shall henceforth assume that $D$ is \emph{indefinite}, 
i.~e. split at $\infty$, and fix an isomorphism 
$\Phi_\infty\colon D_\infty\isom\M_2(\R)$
which will be often left implicit.
There is a unique conjugacy class of maximal orders in $D$. Once for 
all, choose a maximal order $\calR_{1}$ and fix isomorphisms
$\Phi_\ell\colon D_\ell\isom\M_2(\Q_\ell)$ for
$\ell\notin\Si_D$ so that 
$\Phi_\ell(\calR_{1})=\M_2(\Z_\ell)$. For an integer $N$
prime to $\De$ let $\calR_{N}$ be the level $N$ 
Eichler order of $D$ such that
$$
\calR_{N}\otimes_\Z\Z_\ell=\inv{\Phi_\ell}
\left(\left\{\left(
\begin{array}{cc}
	a & b  \\
	c & d
\end{array}\right)
\hbox{$\in\M_2(\Z_\ell)$ such that $c\equiv0\bmod N$}
\right\}\right)
$$
for $\ell\notin\Si_D$, and $\calR_{N}\otimes\Z_\ell$ is the unique
maximal order in $D_\ell$ for $\ell\in\Si_D$.
If $D=\M_2(\Q)$ we take $\calR_{1}=\M_2(\Z)$ and
$\calR_{N}=\left\{\smallmat abcd
\hbox{$\in\M_2(\Z)$ such that $c\equiv0\bmod N$}\right\}$.

There are exactly two homomorphisms
$\orn_\ell^{1},\orn_\ell^{2}\colon\calR_{N}\otimes\F_\ell\mmap\F_{\ell^2}$ 
for each prime $\ell|\De$, and two homomorphisms
$\orn_\ell^{1},\orn_\ell^{2}\colon\calR_{N}\otimes\F_\ell\mmap{\F_{\ell}}^2$ 
for each prime $\ell|N$. These maps are 
called $\ell$-\emph{orientations} and the two $\ell$-orientations are 
switched by the non-trivial automorphism of either $\F_{\ell^2}$ or 
${\F_{\ell}}^{2}$. An orientation for $\calR_{N}$ is the choice of an 
$\ell$-orientation $\orn_\ell$ for all primes $\ell|N\De$.

An involution $d\mapsto\invol d$ in $D$ is \emph{positive} if 
$\tra(d\invol d)>0$ for all $d\in D$. By the Skolem-Noether theorem 
\begin{equation}
	\invol d=\inv t\bar{d}t
	\label{eq:involution}
\end{equation}
where $t\in D$ is some element such that $t^2\in\Q^{{}<0}$ and 
$d\mapsto\bar{d}$ denotes quaternionic conjugation, 
$d+\bar{d}=\tra(d)$. If $t\in D$ is 
such an element, let $B_t$ be the bilinear form on $D$ defined by
\begin{equation}
	B_t(a,b)=\tra(a\bar{b}t)=\tra(at\invol b)\qquad
	\hbox{for all $a, b\in D$}.
	\label{eq:formEt}
\end{equation}
If $\calR\subset D$ is an order, the involution $d\mapsto\invol d$ 
is called $\calR$-\emph{principal} if $\invol\calR=\calR$ and the bilinear form 
$B_t$ is skew-symmetric, non-degenerate and $\Z$-valued on 
$\calR\times\calR$ with pfaffian equal to $1$. 
When $\De>1$ an explicit model for the triple 
$(D,\calR_{N},d\mapsto\invol d)$ can be constructed as follows.
The condition $(n,-N\De)_{\ell}=-1$ for all $\ell\in\Si_D$
on Hilbert symbols defines for $n$  a certain subset 
of non-zero congruence classes modulo $N\De$. 
Passing to classes modulo $8N\De$ 
and taking $n>0$ we may assume that 
$(n,-N\De)_{\infty}=(n,-N\De)_{p}=1$ for all primes $p$ 
dividing $N$ and also $(n,-N\De)_{2}=1$ if $\De$ is odd.
By Dirichlet's theorem of primes in arithmetic progressions there 
exists a prime $p_{o}$ satisfying these conditions and the product 
formula easily implies that
\begin{equation}
    (p_{0},-N\De)_{\ell}=-1\qquad\mbox{if and only if $\ell\in\Si_{D}$.}
    \label{eq:condonHS}
\end{equation}
Let $a\in\Z$ such that $a^2N\De\equiv-1\bmod p_o$.

\begin{thm}[Hashimoto, \cite{Hashim95}]\label{th:Hashimoto}
	Let $D$ be a quaternion algebra over $\Q$ of discriminant $\De$ and let
	$t\in D$ such that $t^2\in\Q^{{}<0}$. Then:
	\begin{enumerate}
		\item  $D$ is isomorphic to the quaternion algebra 
		       $D_H=\Q\oplus\Q i\oplus\Q j\oplus\Q ij$, where 
		       $i^2=-N\De$, $j^2=p_o$ and $ij=-ji$;
		\item  the order $\calR_{H,N}=\Z\ep_1\oplus\Z\ep_2\oplus\Z\ep_3
		       \oplus\Z\ep_4$, where $\ep_1=1$, $\ep_2=(1+j)/2$, 
		       $\ep_3=(i+ij)/2$ and $\ep_4=(aN\De j+ij)/p_o$ is an Eichler 
		       order of level $N$ in $D_H$;
		\item  the skew symmetric form $B_t$ on $D_H$ is $\Z$-valued on 
		       $\calR_{H,N}$ if and only if $ti\in\calR_{H,N}$. 
		       Moreover, it defines a non-degenerate 
		       pairing on $\calR_{H,N}\times\calR_{H,N}$ if and only if 
		       $ti\in\calR_{H,N}^\times$;
		\item  let $t=\inv i$. Then the elements 
		       $\eta_1=\ep_3-\frac12(p_o-1)\ep_4$, $\eta_2=-aN\De-\ep_4$, 
		       $\eta_3=1$ and $\eta_4=\ep_2$ are a symplectic 
		       $\Z$-basis of $\calR_{H,N}$.
	\end{enumerate}
\end{thm}

\noindent We call \emph{Hashimoto model} of a quaternion algebra endowed with 
an Eichler order $\calR$ of level $N$ and a $\calR$-principal positive involution 
the triple $(D_H,\calR_{H,N},\inv i)$ given in the above theorem.
We can fix the isomorphism $\Phi_\infty$ for the Hashimoto model by declaring 
that
$$
\Phi_\infty(i)=
\left(
\begin{array}{cc}
	0 & -1  \\
	N\De & 0
\end{array}
\right),\qquad
\Phi_\infty(j)=
\left(
\begin{array}{cc}
	\sqrt{p_o} & 0  \\
	0 & -\sqrt{p_o}
\end{array}
\right).
$$

\subsection{Moduli spaces.}\label{se:curves}
Fix a $\calR_{1}$-principal positive involution 
$d\mapsto\invol d$ as in \eqref{eq:involution}.
We shall consider the groups
$$
\GaZ(\De,N)=\calR^{1}_{N}=\left\{
\hbox{$\ga\in\calR_{N}$ such that $\nu(\ga)=1$}
\right\}
$$
and
$$
\GaU(\De,N)=\left\{\hbox{$\ga\in\GaZ(\De,N)$ such that 
$\orn^{\ep}_{\ell}(\ga r)=\orn^{\ep}_{\ell}(r)$
 for all $r\in\calR_{N}$, $\ell|N, \ep=1,2$}\right\}.
$$
When $\De=1$, $\GaZ(1,N)$ and $\GaU(1,N)$ are the classical 
congruence subgroup
$$
\Ga_0(N)=\left\{
\left(
\begin{array}{cc}
	a & b  \\
	c & d
\end{array}
\right)
\in\SL_2(\Z)\hbox{ such that $c\equiv0\bmod N$}\right\}
$$
and
$$
\Ga_1(N)=
\left\{
\left(
\begin{array}{cc}
	a & b  \\
	c & d
\end{array}
\right)
\in\SL_2(\Z)\hbox{ such that $a,d\equiv1$ e $c\equiv0\bmod N$}\right\}
$$
respectively. Since $D$ is indefinite $\GaE(\De,N)$ for 
$\vep\in\{0,1\}$ is, via $\Phi_\infty$, a discrete subgroup of 
$\SL_2(\R)$ acting on the complex upper half plane $\gotH$.
When $\De>1$ the quotient $X_{\vep}(\De,N)=\GaE(\De,N)\bs\gotH$ is a 
compact Riemann surface, \cite[proposition~9.2]{ShiRed}. When $\De=1$
let $X_\vep(N)$ be the standard cuspidal 
compactification of $Y_\vep(N)=\GaE(N)\bs\gotH$.

Each of these complete curves $X$ has a canonical model over $\Q$, 
\cite{ShiRed}.
In fact, each $X$ can be reinterpreted as the set of complex 
points of a scheme $\calX$ which is the solution of a moduli 
problem, defined over $\Z[1/{N\De}]$, 
e.g. \cite{BerDar96, DelRap73, DiaIm95, Milne79, Robert89}. 
When $D=M_2(\Q)$ and $N>3$, the functor 
$F_1(N)\colon\hbox{\bf $\Z[1/N]$-Schemes}\map\hbox{\bf Sets}$ defined 
by 
$$
F_1(N)(S)=
\left\{
\sopra{\hbox{Isomorphism classes of generalized elliptic curves $E=E_{|S}$}}
{\hbox{with a section $P\colon S\map E$ of exact order $N$}},
\right\}
$$
is represented by a proper and smooth $\Z[\frac1{N}]$-scheme 
$\calX_1(N)$ such that $\calX_1(N)(\C)=X_{1}(N)$.
The complex elliptic curve with point $P$ of exact order $N$ 
corresponding to $z\in\gotH$ is the torus $E_z=\C/\Z\oplus\Z z$ 
with $P=1/N \bmod\Z$. Denote
\begin{equation}
	\pi_N\colon\calE_{N}\mmap\calX_{1}(N)
	\label{eq:univEC}
\end{equation}
the universal generalized elliptic curve attached to the representable functor 
$F_{1}(N)$.  The scheme $\calX_0(N)$ quotient of $\calX_1(N)$ by 
the action of the group of diamond operators 
$\langle a\rangle\colon\calX_1(N)\map\calX_1(N)$, 
$\langle a\rangle(E,P)=(E,aP)$ for all $a\in(\Z/N\Z)^\times$, 
is the coarse moduli scheme attached to the functor
$$
F_0(N)(S)=
\left\{
\sopra{\hbox{Isomorphism classes of generalized elliptic curves $E=E_{|S}$}}
{\hbox{with a cyclic subgroup $C\subset E$ of exact order $N$}}
\right\}
$$
and a smooth $\Z[1/N]$-model for the curve $X_0(N)$.

When $\De>1$ and $N>3$, $X_1(\De,N)=\calX_{1}(\De,N)(\C)$ for
the proper and smooth $\Z[1/{N\De}]$-scheme $\calX_{1}(\De,N)$ 
representing the functor
$F_1(\De,N)\colon\hbox{\bf $\Z[1/{N\De}]$-Schemes}\map\hbox{\bf Sets}$ 
defined by 
$$
F_1(\De,N)(S)=
\left\{
\sopra{
   \sopra{\hbox{Isomorphism classes of compatibly principally polarized}}
         {\hbox{ abelian surfaces $A=A_{|S}$ with a ring embedding}}}
   {\sopra{\hbox{$\calR_{1}\hmap\End(A)$ and an equivalence class of}}
         {\hbox{$\calR_N$-orientation preserving level $N$ structures}}} 
\right\}.
$$
A level $N$ structure on an abelian surface $A$ with 
$\calR_{1}\subset\End(A)$ is an isomorphism of (left) 
$\calR_{1}$-modules $A[N]\simeq\calR_{1}\otimes(\Z/N\Z)$. Two 
such structures are declared equivalent if they coincide on 
$\calR_N\otimes(\Z/N\Z)$ and induce the same 
$\ell$-orientations on $\calR_N$ for all $\ell|N$. The principal polarization is 
compatible with the embedding $\calR_{1}\subset\End(A)$ if the 
involution $d\mapsto\invol d$ is the Rosati involution. The abelian 
surfaces in $F_1(\De,N)(S)$ are called 
\emph{abelian surfaces with quaternionic multiplications} (QM-abelian surfaces, for short) or 
\emph{false elliptice curves}.
The complex QM-abelian surface corresponding to $z\in\gotH$ is 
\begin{equation}
	A_z=D_\infty^z/\calR_{1},
	\label{eq:QMtori}
\end{equation}
where $D_\infty^z$ is the real vector space $D_\infty$ endowed with the 
$\C$-structure defined by the identification 
$\C^2=\Phi_\infty(D_\infty)\left(\sopra{z}{1}\right)$, i.e. 
$A_z=\C^2/\Phi_\infty(\calR_{1})\vvec{z}{1}$. The complex 
uniformization \eqref{eq:QMtori} defines a level structure 
$\inv{N}\calR_{1}/\calR_{1}=(D/\calR_{1})[N]\isom(A_{z})[N]$ 
and the skew-symmetric form
$\scal{\Phi_\infty(a)\left(\sopra{z}{1}\right)}
{\Phi_\infty(b)\left(\sopra{z}{1}\right)}=B_t(a,b)$ 
for all $a,b\in D$,
where $B_t$ is as in \eqref{eq:formEt}, extended to $\C^{2}$ by 
$\R$-linearity is the unique Riemann form on $A_{z}$ with Rosati 
involution $d\mapsto\invol d$, \cite[lemma~1.~1]{Milne79}. Denote
\begin{equation}
	\pi_{\De,N}\colon\calA_{\De,N}\mmap\calX_{1}(\De,N)
	\label{eq:univQMAV}
\end{equation}
the universal QM abelian surface attached to the representable functor 
$F_{1}(\De,N)$.

As with the split case, a smooth $\Z[1/{N\De}]$-model 
$\calX_{0}(\De,N)$ of $X_{1}(\De,N)$ can be obtained as quotient of 
$\calX_{1}(\De,N)$ by a suitable action of 
$\frac{\GaZ(\De,N)}{\GaU(\De,N)}\simeq(\Z/N\Z)^{\times}$. It is the 
coarse moduli space for the functor
$$
F_0(\De,N)(S)=
\left\{
\sopra{\sopra{\hbox{Isomorphism classes of compatibly principally polarized}}
{\hbox{abelian surfaces $A=A_{|S}$ with a ring embedding $\calR_{1}\hmap\End(A)$}}}
{\hbox{and an $\calR_{N}$-equivalence class 
of level $N$ structures}}
\right\}
$$
where two level $N$ structures are $\calR_{N}$-equivalent if they coincide on 
$\calR_{N}\otimes(\Z/N\Z)$. 

\begin{rem}
    \rm In order to study the reduction of the modular and Shimura 
    curves at primes dividing $N\De$ one has to extend the moduli 
    problems described above to moduli problems defined over $\Z$, 
    see \cite{BoCa91, KatMaz85}. 
    The $\Z$-schemes thus obtained are proper but not smooth. We 
    shall not deal with primes of bad reduction and for the purposes 
    of this paper the above descriptions will suffice. 
\end{rem}

\subsection{Subfields and CM points.}\label{se:CMpts}
Let $\Q\subseteq\pr L\subset L$ be a tower of fields with $[L:\pr L]=2$ 
and assume that $L$ splits $D$, i.~e.  $D\otimes_\Q L\simeq\M_2(L)$ or, 
equivalently, that $L$ admits an embedding in $D\otimes_\Q\pr L$. 
An embedding $\jmath:L\hmap D\otimes_\Q\pr L$ endows $D\otimes_\Q\pr L$ 
with a structure of $L$-vector space. Scalar multiplication by 
$\la\in L$ is left multiplication by $\jmath(\la)$. 
The opposite algebra $\Dop$ acts $L$-linearly on $D$ by right multiplication, 
providing a direct identification 
\begin{equation}
\label{eq:DasEnd}
\Dop\otimes L\stackrel{\sim}{\mmap}\End_{L}(D\otimes\pr L).
\end{equation}
Let $\si$ be the non-trivial element in $\Gal(L/\pr L)$ and $\jmath^\si(\la)=\jmath(\la^\si)$
for all $\la\in L$. By the Skolem-Noether theorem 
there exists $u\in(D\otimes_\Q\pr L)^\times$, well defined up to a 
$L^\times$-multiple, such that $u\jmath(\la)=\jmath^\si(\la)u$ 
for all $\la\in L$ and $u^2\in\pr L$. Thus, with a slight abuse of notation,
the embedding $\jmath$ defines a splitting
\begin{equation}
	D\otimes\pr L=L\oplus L u
	\label{eq:Dsplit}
\end{equation}
which can be more intrinsically seen as the eigenspace decomposition 
under right multiplication by $\jmath(L^\times)$.
Also, there is an isomorphism
\begin{equation}
	D\iisom\Dop,\qquad
	\la_1+\la_2u\mapsto\la_1+\la_2^{\si}u.
	\label{eq:isoDDop}
\end{equation}
Let $L=\pr L(\al)$ with $\al^2=A\in\pr L$. The element
\begin{equation}
\label{eq:idempotent}
e_\jmath=\frac{1}{2}\left(1\otimes1+
\frac{1}{A}\jmath(\al)\otimes\al\right)\in D\otimes\pr L
\end{equation}
is an idempotent which is easily seen to be, under 
\eqref{eq:DasEnd}, \eqref{eq:Dsplit} and \eqref{eq:isoDDop}, the projection onto $L$ 
with kernel $Lu$. If $L\subseteq\C$ the idempotent 
$e_\jmath$ defines a projector in $D_\infty^z$ for all $z\in\gotH$ 
by scalar extension.

An involution $d\mapsto\invol d$ in $D$ extends by linearity to 
$D\otimes\pr L$. If $\invol\jmath$ is the embedding 
$\invol\jmath(\la)=\jmath(\la)\invol{}$, the explicit description 
\eqref{eq:idempotent} implies at once that 
$e_{\invol\jmath}=\invol{e_\jmath}$ and in particular
$$
\hbox{$\invol{e_\jmath}=e_\jmath$ if and only if 
$\invol{\jmath(L)}=\jmath(L)$ pointwise.}
$$
When the involution is positive a fixed idempotent can be constructed as follows. 
As an element of $\End(D)$ the involution \eqref{eq:involution} 
has determinant $-1$. Since $\invol1=1$ and $\tra(\invol d)=\tra(d)$ 
for all $d\in D$ its $(-1)$-eigenspace is a subspace of trace 
$0$ elements of dimension either $1$ or $3$. 
If the dimension is $3$ then the involution is the quaternionic conjugation, 
contradicting the positivity assumption. 
Therefore there exist a non-zero element $d$ of trace $0$ fixed by the involution.
The  subalgebra $F=\Q(d)\subset D$ 
is a quadratic field fixed by the involution and the corresponding 
idempotent $e\in D\otimes_\Q F$ has the desired property. 
Note that the positivity of the involution implies further 
that $F$ is real quadratic.

The conductor of an embedding $\jmath\colon L\map D$ of the quadratic field $L$
relative to the order $\calR_N$ is the integer $c=c_N>0$ such that 
$\jmath(\calO_{L,c})=\jmath(L)\cap\calR_{N}$.
Denote $\bar c$ the \textit{minimal conductor}, 
i.e. the conductor relative to the maximal order 
$\calR_{1}$. It is clear that $c$ is a multiple of $\bar c$, in fact 
$c/\bar c$ is a divisor of $N$ because $\calO_{L,c}/\calO_{L,\bar c}$ injects 
into $\calR_{1}/\calR_{N}\simeq\Z/N\Z$. In the following result the embedding is left implicit to simplify the notation.

\begin{pro}\label{prop:decomporder}
   Let $L\subset D$ be a quadratic subfield with associated decomposition $D=L\oplus Lu$. 
   Let $\La=L\cap\calR_{N}$ and $\La^\prime=Lu\cap\calR_{N}$. Then:
   \begin{enumerate}
  \item $D$ is split at the prime $p$ if and only if $(u^2,\de_L)_p=1$;
  \item if $p$ is unramified in $L$ and $\mcd pc=1$ then 
            $\calR_N\otimes\Z_p=\La\otimes\Z_p\oplus\La^\prime\otimes\Z_p$. Moreover, 
            $\La^\prime\otimes\Z_p=\calJ u$ for some fractional ideal $\calJ\subset L\otimes\Q_p$ 
            such that ${\rm N}(\calJ)\nu(u)=(q^\ep)$ with $\ep=1$ if $q|N\De$ and $\ep=0$ otherwise. 
\end{enumerate}
\end{pro}

\pf Let $L=\Q(\sqrt d)$. Then $\{1,\sqrt d,u, \sqrt du\}$ is a $\Q$-basis of $D$ and the local invariants of the norm form are $\det=1$ and $\ep_p=(-1,-1)_p(u^2,d)_p=(-1,-1)_p(u^2,\de_L)_p$, thus proving the first part.

For the second part, choose $u$ so that $u^2\in\Z$. Then there is an inclusion of orders 
$\calR^\prime=\calO_{L,c}\oplus\calO_{L,c}u\subseteq\La\oplus\La^\prime\subseteq\calR_N$. 
The elements $\{1,c\om_d,u,c\om_du\}$ are a $\Z$-basis of $\calR^\prime$, so that $\calR^\prime$ has reduced discriminant $\de_Lcu^2$. We are thus reduced to check that when $p|u^2$ and 
$\mcd p{c\de_L}=1$ then there is no element $x\in\calR_N$ of the form $x=(r+r^\prime u)/p$ with $r$, 
$r^\prime\in\calO_{L,c}-p\calO_{L,c}$. For such an element $x$ one must have $p|\tra(r)$ and $p|{\rm N}(r)$ from which one derives quickly a contradiction. 

The last claim follows from the very same discriminant computation since $\calR_N$ has reduced discriminant $N\De$.\qed

Fix a quadratic imaginary field $K$ that splits $D$. 
Exactly one of the two embeddings $\jmath,\jmath^{\si}$ is normalized 
in the sense of \cite[(4.~4.~5)]{ShiRed}.
The normalized embeddings correspond bijectively to a special subset 
of points $\tau\in\gotH$. More precisely, there is a bijection
$$
	\left\{
	\sopra{\displaystyle\hbox{normalized embeddings}}
	{\displaystyle\jmath\colon K\hmap D}
	\right\}
	\longleftrightarrow
	\CM_{\De,K}=\left\{
	\sopra{\displaystyle\hbox{$\tau\in\gotH$ such that 
	        $\Phi_\infty(\jmath(K^\times))=$}}
	{\displaystyle\{\ga\in\Phi_\infty(D^\times)
	\cap\GL_2^+(\R)~|~\ga\cdot\tau=\tau\}}
	\right\}.
$$
The bijection is $\GaZ(\De,N)$-equivariant where $\GaZ(\De,N)$ acts by 
conjugation on the left set and on $\CM_{\De,K}$ via its action on 
$\gotH$.
Also, the correspondence $\jmath\leftrightarrow\tau$ is characterized 
by the fact that the complex structure on $D_\infty$ induced by the 
embedding $\jmath$ coincides with that of $D_{\infty}^{\tau}$.
In the split case $\CM_{1,K}=K\cap\gotH$. 

We shall denote $c_{\tau}=c_{\tau,N}$ the conductor relative to the order 
$\calR_{N}$ of the embedding associated to the point $\tau\in\CM_{\De,K}$ 
and $\bar c_{\tau}$ its minimal conductor.

\begin{pro}
   Let $\tau$ and $\pr\tau\in\CM_{\De,K}$ such that $\pr\tau=\ga\cdot\tau$ 
   for some $\ga\in\GaZ(\De,N)$. 
   Then $c_{\pr\tau,N}=c_{\tau,N}$.
\end{pro}

\pf Let $\jmath$ and $\pr\jmath$ be the embeddings corresponding 
to $\tau$ and $\pr\tau$ respectively. 
Then $\pr\jmath=\ga\jmath\inv{\ga}$ and so
$\pr\jmath(\calO_{c_{\pr\tau,N}})=\pr\jmath(K)\cap\calR_{N}=
\ga\jmath(K)\inv{\ga}\cap\calR_{N}=
\ga(\jmath(K)\cap\calR_{N})\inv\ga=\ga\jmath(\calO_{c_{\tau,N}})\inv\ga=
\pr\jmath(\calO_{c_{\tau,N}})$. 
\qed

\begin{dfn}
	A point $x\in X_{0}(\De,N)$ is a CM point of type $K$ and conductor 
	$c=c_x$ if it is represented by a $\tau\in\CM_{\De,K}$ with 
	$c_{\tau,N}=c$. Denote
	$$
         \CM(\De,N;\calO_{K,c})=\{\mbox{CM points of $X_{0}(\De,N)$ of type $K$ and 
         conductor $c$}\}.
          $$
\end{dfn}

The following result is \cite[Lemma 4.17]{Dar04}

\begin{pro}\label{teo:existCM}
    Let $c>0$ be an integer such that $\mcd c{N\De}=1$. Then the set
    $\CM(\De,N;\calO_{K,c})$ is non-empty if and only if
    \begin{itemize}
        \item  all primes $\ell|\De$ are inert in $K$, and
        \item  all primes $\ell|N$ are split in $K$.
    \end{itemize}
\end{pro}

For $\tau\in\CM_{1,K}$ the elliptic curve $E_{\tau}$ has complex 
multiplications in the field $K$. When $\De>1$ and $\tau\in\CM_{\De,K}$
the QM abelian surface $A=A_{\tau}=\calA(\C)$ contains the elliptic curve
$E=K\otimes\R/\calO_{K,\bar c}$  and in fact is 
isogenous to the product $E\times E$. In particular there is 
an identification $\End^{o}(A)\simeq D\otimes K$. Consider the 
left ideal $\gote=\End(A)\cap\End^{o}(A)(1-e_{\jmath})$ where $e_{\jmath}$ is the 
idempotent \eqref{eq:idempotent} 
attached to the embedding $\jmath:K\hmap D$ associated to $\tau$ 
and let $\calE=\calA[\gote]^{o}$ be the connected component of the subgroup scheme 
of $\calA$ killed by $\gote$. Note that since 
$\jmath(\calO_{K,\bar c})$ and $e_{\jmath}$ commute, the order $\calO_{K,\bar c}$ acts 
on $\calE$.

\begin{pro}\label{teo:idgrsch}
    $\calA=\calE\otimes_{\calO_{K,\bar c}}\calR_{1}$ as group schemes.
\end{pro}

\pf Let $S$ be any scheme of definition for $A$. Over any $S$-scheme 
$T$ there is an obvious map 
$(\calE\otimes_{\calO_{K,\bar c}}\calR_{1})(T)\map\calA(T)$ which is surjective 
because $\calE\otimes_{\calO_{K,\bar c}}\calR_{1}$ contains two independent 
abelian schemes of dimension 1, $\calE$ and any translate of it by an 
$r\in{\cal R}_{1}-\calO_{K,\bar c}$. To show that the map is injective, it is 
enough to do so over an algebraically closed field. Over $\C$ we have 
$\calE(\C)=E$ and thus
$(\calE\otimes_{R_{\bar c}}\calR_{1})(\C)=E\otimes_{\calO_{K,\bar c}}\calR_{1}=
(K\otimes\R\otimes_{\calO_{K,\bar c}}\calR_{1})/\calR_{1}=A$. \qed

\begin{dfn}\label{th:testpair}
   Let $p$ be an odd prime number, $\mcd p{N\De}=1$. 
   A \emph{$p$-ordinary test triple} for $\GaE(\De,N)$ is a 
   triple $(\tau, v, e)$, where $\tau\in\CM_{\De,K}$, 
   $v$ is a finite place dividing $p$ in a finite extension $L\supseteq\Q$
   and $e\in D\otimes F$ is the idempotent associated to a real 
   quadratic subfield $F\subset D$ pointwise fixed by the positive 
   involution, such that
     \begin{enumerate}
       \item $FK\subseteq L$; 
       \item the CM curve $E_\tau$ or QM-abelian surface $A_\tau$ 
	       has ordinary good reduction modulo $\gotp_v$;
       \item if $w$ is the restriction of $v$ to  $F$ then 
                $e\in\calR_{1}\otimes_{\Z}\calO_{(w)}$. 
    \end{enumerate}
    Furthermore, a $p$-ordinary test triple $(\tau, v, e)$ is  said 
    \emph{split} if $p$ splits in $F$.
\end{dfn}

\noindent Let us observe that:
\begin{enumerate}
    \item  the ordinarity hypothesis implies that $p$ splits in $K$;
    \item the idempotent $e$ plays no role in the split case and can be omitted in that case;
    \item  the explicit description \eqref{eq:idempotent} of $e$ 
           shows that the third condition above is equivalent to 
           $\mcd p{\bar c\de_{F}}=1$ where $\bar c$ is the minimal 
           conductor of $F$; 
    \item  for a $p$-ordinary triple $(\tau,v,e)$ for $\GaU(\De,N)$ the point 
           $x\in X_1(\De,N)$ represented by $\tau$ is a smooth point in 
	   $\calX_1(\calO_{(v)})$. This is clear for $D$ split and follows for 
	   instance from $\cite[Theorem~1.1]{Jord86}$ in the non-split case.
\end{enumerate}

\begin{pro}\label{teo:anypworks}
    Let $p$ be an odd prime number, $\mcd p{N\De}=1$. There exist split 
    $p$-ordinary triples for $\GaE(\De,N)$.
\end{pro}

\pf Since any two positive involutions \eqref{eq:involution} are
conjugated in $D$, up to a different choice of maximal order we are 
reduced to the Hashimoto model. Up to replacing $p_{o}$ in \eqref{eq:condonHS} 
in its congruence class modulo $8N\De p$, we may assume also that 
$\vvec{p_{o}}p=1$. Thus the subfield $F=\Q\oplus\Q j\subset D_{H}$ is 
pointwise fixed by the involution, has discriminant prime to $p$ and $p$ 
splits in it. Finally, the minimal conductor of the embedding 
$\sqrt{p_o}\mapsto j\in F$ is prime to $p$ since $j\in\calR_{H,N}$.\qed

The decomposition $D=K\oplus Ku$ associated to a choice of 
$\tau\in\CM_{\De,K}$ is also an orthogonal decomposition under 
the non-degenerate pairing $(x,y)_{D}=\tra(x\bar{y})$. 
Note that here $u^{2}>0$ since the norm is indefinite.
We shall be concerned with the algebraic group of 
similitudes of $(\cdot,\cdot)_{D}$, i.e.
$$
\GO(D)=\left\{\hbox{$g\in\GL(D)$ such that $(gx,gy)_{D}=\nu_{0}(g)(x,y)_{D}$
for all $x,y\in D$}\right\}.
$$
The structure of the group $\GO(D)$ is well understood, e.g. 
\cite[\S1.1]{Harris93}, \cite[\S7]{HaKu91}. Let $\mathbf{t}\in\GO(D)$ 
be the involution $\mathbf{t}(d)=\bar{d}$. Then 
$\GO(D)=\GO^{o}(D)\ltimes<\mathbf{t}>$, where $\GO^{o}(D)$ is the 
Zariski connected component described by the short exact sequence of 
algebraic groups
\begin{equation}
    1\mmap\G_{m}\mmap
    D^{\times}\times D^{\times}\stackrel{\varrho}{\mmap}
    GO^{o}(D)\mmap1
    \label{eq:sesGOD}
\end{equation}
where $\G_{m}$ is embedded diagonally and 
$\varrho(d_{1},d_{2})(x)=d_{1}xd_{2}^{-1}$. The norm $\nu$ restricts 
to $N_{K/\Q}$ and $-u^{2}N_{K/\Q}$ on $K$ and  $Ku$ respectively, and 
$\GO(K)^o\simeq\GO^o(Ku)\simeq R_{K/\Q}\G_{m,K}$ where the isomorphism is 
given by left multiplication. Thus, the subgroup of $\GO^{o}(D)$ that 
preserves the splitting $D=K\oplus Ku$ can be identified with the group
$$
G(O(K)\times O(Ku))^o=
\left\{\hbox{$(k_{1},k_{2})\in(R_{K/\Q}\G_{m,K})^2$
such that $N_{K/\Q}(k_{1}\inv k_{2})=1$}\right\}
$$
and there is a commutative diagram
\begin{equation}
\begin{CD}
	K^{\times}\times K^{\times} @>{\al}>> G(O(K)\times O(Ku))^o  \\
	@V{\jmath\times\jmath}VV @VVV \\
	D^{\times}\times D^{\times} @>{\varrho}>> \GO^{o}(D)\\
\end{CD}
	\label{eq:similitudes}
\end{equation}
where $\al(k_{1},k_{2})=(k_{1}k_{2}^{-1},k_{1}\bar{k}_{2}^{-1})$.

We will normalize the complex coordinates in 
$D_{\infty}^{\tau}=(K\oplus Ku)\otimes\R$ as follows. The standard normalized 
embedding 
$\jmath^{\rm st}:\Q(\sqrt{-1})\hookrightarrow\M_2(\Q)$ with fixed point $i\in\gotH$ 
defines a splitting $M_{2}(\R)^{i}=\C\oplus\C^\perp$ with
$\C=\R\smallmat 1{}{}1\oplus\R\smallmat{}{-1}1{}$ and 
$\C^\perp=\C\smallmat {}11{}=\R\smallmat {}11{}\oplus\R\smallmat{-1}{}{}1$.
Define standard complex coordinates $z_1^{\rm st}$, $z_2^{\rm st}$ in 
$D_{\infty}$ by the identity
\begin{equation}
\Phi_\infty(d)=z_1^{\rm st}+z_2^{\rm st}\smallmat {}11{}.
    \label{eq:standardcoord}
\end{equation}
The $\R$-linear extensions of the embeddings $\jmath^{\rm st}$ and 
$\Phi_\infty\circ\jmath$ are conjugated in $M_{2}(\R)$, namely 
$\Phi_{\infty}\circ\jmath=d_\infty\jmath^{\rm st}d_\infty^{-1}$ where
$d_\infty=\smallmat{y^{1/2}}{sy^{1/2}}{}{y^{-1/2}}$
and $\tau=s+iy$. So we define normalized coordinates $z_{1}$ and $z_{2}$ in 
$D_\infty^\tau$ by the identity
$$
z_{i}(d)=z_{i}^{\rm st}
(\Phi_\infty^{-1}(d_\infty^{-1})d\Phi_\infty^{-1}(d_\infty)),
\qquad\hbox{for all $d\in D_{\infty}$,\quad $i=1$,$2$}.
$$


\section{Some differential operators}

\subsection{Preliminaries.}\label{se:KSprel}
We briefly review some basic facts about the Kodaira-Spencer map and
the Gau{\ss}-Manin connection. For more details see \cite{Katz70, KatOda68}. 
The \emph{Kodaira-Spencer class} of a composition of smooth morphisms of schemes 
$X\stackrel{\pi}{\map}S\map T$ is the element in 
$H^1(X,\dual{(\Om^1_{X/S})}\otimes \pi^*\Om^1_{S/T})$ 
arising from the canonical exact sequence. 
\begin{equation}
0\mmap\pi^*\Om^1_{S/T}\mmap\Om^1_{X/T}\mmap\Om^1_{X/S}\mmap 0
\label{eq:canexseq}
\end{equation}
by local freeness of the sheaves $\Om^1$.
The \emph{Kodaira-Spencer map} is the boundary map
$$
\KS:\pi_*\Om^1_{X/S}\mmap R^1\pi_*(\pi^*\Om^1_{S/T})\simeq
\Om^1_{S/T}\otimes R^1\pi_*\calO_X
$$
in the long exact sequence of derived functors obtained from \eqref{eq:canexseq} by pushing down. 
Under the natural maps
$H^1(X,\dual{(\Om^1_{X/S})}\otimes \pi^*\Om^1_{S/T})\map
H^0(S,R^1\pi_*(\dual{(\Om^1_{X/S})}\otimes \pi^*\Om^1_{S/T}))\map
H^0(S,\Om^1_{S/T}\otimes R^1\pi_*\calO_X\otimes\dual{(\pi_*\Om^1_{X/S})})$ the Kodaira-Spencer class maps to the Kodaira-Spencer map.

The $q$-\emph{th relative de Rham cohomology sheaf} of $X/S$ is defined as
$\derham{q}(X/S)=\R^q\pi_*(\Om^\bullet_{X/S})$ (hypercohomology). Following \cite{KatOda68}, the \emph{Gau{\ss}-Manin connection}
$$
\nabla\colon\derham{q}(X/S)\mmap\Om^1_{S/T}\otimes_{\calO_S}\derham{q}(X/S).
$$
can be seen as the differential
$d_1^{0,q}\colon E_1^{0,q}\map E_1^{1,q}$ in the spectral sequence
defined by the finite filtration
$F^i\Om^\bullet_{X/T}=\Imm(\Om^{\bullet-i}_{X/T}
\otimes_{\calO_X}\pi^*\Om^i_{S/T}\mmap\Om^\bullet_{X/T})$,
with associated graded objects
$\gr^i(\Om^\bullet_{X/T})=\Om^{\bullet-i}_{X/T}
\otimes_{\calO_X}\pi^*\Om^i_{S/T}$.

If $X/S=\calA/S$ is an abelian scheme with $0$-section $e_{0}$
and dual $\calA^t/S$ , 
denote
$\oom=\oom_{\calA/S}=\pi_*\Om^1_{\calA/S}={e_0}^*\Om^1_{\calA/S}$ 
the sheaf on $S$ of translation invariant relative $1$-forms on $A$. 
The first de Rham sheaf $\derham{1}=\derham{1}(\calA/S)$ is the
central term in a short exact sequence
\begin{equation}
0\mmap\oom\mmap\derham{1}\mmap R^1\pi_*\calO_{\calA}\mmap0
\label{eq:Hodgeseq}
\end{equation} 
(called the \emph{Hodge sequence}). By Serre duality 
\begin{equation}
\Hom_{\calO_S}(\pi_*\Om^1_{\calA/S}, R^1\pi_*(\pi^*\Om^1_{S/T}))
\simeq
\Hom_{\calO_S}(\oom_{\calA/S}\otimes\oom_{\calA^t/S},\Om^1_{S/T})
\label{eq:KSmap}
\end{equation}
and the Kodaira-Spencer map can be seen as an element of
the latter group. It can be reconstructed from the Gau{\ss}-Manin connection as the composition
\begin{equation}
\oom_{\calA/S}\hmap\derham{1}\stackrel{\nabla}{\mmap}\derham{1}\otimes\Om^1_{S/T}
\mmap\dual{\oom_{\calA^t/S}}\otimes\Om^1_{S/T}.
\label{eq:KSfromGM}
\end{equation}
In fact, when $\calA/S\simeq\calA^t/S$ is principally polarized,
the Kodaira-Spencer map becomes a \emph{symmetric} map
$\KS\colon\Sym^2(\oom)\map\Om^1_{S/T}$, \cite[Section III.~9]{FalCha90}.

Let $(\pr{S},i_o)$ be a smooth closed reduced subscheme of $S$
and consider the commutative pull-back diagram of $T$-schemes
$$
	\begin{CD}
	\pr{X} @>i>> X \\
	@VV{\pr{\pi}}V  @VV{\pi}V \\
	\pr{S} @>{i_o}>> S \\
	\end{CD}.
$$
Since also $\pr{\pi}$ is smooth, we can consider the Kodaira-Spencer class, or map, $\pr{\KS}$ attached to the morphisms $\pr{X}\stackrel{\pr{\pi}}{\map}\pr{S}\map T$.
When $X=\calA$ is a principally polarized abelian scheme,
$\pr{\KS}\in\Hom_{\calO_{\pr{S}}}({\pr{\oom}}^{\otimes 2},\Om^1_{\pr{S}/T})$
as in \eqref{eq:KSmap}, where
${\pr{\oom}}=\oom_{\pr{\calA}/\pr{S}}=\pi_*\Om^1_{\pr{\calA}/\pr{S}}=
{e^\prime_0}^*\Om^1_{\pr{\calA}/\pr{S}}$.
Since $i^*\pi^*\Om^1_{S/T}={\pi^\prime}^*i_0^*\Om^1_{S/T}$ and 
$i^*\Om_{X/S}\simeq\Om_{X^\prime/S^\prime}$ canonically, applying $i^{*}$ to \eqref{eq:canexseq}  yields an exact sequence
$$
0\mmap{\pr{\pi}}^*i_0^*\Om^1_{S/T}\mmap
i^*\Om^1_{X/T}\mmap\Om^1_{\pr{X}/\pr{S}}\mmap 0,
$$
hence an element
$\KS^*\in\Ext^1_{\calO_{\pr{X}}}(\Om^1_{\pr{X}/\pr{S}},
{\pr{\pi}}^*i_0^*\Om^1_{S/T})$. The composition
$\pr{S}\stackrel{i_0}{\map}S\map T$ defines a canonical surjective map
${\pr{\pi}}^*i_0^*\Om^1_{S/T}\map{\pr{\pi}}^*\Om^1_{\pr{S}/T}$.
In the same way, we get a surjective map $i^*\Om^1_{X/T}\map\Om^1_{\pr{X}/T}$.
These data define a commutative diagram of $\calO_{X^\prime}$-modules
$$
	\begin{CD}
	0 @>>> {\pr{\pi}}^*i_0^*\Om^1_{S/T} @>>> i^*\Om^1_{X/T} @>>>
	\Om^1_{\pr{X}/\pr{S}} @>>> 0\\
	@. @VVV @VVV @| \\
	0 @>>> {\pr{\pi}}^*\Om^1_{\pr{S}/T} @>>> \Om^1_{\pr{X}/T} @>>>
	\Om^1_{\pr{X}/\pr{S}} @>>> 0\\
	\end{CD}
$$
Standard diagram-chasing shows that
$\KS^*\mapsto\pr{\KS}$ under the canonical map of $\Ext^1$ groups.
The following result follows easily from the definitions.

\begin{pro}\label{th:KSforpullb}
	Let $\calA/S$ be an abelian scheme with $S$ smooth
	over $T$, $(\pr{S},i_0)$ a
	closed $T$-smooth subscheme of $S$ and
	$\pr{\calA}=\calA\times_S\pr{S}$. Let
	$\KS\colon\oom^{\otimes2}\map\Om^1_{S/T}$ and
	$\KS^\prime\colon{\pr{\oom}}^{\otimes2}\map\Om^1_{\pr{S}/T}$
	be the corresponding Kodaira-Spencer maps. Then
	$\pr{\KS}=\iota_0\circ i_0^*\KS$, where
	$\iota_0\colon i_0^*\Om^1_{S/T}\map\Om^1_{\pr{S}/T}$
	is the canonical pull-back map.
\end{pro}

Let again $X=\calA$ be an abelian scheme, and let $\phi\colon\calA\map\calA$
be an $S$-isogeny (i.~e., a surjective endomorphism such that $\pi\phi=\pi$). The pull-back
$\phi^*\Om^\bullet_{\calA/T}\map\Om^\bullet_{\calA/T}$ respects filtrations.
Thus we have maps $\phi^*(F^i/F^j)\map F^i/F^j$ for all $i\leq j$
because the sheaves $F^i$ are locally free. In particular, there is a
map of short exact sequences
$$
\begin{CD}
	0 @>>> \phi^*\gr^{p+1} @>>> \phi^*(F^p/F^{p+2}) @>>>
	\phi^*\gr^p @>>> 0  \\
	@.  @VVV  @VVV  @VVV  \\
	0 @>>> \gr^{p+1} @>>> F^p/F^{p+2} @>>> \gr^p @>>> 0 \\
\end{CD}
$$
where the bottom row is the tautological exact sequence of graded objects and the top
row is obtained applying $\phi^*$ to it (again, it remains exact because the
sheaves are locally free).
Since $\phi$ is surjective, $\pi_*\phi^*=\pi_*$ as
functors and the previous diagram yields a map of derived functors long exact sequences
\begin{equation}
	\begin{CD}
		\ldots @>>> R^{p+q}\pi_*\gr^p @>>> R^{p+q+1}\pi_*\gr^{p+1} @>>>
\ldots \\
		@. @VV{[\phi]_{p,q}}V @VV{[\phi]_{p+1,q}}V \\
		\ldots @>>> R^{p+q}\pi_*\gr^p @>>> R^{p+q+1}\pi_*\gr^{p+1} @>>>
\ldots \\
	\end{CD}
	\label{cd:four}
\end{equation}

\begin{pro}\label{th:GMcomm}
	Let $\calA/S$ be an abelian scheme with $S$ 	smooth over $T$.
	The algebra $\End_S(\calA)$ acts linearly on the sheaves $\derham{q}(\calA/S)$. 
	If $\phi\in\End_S(\calA)$ acts as $[\phi]$, then
	$$
	\nabla\circ[\phi]=(1\otimes[\phi])\nabla.
	$$
\end{pro}

\pf Let $\phi\in\End_S(\calA)$ be an isogeny. The endomorphism $[\phi]$ of
$\derham{q}(\calA/S)$ attached to $\phi$ is the vertical map $[\phi]_{0,q}$ in diagram
\eqref{cd:four} at $R^q\pi_*\gr^0$.  Under the identification 
$R^{q+1}\pi_*\gr^1=R^{q+1}\pi_*(\pi^*\Om^1_{S/T}\otimes_{\calO_{\calA}}\Om^{\bullet-1}_{\calA/T})=
\Om^1_{S/T}\otimes_{\calO_S}R^{q+1}\pi_*(\Om^{\bullet-1}_{\calA/T})=
\Om^1_{S/T}\otimes_{\calO_S}\derham{q}(\calA/S)$
the Gau{\ss}-Manin connection is the connecting homomorphism for the
tautological exact sequence of graded objects, i.÷e. either horizontal connecting homomorphism in 
\eqref{cd:four} at $p=0$ and also $[\phi]_{1,q}=1\otimes[\phi]_{0,q}$ since $\phi$ acts trivially on $S$. The formula follows.

Let $s\in S$ be a geometric point and $A_s$ the fiber at $s$.
Without loss of generality we may assume that $S$ is connected and Grothendieck's rigidity lemma \cite{Mum65} implies that the canonical map 
$\End_S(\calA)\map\End(A_s)$ is injective. It follows that there exist 
division algebras $D_1,\ldots,D_t$ such that 
$\End_S(\calA)$ is identified to a subring of 
${\rm M}_{n_1}(D_1)\times\cdots\times{\rm M}_{n_t}(D_t)$. The latter 
algebra is spanned over $\Q$ by the invertible elements, so the result 
follows by linearity. \qed

\subsection{Computations over $\C$}\label{se:KSoverC}
In order to compute explicitely the Kodaira-Spencer map for a complex family, i.e. when $T=\Spec(\C)$, it is more convenient to appeal to GAGA principles, work in the analytic category and follow \cite{Katz76, Harris81}. If $\calA/S$ is a principally polarized family of abelian varieties
over the smooth complex variety $S$ and $U\subset S$ is an open set, the choice of a section  
$\si\in H^0(U,\dual{(\Om^1_S)})$ defines a map 
$\varrho_\si\colon H^0(U,\oom)\map H^0(U,\dual\oom)$
by the composition
$$
H^0(U,\oom)\hmap
H^0(U,\derham{1})\stackrel{\nabla}{\map}
H^0(U,\derham{1}\otimes\Om^1_{S/T})\stackrel{1\otimes\si}{\map}
H^0(U,\derham{1})\stackrel{\varrho}{\map}
H^0(U,\dual\oom),
$$
where $\varrho$ is induced by the polarization pairing 
$\scaldR\cdot\cdot\colon\derham{1}\otimes\derham{1}\map\calO_S$.
The association $\si\mapsto\varrho_\si$ defines a map
$\dual{(\Om^1_{S/T})}\map\Hom(\oom,\dual\oom)$
whose dual is the Kodaira-Spencer map $\KS$.
By \'{e}tale-ness, the actual computation of the map $\KS$ can be
obtained applying the above procedure to the pullback of the family
$\calA/S$ on the universal cover of $S$. For instance, for the universal family \eqref{eq:univEC}
\begin{equation}
    \KS(d\ze^{\otimes 2})=\frac{1}{2\pi i}dz.
    \label{eq:Kodforell}
\end{equation}
where $\ze$ is the standard complex coordinate in the elliptic curve $E_z=\C/\Z\oplus\Z z$, $z\in\gotH$.

We follow this approach to compute the Kodaira-Spencer map for
the universal complex family of QM abelian surfaces over $X_\Ga$
using the Shimura family $\calA^{\rm Sh}/\gotH$ of
\eqref{eq:QMtori} in terms of the arithmetic of the maximal order 
$\calR_{1}$.
Let $\underline{r}=\{r_1,\dots,r_4\}$ be a symplectic basis of $\calR_{1}$.
By linear extension, the real dual basis  $\{\dual{r_1},\dots,\dual{r_4}\}$
of $\dual{D_\infty}$ is a basis of
$\Hom(\calR_{1}\otimes_\Z\C,\C)\simeq\Derham{1}(D_\infty^z/\calR_{1})$.
Thus, the elements $\dual{r_1},\dots,\dual{r_4}$ define global
$C^\infty$-sections of $\derham{1}(\calA^{\rm Sh}/\gotH)$ with
constant periods, hence $\nabla$-horizontal. If $H$
denotes the $\C$-span of these sections, there is an isomorphism
$\derham{1}(\calA^{\rm Sh}/\gotH)=H\otimes_\C\calO_{\gotH}$.
In terms of this trivialization, $\nabla=1\otimes d$, where
$d$ is the exterior differentiation. Also,
\begin{equation}
  \scaldR{\dual{r_i}}{\dual{r_j}}=\frac{1}{2\pi i}B_t(r_j,r_i),
\qquad i,j=1,\ldots,4,  
    \label{eq:derhampair}
\end{equation}
where the $2\pi i$ factor  accounts for the difference of Tate twists
between singular and algebraic de Rham cohomology, e.~g. \cite[\S1]{Del82}.
Let $\ze_1$ and $\ze_2$ denote the standard coordinates in $\C^2$.

\begin{pro}\label{th:KSoverC}
	$\left(\KS(d\ze_i\otimes d\ze_j)\right)_{i,j=1,2}=
	\frac1{2\pi i}\smallmat100{\De}\,dz.$
\end{pro}

\pf Write
$$
\left(\begin{array}{c}
	d\ze_1  \\
	d\ze_2
\end{array}\right)=
\Pi_{\underline{r}}(z)
\left(
\begin{array}{c}
	\dual{r_1}  \\
	\vdots  \\
	\dual{r_4}
\end{array}
\right)
$$
where $\Pi_{\underline{r}}(z)$ is the period matrix computed in terms of the
basis $\underline{r}$. Using \eqref{eq:derhampair} and the definitions we first obtain
$$
\varrho_\si\left(\begin{array}{c}
	d\ze_1  \\
	d\ze_2
\end{array}\right)=
\frac1{2\pi i}\frac{d\Pi_{\underline{r}}(z)}{dz}
\left(\begin{array}{cc}
	0 & I_2  \\
	-I_2 & 0
\end{array}\right)
\left(\begin{array}{c}
	r_1  \\
	\vdots  \\
	r_4
\end{array}\right)\si(dz),
$$
and finally
$$
	\left(\KS(d\ze_i\otimes d\ze_j)\right)_{i,j=1,2}=
	\frac{1}{2\pi i}\frac{d\Pi_{\underline{r}}(z)}{dz}
	\left(\begin{array}{cc}
		0 & I_2  \\
		-I_2 & 0
	\end{array}\right)
	{}^t\Pi_{\underline{r}}(z)\,dz.
$$
To obtain the final formula, we make use of the Hashimoto model with 
$N=1$.
In terms of the symplectic basis
$\underline{\eta}=\{\eta_1,\dots,\eta_4\}$
of theorem \ref{th:Hashimoto}
\begin{equation}
	\Pi_{\underline{\eta}}(z)=\left(
	\begin{array}{cccc}
	\frac{\varpi^-}{2\srh}(\al^+a\De z+1) &
	-\frac{1}{\srh}(\al^{+}a\De z+1)
	& z & \frac12\al^+z \\
		\\
	\frac{\al^+}{2\srh}a\De (z-\al^-) &
	\frac{1}{\srh}\De(-z+\al^-a) & 1 & \frac12\al^- \\
	\end{array}\right)
	\label{eq:permat}
\end{equation}
where $\al^{\pm}=1\pm\srh$.
Plugging these values into the previous formula yields the result. \qed

\subsection{Maass operators.}\label{ss:maass}
When $D$ is split. the universal family 
\eqref{eq:univEC} defines the
line bundle  $\oom=\oom_{\calE_N/\calY_1(N)}$
on the Zariski open set $\calY_1(N)$ complement of the
cusp divisor $C$ in $\calX_1(N)$. The Kodaira-Spencer map
$\KS\colon\oom^{\otimes2}\isom\Om^1_{\calY_1(N)}$
is an isomorphism.

\begin{thm}\label{th:KSextended}
	The line bundle $\oom$ extends uniquely to a line bundle, still
        denoted
	$\oom$, on the complete curve $\calX_1(N)$ and the Kodaira-Spencer
	isomorphism extends to an isomorphism
	$$
	\KS\colon\oom^{\otimes2}\iisom\Om^1_{\calX_1(N)}(\log C).
	$$
\end{thm}

\pf See \cite{Katz73} and also \cite[section~10.~13]{KatMaz85} where the
extension property is discussed for a general representable
moduli problem. \qed

If $D$ is not split, the universal family \eqref{eq:univQMAV}
of QM-abelian surfaces defines the sheaf
$\oom=\oom_{\calA_{\De,N}/\calX_1(\De,N)}$
and the Kodaira-Spencer map is a surjective map
$\KS\colon\Sym^2\oom\map\Om^1_{\calX_1(\De,N)}$.

Let $p$ be a prime such that $(p,N\De)=1$ and let $v$ be a
place of a number field $L$ dividing $p$. The algebra 
$\calR_{1}\otimes_\Z\calO_v$ acts contravariantly and 
$\calO_{v}$-linearly on $\oom_{v}=\oom\otimes\calO_v$ by pull-back. 
For any geometric point $s\in\calX_1(\De,N)\otimes\calO_v$ 
and any non-trivial idempotent $e\in\calR_{1}\otimes_\Z\calO_v$ 
there is a non-trivial decomposition $H^0(A_s,\Om^1_{A_s/k(s)})=
eH^0(A_s,\Om^1_{A_s/k(s)})\oplus(1-e)H^0(A_s,\Om^1_{A_s/k(s)})$. 
Therefore the subsheaf $e\oom_{v}$ is a line subbundle.
Let $e\oom_v\circ\invol{e}\oom_v\subseteq\Sym^2\oom_v$ be the line 
bundle image
of $e\oom_v\otimes\invol{e}\oom_v$ under the natural
map $\oom^{\otimes 2}_v\map\Sym^2{\oom_v}$.

\begin{thm}\label{th:Ltbundles}
	If $p$, $v$ and $e$ are as above, then the Kodaira-Spencer map defines an 
	isomorphism
	$$
	\KS\colon e\oom_v\otimes\invol{e}\oom_v\mmap
	\Om^1_{\calX_1(\De,N)/\calO_v}
	$$
	of line bundles on $X_{1}(\De,N)$ defined over $\calO_{v}$.
\end{thm}

\pf 
We claim that the action of $r\otimes\la\in\calR_{1}\otimes\calO_v$ on the universal family \eqref{eq:univQMAV} base-changed to $\calO_v$ gives rise to a commutative diagram
\begin{equation}
	\begin{CD}
	 	\oom _v @>>> \Om^1_{\calX_1(\De,N)/\calO_v}\otimes\dual{\oom_v} \\
	 	@VV{r\otimes\la}V @VV{1\otimes\invol r\otimes\la}V \\
	 	\oom _v @>>> \Om^1_{\calX_1(\De,N)/\calO_v}\otimes\dual{\oom_v}
	\end{CD}
	\label{cd:five}
\end{equation}
Indeed, under the Serre duality identification
$R^1\pi_*\calO_{\calA}\simeq\dual{\oom_{\calA/S}}$ for a principally polarized abelian scheme $\calA/S$ the actions of $\End_S(\calA)$ correspond up to Rosati involution. The
commutativity of the diagram \eqref{cd:five} follows from proposition
\ref{th:GMcomm} and \eqref{eq:KSfromGM}.

For an idempotent $e\in\calR_{1}\otimes_\Z\calO_v$, diagram
\eqref{cd:five} defines a map
$e\oom_v\map\Om^1\otimes\invol{e}(\dual{\oom_v})$
which can be shown to be an isomorphism by the same deformation theory
argument in \cite[Lemma 6]{DiaTay94}. This is enough to conclude,
because the sheaves $\invol{e}(\dual{\oom_v})$
and $\invol{e}\oom_v$ are dual of each other.\qed

\begin{rems}\label{re:LoverC}
	\rm
	\begin{enumerate}
	 \item We proved theorem \ref{th:Ltbundles} for $p$-adically complete 
	       rings of scalars. In fact the projectors onto the quadratic 
	       subfields of $D$ are defined over the $p$-adic localizations 
	       of their rings of integers for almost all $p$.
               Thus, in these cases, $e\oom$ and the Kodaira-Spencer
               isomorphism are defined over the subrings $\calO_{(v)}\subset\C$.
        \item  If the $\pr e=ded^{-1}$ are conjugated in $\calR_1\otimes B$ 
	       for some ring $B$,  then the action of $d$ on 
	       $\oom\otimes B$ defines an isomorphism of      
                $e\oom$ with $\pr e\oom$ over $B$.
        \item  In the complex case the isomorphism of theorem \ref{th:Ltbundles}
	       can be checked by a straightforward application of the 
	       computation in section \ref{se:KSoverC}. For instance, in the 
	       Hashi\-mo\-to model for $N=1$ of theorem \ref{th:Hashimoto} let 
	       $d=ai+bj+cij\in D_{H}$ with 
	       $\de=d^{2}=-a^{2}\De+b^{2}p_{o}+c^{2}\De p_{o}\in\Q$ and let 
	       $e\in D\otimes_{\Q}\Q(\sqrt{\de})$ be the idempotent giving 
	       the projection onto $\Q(d)$. Then
	       $$
	       e=\frac{1}{2\sqrt{\de}}
	       \left(
	       \begin{array}{cc}
	          \sqrt{\de}+b\sqrt{p_{o}}  & -a+c\sqrt{p_{o}}  \\
	          (a+c\sqrt{p_{o}})\De & \sqrt{\de}-b\sqrt{p_{o}}
	       \end{array}
	       \right)
	       $$
	       and
	       $$
	       \invol e=\frac{1}{2\sqrt{\de}}
	       \left(
	       \begin{array}{cc}
	          \sqrt{\de}+b\sqrt{p_{o}}  & a+c\sqrt{p_{o}}  \\
	          (-a+c\sqrt{p_{o}})\De & \sqrt{\de}-b\sqrt{p_{o}}
	       \end{array}
	       \right).
	       $$
	       Therefore $e\oom\circ\invol e\oom$ is generated over $\gotH$
	       by the global section
	       $$
	       (\sqrt{\de}+b\sqrt{p_{o}})^{2}d\ze_{1}\circ d\ze_{1}+
	       (c^{2}p_{o}-a^{2})d\ze_{2}\circ d\ze_{2}+
	       2(\sqrt{\de}+b\sqrt{p_{o}})c\sqrt{p_{o}}d\ze_{1}\circ d\ze_{2}
	       $$
	       whose image under the Kodaira-Spencer map $\KS$ is, 
	       by proposition \ref{th:KSoverC}, 
	       \begin{equation}
	           \frac{1}{\pi i}(\sqrt{\de}+b\sqrt{p_{o}})dz
		   \in\Ga(\gotH,\Om^{1}_{\gotH}).
	           \label{eq:imkasec}
	       \end{equation}
                Since $\de\neq p_{o}b^2$ (else $p_o=(a/c)^2\in\Q$ which is impossible) 
                the section \eqref{eq:imkasec} does not vanish and the Kodaira-Spencer map is an isomorphism.
	   \end{enumerate}
\end{rems}

\begin{notat}
    \rm We will denote $\calL$ either the line bundle $\oom_v$ on 
    $\calY_{1}(N)$ or the line bundle $e\oom_v$ on $\calX_{1}(\De,N)$ 
    for some choice of idempotent $e$ satisfying the hypotheses of 
    theorem \ref{th:Ltbundles} and such that $\invol e=e$. In either 
    case the Kodaira-Spencer map gives an isomorphism
    $$
    \KS\colon\calL^{\otimes 2}\stackrel{\sim}{\mmap}\Om^{1}.
    $$
    With an abuse of notation we will denote also $\calL$ the pullback of 
    the complexified bundle to $\gotH$ under the natural quotient maps.
\end{notat}

If $\ga\in\GaU(\De,N)$ the identities
$\Z^{2}\vvec{\ga\cdot z}{1}=\inv{j(\ga,z)}\Z^{2}\vvec{z}{1}$ and
$\Phi_\infty(\calR_{1})\vvec{\ga\cdot z}{1}=
\inv{j(\ga,z)}\Phi_\infty(\calR_{1})\vvec{z}{1}$
as subsets of $\C$ (in the split case) and of $\C^{2}$ (in the non-split 
case) respectively, show that
the natural action of $\GaU(\De,N)$ on $\oom$ over $\gotH$ is scalar
multiplication by the automorphy factor. Thus the $\GaU(\De,N)$-action 
extends to a $\SL_2(\R)$-homogeneous structure on $\oom$, 
and on $\Sym^2(\oom)$ as well.
Also, in the non-split case the fiber identifications
induced by the action are $D\otimes\C$-contravariant and
since the line bundle $\calL$ is defined using the $D$ action on $\oom$, 
it is an homogeneous
line subbundle of $\Sym^2(\oom)$.

Let $n\in\Z$ and let $V_n$ be the $1$-dimensional representation of
$\C^\times$ given by the character $\chi_n(z)=z^n$. Let
$\calV_n=V_n\times\gotH$
the homogeneous line bundle on $\gotH$ with action
$g\cdot(v,z)=(\chi_n(j(g,z))v,g\cdot z)$. Since $-1\notin\GaU(\De,N)$ and
$\GaU(\De,N)$ has no elliptic elements, the
quotient $\GaU(\De,N)\bs\calV_n$ is a line bundle on $X_{1}(\De,N)$
which we shall
denote $\calV_n$ again. Pick $v_n\in V_n$, $v_n\neq0$, and
let $\tilde{v}_n=(v_n,z)$ be the corresponding global constant
section of $\calV_n$ over $\gotH$. Also, let $s(z)$ be the global
section of $\calL$ over $\gotH$, defined up to a sign, 
normalized so that
\begin{equation}
    \KS(s(z)^{\otimes 2})=2\pi i\,dz.
    \label{eq:kanormsec}
\end{equation}
Then $\calL^{\otimes k}=\calO_{\gotH}s(z)^{\otimes k}$ for all 
$k\geq1$  and 
there are identifications of homogeneous complex line bundles over
both $\gotH$ and $X_1(\De,N)$
\begin{equation}
	\calV_2\iisom\Om^1,
	\ \tilde{v}_2\mapsto 2\pi i\,dz\qquad
	\hbox{and}\qquad
	\calV_{k}\iisom{\calL}^{\otimes k},
	\ \tilde{v}_{k}\mapsto{s(z)}^{\otimes k}.
	\label{eq:Eknsplit}
\end{equation}
These identifications preserve holomorphy
and are compatible with tensor products and 
the Kodaira-Spencer isomorphisms. Note that $s(z)=\pm 2\pi i\,d\ze$ in the split 
case by \eqref{eq:Kodforell}, and see remark \ref{re:LoverC}.3 for the non-split case.

Following \cite{Katz76}, we shall define
differential operators associated to splittings of the Hodge seguence
\eqref{eq:Hodgeseq} where $\calA/S$ is either the universal elliptic
curve $\calE_N/\calY_1(N)$ or the universal QM-abelian surface
$\calA_{\De,N}/\calX_1(\De,N)$.

Over the associated differentiable manifold, which amounts to tensoring with the
sheaf of $\calO_S$-algebras $\calO^\infty_S=\calC^\infty(S^{\rm an})$
and which will be denoted with an $\infty$ subscript, the Hodge
decomposition
$\calH^1_\infty=\oom_\infty\oplus\overline{\oom}_\infty$
is a splitting of the Hodge sequence with projection
$\PH\colon\calH^1_\infty\map\oom_\infty$. For each $k\geq1$, let 
$\Th_{k,\infty}^o$ be the operator defined by the composition
\begin{equation}
	\begin{CD}
	  	\Sym^k(\oom_\infty)\subset\Sym^k(\derham1)_\infty @>\nabla>>
	  	\Sym^k(\derham1)_\infty\otimes\Om^1 @>{1\otimes\inv{\KS}}>>
	  	\Sym^k(\derham1)_\infty\otimes\calL_{\infty} \\
	  	@. @. @VV{\PH^{\otimes k}\otimes1}V \\
	  	@. @. \Sym^k(\oom)_\infty\otimes\calL_{\infty} \\
	\end{CD}
	\label{eq:algmaassnsp}
\end{equation}
where the \Gauss-Manin connection $\nabla$ extends to $\Sym^k$ by the
product rule. The composition 
${\calL}^{\otimes k}\subset\oom^{\otimes k}\map\Sym^{k}(\oom)$ is injective and
let $\Th_{k,\infty}$ be the restriction of $\Th_{k,\infty}^o$ to
${\calL}^{\otimes k}_{\infty}$.

\begin{pro}\label{th:operatorrestricts}
	$\Th_{k,\infty}$ is an operator
	${\calL}^{\otimes k}_\infty\map{\calL}^{\otimes k+2}_\infty$
\end{pro}

\pf If $D$ is split then $\calL^{\otimes k}=\oom^{k}=\Sym^k(\oom)$ and 
there is nothing to prove. If $D$ is non-split, the element 
$e^{\otimes k}\in(\calR_{1}\otimes_{\Z}\calO_{(v)})^{\otimes k}$ acting 
componentwise defines a projection 
$\oom^{\otimes k}\map\calL^{\otimes k}$ which factors through 
$\Sym^k(\oom)$. By proposition \ref{th:GMcomm} the \Gauss-Manin connection
$\nabla$ commutes with $e^{\otimes k}$ and also the Hodge projection
$\PH$ is the identity on $\calL_{\infty}$ (in fact on $\om_{\infty}$). 
The result follows. \qed

The operators $\Th_{k,\infty}$  can be computed in terms of the complex
coordinate $z=x+iy\in\gotH$ and the identifications \eqref{eq:Eknsplit}. 
For any (say $\calC^\infty$) function $\phi$ on $\gotH$ let
$$
\de_k(\phi)=\frac1{2\pi i}\left(\frac{d}{dz}+\frac{k}{2iy}\right)\phi.
$$
The operator $\de_k$ was introduced, together with its higher 
dimensional analogues by Maass \cite{Maass53} and later extensively studied 
by Shimura (see \cite[Ch.~10]{Hida93} and the references cited therein).

\begin{pro}\label{th:maassincoord}
	There are commutative diagrams of $C^\infty$-bundles and
	differential operators
	$$
	\begin{CD}
		\calV_k @>\sim>> \calL^{\otimes k} \\
		@VV{\widetilde{\de}_k}V @VV{\Th_{k,\infty}}V \\
		\calV_{k+2} @>\sim>> \calL^{\otimes k+2} \\
	\end{CD}
	$$
	where $\widetilde{\de}_n(\phi\tilde{v}_n)=\de_n(\phi)\tilde{v}_{n+2}$.
\end{pro}

\pf The diagram for $D$ split is but the simplest case (dimension
$1$) of \cite[theorem~6.5]{Harris81}. The computation in the non-split case
is very similar. Let $s$ be the $\KS$-normalized section of $\calL$ as 
in \eqref{eq:kanormsec},
$\underline{\eta}=\{\eta_1,\dots,\eta_4\}$ be Hashimoto's symplectic
basis of theorem \ref{th:Hashimoto} 
and $\Pi=\Pi_{\underline{\eta}}(z)$ the period matrix as in \eqref{eq:permat}.
Since the sections $\dual{\eta_1},\ldots,\dual{\eta_4}$ are
$\nabla$-horizontal,
$$
\nabla\left(
\begin{array}{c}
	d\ze_1  \\
	  \\
	d\ze_2
\end{array}
\right)=d\Pi
\left(\begin{array}{c}
	\dual{\eta_1}  \\
	\vdots  \\
	\dual{\eta_4}
\end{array}\right)=
d\Pi\inv{\left(
\begin{array}{c}
	\Pi  \\
	  \\
	\overline{\Pi}
\end{array}
\right)}
\left(\begin{array}{c}
	d\ze_1  \\
	d\ze_2  \\
	d\bar{\ze}_1 \\
	d\bar{\ze}_2
\end{array}\right)=
\frac{\left(I_2,-I_2\right)}{z-\bar{z}}
\left(\begin{array}{c}
	d\ze_1  \\
	d\ze_2  \\
	d\bar{\ze}_1 \\
	d\bar{\ze}_2
\end{array}\right)\otimes dz.
$$
Since $s$ is in the $\C$-span of $d\ze_1$ and $d\ze_2$, 
$\nabla(s)=\left(\frac1{z-\bar{z}}s+s_0\right)\otimes dz$
with $\PH(s_0)=0$. Plugging this into
$\Th_{k,\infty}(\phi s^{\otimes k})=\PH(1\otimes\inv{\KS})
\left(\frac{d\phi}{dz}s^{\otimes k}\otimes dz+
k\phi s^{\otimes k-1}\nabla(s)\right)$ yields the result. \qed

Let $B$ be a $p$-adic algebra with $(p,\De N)=1$ and such that the 
$e$ is defined over $B$ and the isomorphism of theorem \ref{th:Ltbundles} 
holds for the sheaves base-changed to $B$.
Let $\calO^{(p)}$ be the structure sheaf of the formal scheme
$S^{(p)}=\limproj{n}(S\otimes B/p^nB)^{p\mathrm{-ord}}$ obtained
taking out the non-ordinary points in characteristic $p$. Denote
$\calM^{(p)}$ the tensorization with $\calO^{(p)}$ of the restriction to
$S^{(p)}$ of a sheaf $\calM$.

In the split case the Dwork-Katz construction \cite[\S A2.3]{Katz73}
of the unique Frobenius-stable $\nabla$-horizontal
submodule $\calU\subset\derham{1}\otimes B$ defines a splitting
$(\derham{1})^{(p)}=\oom^{(p)}\oplus\calU$ with projection
$\PFr\colon(\derham{1})^{(p)}\map\oom^{(p)}$. The construction can
be carried out in the non split case as well. If $\pr{B}$ is a 
$B$-algebra and $A_{\pr{B}}$ is  is a QM-abelian surface with ordinary reduction
and canonical subgroup $H$ (which, by ordinarity, is simply the
Cartier dual of the lift of the kernel of Verschiebung--an {\`e}tale group),
then $H\subset A[p]$ with $A[p]/H$ lifting an {\`e}tale group and for every
$\phi\in\End(A)$, $\phi(H)\subseteq H$ by connectedness. 
Thus $A/H$ is a QM-abelian
surface with a canonical embedding $\calR_1\hookrightarrow\End(A/H)$ and the
construction of the Frobenius endomorphism of $(\derham{1})^{(p)}$ and
its splitting follows. Assuming that the line bundle $\calL$ is 
defined over $B$ and following the same procedure as in
\eqref{eq:algmaassnsp} with the projection $\PH$ replaced by $\PFr$ 
yields a differential operator 
$$
\Th_{k,p}^o\colon\Sym^k(\oom^{(p)})\mmap\Sym^k(\oom^{(p)})\otimes\calL^{(p)}.
$$
Let $\Th_{k,p}$ be its restriction to $(\calL^{\otimes k})^{(p)}$.

\begin{pro}
   $\Th_{k,p}$ is an operator
   $({\calL}^{\otimes k})^{(p)}\map({\calL}^{\otimes k+2})^{(p)}$.
\end{pro}

\pf The argument is the same as in the proof of proposition
\ref{th:operatorrestricts}.
The action of the endomorphisms commutes with the pullback of forms
in the quotient $A\map A/H$ and so with the Frobenius endomorphism.
Since $\calU$ is Frobenius-stable, the endomorphisms commute with the
projection $\PFr$.\qed

Let $\ast\in\{\infty,p\}$. The operators $\Th_{k,\ast}$ can be iterated.
For all $r\geq1$ let
$$
\Th_{k,\ast}^{(r)}=
\Th_{k+2r-2,\ast}\circ\cdots\circ\Th_{k,\ast}.
$$
Since the kernel of the projecton ${\rm Pr}_\ast$ is $\nabla$-horizontal one has in fact
\begin{equation}
\label{eq:onlyonepr}
\Th_{k,\ast}^{(r)}={\rm Pr}_\ast\left((1\otimes\inv\KS)\nabla\right)^r.
\end{equation}
The operators $\Th_{k,\infty}^{(r)}$ do not preserve holomorphy
because the Hodge projection $\PH$ is not holomorphic. Similarly, the
operators $\Th_{k,p}^{(r)}$ are only defined over
$p$-adically complete ring of integers.
Nonetheless, the operators $\Th_{k,\ast}^{(r)}$ are algebraic over the
CM locus, in the following sense.
Let $x\in\calX_{1}(\De,N)(\calO_{(v)})$ be represented by a $\tau\in\gotH$ belonging to a 
$p$-ordinary test triple $(\tau,v,e)$.
Let ${\calL}(x)=x^*{\calL}$ be the algebraic fiber at $x$. 
The choice of an invariant form $\om_o$ on 
$A_x$ which generates either $H^0(A_x,\Om^1\otimes\calO_{(v)})$ (in the split case) or
$eH^0(A_x,\Om^1\otimes \calO_{(v)})$ (in the non-split case) over $\calO_{(v)}$
identifies ${\calL}(x)$ with a copy of $\calO_{(v)}$.

\begin{pro}\label{teo:Thkalgebraic}
	Let $x\in\calX_{1}(\De,N)(\calO_{(v)})$ be a point represented by a 
	$p$-ordinary test triple and let $\om_o$ be an invariant form on 
	$A_x$ as above. Then, for all         
	$r\geq1$, the operators $\Th_{k,\ast}^{(r)}$ define maps
	$$
	\Th_{k,\ast}^{(r)}(x)\colon
	H^0(\calX_{1}(\De,N)\otimes\calO_{(v)},{\calL}^{\otimes k})\mmap
	{\calL}^{\otimes k+2r}(x)\simeq
	\calO_{(v)}{\om_o}^{\otimes k+2r}.
	$$
	Moreover $\Th_{k,\infty}^{(r)}(x)=\Th_{k,p}^{(r)}(x)$.
\end{pro}

\pf The result follows, as in \cite[theorem~2.4.5]{Katz78},
from the following observation. Let $A$ be an abelian
variety isogenous over $\calO_{(v)}$ to the $g$-fold product of elliptic curves with
complex multiplications in the field $K$ and ordinary good reduction
modulo $v$. The CM splitting of the first de Rham group of $A$ is the splitting
$\Derham{1}(A/\calO{(v)})=H_{\si_1}\oplus H_{\si_2}$ where
$H_{\si_i}$ is the $\si_i$-eigenspace under the action of complex
multiplications, $I_K=\{\si_1,\si_2\}$. The Hodge
decomposition $\Derham{1}(A)\otimes\C=H^{1,0}\oplus H^{0,1}$ and
the Dwork-Katz decomposition
$\Derham{1}(A)\otimes B=H^0(A\otimes B,\Om^1)\oplus U$ for some
$p$-adic $\calO{(v)}$-algebra $B$ are both obtained from the CM splitting by
a suitable tensoring. The result follows from the algebraicity of the
Gau\ss-Manin connection and the Kodaira-Spencer map, using the expression
\eqref{eq:onlyonepr}.\qed

For all $r\geq1$ write
$$
\de_k^{(r)}=\de_{k+2r-2}\circ\cdots\circ\de_k=
\left(\frac{1}{2\pi i}\right)^r\left(\frac{d}{dz}+\frac{k+2r-2}{2iy}\right)
\circ\cdots\circ\left(\frac{d}{dz}+\frac{k}{2iy}\right)
$$
and set $\de_k^{(0)}(\phi)=\phi$.


\section{Expansions of modular forms}

\subsection{Serre-Tate theory.}\label{se:STtheory}
Let $\kk$ be any field, $(\La,\gotm)$ a complete local noetherian ring
with residue field $\kk$ and $\calC$ the category of artinian local
$\La$-algebras with residue field $\kk$. Let $\wtA$ be an abelian
variety over $\kk$ of dimension $g$. By a fundamental result of
Grothendieck \cite[2.2.1]{Oort71}, the \emph{local moduli functor}
$\calM\colon\calC\map\mbox{\bf Sets}$ which associates to each
$B\in{\rm Ob}\,\calC$ the set of deformations of $\wtA$ to $B$, is
pro-represented by $\La[[t_1,\ldots,t_{g^2}]]$.
When $\kk$ is perfect of characteristic $p>0$ and $\La=W_\kk$ is the 
ring of Witt vectors of $\kk$,
deforming $\wtA$ is equivalent to deforming its formal group, as 
precised by the Serre-Tate theory \cite[\S2]{Katz81}. 
If $\kk$ is algebraically closed and $\wtA$ is
\emph{ordinary}, an important consequence
of the Serre-Tate theory is that there is a canonical isomorphism of functors
$$
\calM\iisom\Hom(T_p\wtA\otimes T_p\wtA^t,\widehat{\G}_m),
$$
\cite[theorem 2.1]{Katz81}. 
Write $\calM=\Spf(\gotR^u)$ with universal formal deformation
$\calA^u$ over $\gotR^u$. The isomorphism endows $\calM$ with a canonical
structure of formal torus and identifies its group of characters
$X(\calM)=\Hom(\calM,\widehat{\G}_m)\subset\gotR^u$
with the group $T_p\wtA\otimes T_p\wtA^t$. Denote $q_S$ the character corresponding to 
$S\in T_p\wtA\otimes T_p\wtA^t$.
For a deformation $\calA_{/B}$ of $\wtA$ with
$(B,\gotm_B)\in{\rm Ob}\,\widehat{\calC}$, let
$$
q(\calA_{/B};\cdot,\cdot)\colon
T_p\wtA\times T_p\wtA^t\mmap\widehat{\G}_m(B)=1+\gotm_B
$$
be the corresponding bilinear form.
When $\kk$ is not algebraically
closed, the group structure on $\calM\otimes\overline{\kk}$
descends to a group structure on $\calM$, for the details see
\cite[1.1.14]{Noot92}.

Let $\calN\subset\calM$ be a formal subgroup and
$\rho\colon X(\calM)\map X(\calN)$ the restriction map. The $\Z_p$-module
$N=\ker(\rho)$ is called the \emph{dual} of $\calN$. Via Serre-Tate
theory, $N\subseteq T_pA\otimes T_pA^t$. Then
$$
\calN\iisom\Hom\left(\frac{T_p\wtA\otimes T_p\wtA^t}N,\widehat{\G}_m\right)
$$
and
\begin{eqnarray*}
	\mbox{$\calN$ is a subtorus of $\calM$} & \Longleftrightarrow &
	\mbox{$X(\calN)\simeq X(\calM)/N$ is torsion-free}  \\
	 & \Longleftrightarrow & \mbox{$N$ is a direct summand of $T_p\wtA\otimes T_p\wtA^t$.}
\end{eqnarray*}

\noindent To simplify some of the next statements, we shall henceforth assume that $p>2$.

\begin{pro}\label{th:maplift}
	Let $\wtf\colon\wtA\map\wtB$ be a morphism of ordinary abelian
        varieties
	over $\kk$. The morphism $\wtf$ lifts to a morphism
	$f\colon\calA\map\calB$
	of deformations over $B$ if and only if
	$$
	q(\calA_{/B};P,\wtf^t(Q))=q(\calB_{/B};\wtf(P),Q)\quad
	\mbox{for all $P\in T_p\wtA$ and $Q\in T_p\wtB$}.
	$$
	In particular, if $(\wtA,\wt{\la})$ is principally polarized,
	the formal subscheme $\calM^{\rm pp}$ that classifies
	deformations of $\wtA$ with a
	lifting $\la$ of the principal polarization is a subtorus whose group
	of characters is
	$$
	X(\calM^{\rm pp})=\Sym^2(T_p\wtA).
	$$
\end{pro}

\pf The first part of the statement is \cite[2.1.4]{Katz81}. For the
second part, the principal polarization $\wt{\la}$ identifies
$T_p\wtA\simeq T_p\wtA^t$. For a deformation $\calA_{/B}$ let
$\pr{q}(\calA_{/B};P,\pr{P})=q(\calA_{/B};P,\wt{\la}(\pr{P}))$.
Then $\wt{\la}$ lifts to $\calA$ if and only if $\pr{q}$ is symmetric,
and the submodule of symmetric maps is a direct summand.
\qed

The last part of the proposition can be rephrased by saying that there
is a commutative diagram
$$
\begin{CD}
	T_p\wtA\otimes T_p\wtA^t @>\sim>> X(\calM) \\
	@V{\sopra{\hbox{polarization}}{\hbox{$+$ quotient}}}VV
	@VV\mbox{restriction}V \\
	\Sym^2(T_p\wtA) @>\sim>> X(\calM^{\rm pp})
\end{CD}
$$
Concretely, if $\{P_1,\ldots,P_g\}$ and $\{P_1^t,\ldots,P_g^t\}$ are
$\Z_p$-bases of $T_p\wtA$ and of $T_p\wtA^t$ respectively,
the $g^2$ elements
$q_{i,j}=q(\calA^u_{/\gotR^u};P_i,P_j^t)-1$ define an isomorphism
$\gotR^u\simeq W_{\kk}[[q_{i,j}]]$. If $\wtA$ is principally polarized
we may take $P_i=P_i^t$ under the identification
$T_p(\wtA)\simeq T_p(\wtA^t)$. Then $q_{i,j}=q_{j,i}$ on
$\calM^{\rm pp}=\Spf(\gotR^{\rm pp})$ and
$$
\gotR^{\rm pp}\simeq W_{\kk}[[q^{\rm pp}_{i,j}]],\qquad
\hbox{with $q^{\rm pp}_{i,j}={q_{i,j}}_{|\calM^{\rm pp}}$, $1\leq i\leq j\leq g$.}
$$
More generally, if $\calN=\Spf(\gotR_\calN)$ is a subtorus with
$n={\rm rk}_{\Z_p}(N)$, a $\Z_p$-basis
$\{S_1,\ldots,S_n\}$ of $N$ can be completed to a basis
$\{S_1,\ldots,S_n,S_{n+1},\ldots,S_{g^2}\}$ of
$T_p\wtA\otimes T_p\wtA^t$. If $q_i=q_{S_i}(\calA^u_{/\gotR^u})-1$ and
$q_i^\calN={q_i}_{|\calN}$, then $q^\calN_{1}=\ldots=q^\calN_{n}=0$ by construction
and $\gotR_\calN=W_{\kk}[[q^\calN_{n+1},\ldots,q^\calN_{g^2}]]$.

Since $\wtA$ is ordinary, there is a canonical isomorphism
$T_p\wtA^t\isom\Hom_B(\widehat{\calA},\widehat{\G}_m)$
for any deformation $\calA_{/B}$ of $\wtA$.
Composition with the pulling back the standard invariant form $dT/T$ on
$\widehat{\G}_m$ yields a functorial $\Z_p$-linear homomorphism
$\om\colon T_p\wtA^t\map\oom_{\calA/B}$ which is compatible with morphisms
of abelian schemes, in the sense that if the morphism
$f\colon\calA\map\calB$ lifts the morphism $\wtf\colon\wtA\map\wtB$ of
abelian varieties over $\kk$ then, \cite[lemma 3.5.1]{Katz81},
\begin{equation}
	f^*(\om(P^t))=\om(\wtf^t(P^t)),
	\qquad\mbox{for all $P^t\in T_p\wtB^t$.}
	\label{eq:omcomp}
\end{equation}
By functoriality, the maps $\om$ extend
to a well-defined $\Z_p$-linear homomorphism
$$
\om^u\colon T_p\wtA^t\map\oom_{\calA^u/\calM}.
$$
whose $\gotR^u$-linear extension
$T_p\wtA^t\otimes\calO_\calM\isom\oom_{\calA^u/\calM}$ is an isomorphism.
Thus, a choice of a $\Z_p$-basis $\{P_1^t,\ldots,P_g^t\}$ of
$T_p\wtA^t$ yields an identification
$\oom_{\calA^u/\calM}=\left(\bigoplus_{i=1}^{g}\gotR^u\om_i\right)^{\rm sh}$
where $\om_i=\om^u(P^t_i)$, $i=1,\ldots,g$ and the superscript 
$(\ )^{\rm sh}$ denotes the sheafified module.

Suppose that $\wtA$ is principally polarized. Let
$\calN\subset\calM^{\rm pp}$ be a subtorus with dual $N$
and let $\calA_\calN$ be the restriction over $\calN$ of the
universal deformation $\calA^u/\calM$.
Let $\{S_1,\ldots,S_{\frac{g(g+1)}{2}}\}$ be a $\Z_p$-basis of
$\Sym^2(T_p\wtA^t)=\Sym^2(T_p\wtA)$
such that $N=\bigoplus_{j=1}^n\Z_pS_j$. Let
$\om^{(2)}_i$ the pullback to $\calA_\calN$ of $\Sym^2(\om^u)(S_i)$, 
$q_i^{\rm pp}$ and $q_i^\calN$be the restriction of the local
parameters $q_{S_i}$ constructed above, $i=1,\ldots,g(g+1)/2$.
\begin{pro}\label{th:KSforsub}
  	$$
  	\KS_\calN(\om^{(2)}_i)=
  	\left\{
  	\begin{array}{ll}
  	 	0 & \mbox{for $i=1,\ldots,n$}  \\
  	 	d\log(q_i^{\rm pp}+1)_{|\calN} & \mbox{for $i=n+1,\ldots,g(g+1)/2$}
  	 \end{array}\right.
  	 $$
\end{pro}

\pf It is an immediate application of proposition
\ref{th:KSforpullb} to Katz's computations \cite[theorem 3.7.1]{Katz81}, because there are 
identifications
$\Sym^2(\oom_{\calA_\calN/\calN})=\left(\bigoplus_{i=1}^{g(g+1)/2}\gotR_\calN\om_i^{(2)}\right)^{\rm sh}$,
$\Om^1_{\calM^{\rm pp}/W}=\left(\bigoplus_{i=1}^{g(g+1)/2}\gotR^{\rm pp}dq_i^{\rm pp}\right)^{\rm sh}$,
$\Om^1_{\calN/W}=\left(\bigoplus_{i=n+1}^{g(g+1)/2}\gotR_\calN dq_i^\calN\right)^{\rm sh}$ 
and $dq_i^\calN=0$ for $i=1,...,n$. \qed

\paragraph{Elliptic curves and false elliptic curves.}
Let $(\tau,v,e)$ be a split $p$-ordinary test triple and denote 
$A$ the corresponding CM curve (if $D$ is split) or
principally polarized QM-abelian surface (if $D$ is non split) with $\wtA$ 
its reduction modulo $p$.
Of the two embeddings $K\hookrightarrow\Q_p$, fix the one that makes the $K$ action on 
$T_p\wtA\otimes\Q_p$ coincide with its natural $\Q_p$-vector space structure.

Let $\calM$ be the local moduli functor corresponding to $\wtA$.
If $D$ is split, everything has been already implicitely described:
$\calM^{\rm pp}=\calM$
is a $1$-dimensional torus and $\gotR^{\rm pp}=\gotR^u=W[[q-1]]$ with
$q=q_{P\otimes P}$ for a $\Z_p$-generator $P$ of $T_p\wtA$. Also,
if $\om_u=\om^u(P)$ then $\om_u^{\otimes 2}=\Sym^2(\om^u)(P\circ P)$
and $\KS(\om_u^{\otimes 2})=d\log(q)$.

If $D$ is non-split,
let $\calN=\calN_D$ be the subfunctor 
$$
\calN_D(B)=
\left\{
\sopra{\hbox{principally polarized deformations $\calA_{/B}$ of $\wtA$ 
with a lift of}}
{\hbox{the endomorphisms given by elements of the
maximal order $\calR_{1}$}}
\right\}
$$
The maximal order $\calR_{1}$ acts naturally on $T_p\wtA$ and since 
$e\in\calR_{1}\otimes\Z_p$ we can find a $\Z_p$-basis $\{P,Q\}$ of $T_p\wtA$
such that $eP=P$ and $eQ=0$.

\begin{pro}\label{th:NR}
	\begin{enumerate}
	 \item  $\calN=\Spf(\gotR_{\calN})$ is a $1$-dimensional subtorus  of $\calM^{\rm pp}$;
	 \item  $\gotR_{\calN}=W_{\kk}[[q-1]]$, where $q=q^\calN_{P\circ P}$;
	 \item  if $\om_u$ denotes the pullback of $\om^u(P)$ to $\calA_{\calN}$, then 
	            $\KS(\om_u^{\otimes 2})=d\log(q)$.
	\end{enumerate}
\end{pro}

\pf It follows from proposition \ref{th:maplift} that $\calN(B)$ is
identified with the set of the symmetric bilinear forms
$q\colon T_p\wtA\times T_p\wtA\map\widehat{\G}_m(B)$ such that
$q(P,\invol{r}Q)=q(rP,Q)$ for all $r\in\calR_{1}$.
This makes clear that $\calN$ is a subgroup, and that its
dual $N$ is the $\Z_p$-submodule generated by the elements
$$
    \left\{
    \begin{array}{l}
        P_1\otimes P_2-P_2\otimes P_1 \\
        rP_1\otimes P_2-P_1\otimes\invol{r}P_2 
    \end{array}\right.
    \quad\mbox{for all $P_1,P_2\in T_p(\wtA)$ and $r\in\calR_{1}$.}
$$
Choose $u$ in the decomposition \eqref{eq:Dsplit} for the subfield 
$F\subset D$ so that $u\in\calR_{1}$, $\invol u=-u$ and 
$\vass{\nu(u)}$ is minimal. In particular $u^{2}=-\nu(u)$ is a 
square-free integer. Pick a basis of $D$ in $\calR_{1}$ of 
the form $\{1,r,u,ru\}$ with $r^{2}\in\Z$. From our choice of test 
triple we can assume that $\Z[r]$ is an order in $F$ of conductor 
prime to $p$, in particular $p$ does not divide $r^{2}$.
The submodule $\calR=\Z\oplus\Z r\oplus\Z u\oplus\Z ru$ is actually an 
order of discriminant $-16r^{4}u^{4}$ such that $\invol\calR=\calR$.

Suppose that $p|u^{2}$ and let 
$y=\frac1p\left(a+br+cu+dru\right)$ with $a$, $b$, $c$ and $d\in\Z$
be an element in $\calR_{1}-\calR$ such that $\bar y=y+\calR$ 
generates the unique subgroup of order $p$ in $\calR/\calR_{1}$.
The conditions
$\tra(y)\in\Z$, $\nu(y)\in\Z$, and $\invol{\bar y}=\pm\bar y$
easily imply that the coefficients $a$, $b$, $c$ and $d$ are all 
divisible by $p$ and this is a contradiction. Thus 
$\calR\otimes\Z_{p}=\calR_{1}\otimes\Z_{p}$ and this reduces
the set of generators for $N$ to
$$
\begin{array}{ccc}
    P\otimes Q-Q\otimes P & rP\otimes Q-P\otimes rQ & uP\otimes 
    P+P\otimes uP  \\
    uP\otimes Q+P\otimes uQ & uQ\otimes Q+Q\otimes uQ & 
    ruP\otimes Q-P\otimes ruQ
\end{array}.
$$
From the relations $re=er$ and $ue=(1-e)u$ in $\calR_{1}\otimes\Z_p$ 
we get that the elements $r$ and $u$ act on the basis $\{P,Q\}$ as the 
matrices $\smallmat{\al}00{-\al}$ and $\smallmat 0{\be}{\ga}0$ 
respectively. Since $\al$, $\be$ and $\ga$ are $p$-units, finally $N$ 
turns out to be the $\Z_{p}$-module generated by
$$
P\otimes Q-Q\otimes P,\quad
\ga P\otimes P+\be Q\otimes Q,\quad
P\otimes Q+Q\otimes P
$$
and
$$
T_{p}\wtA\otimes T_{p}\wtA=N\oplus\Z_{p}(P\otimes P).
$$
This proves points 1 and 2, and the last part follows at once from 
proposition \ref{th:KSforsub}.\qed

\begin{rem}\label{rm:OK}
\rm The fact that $\calN$ is actually a subtorus can be reinterpreted, 
as in  \cite{Noot92}, in the more general context of Hodge and Tate
classes.
\end{rem}

\subsection{Power series expansion}
We shall now use the Serre-Tate theory to write a power 
series expansion around an ordinary CM point of a modular form 
$f\in M_{k,1}(\De,N)$ and compute the coefficients of this expansion in 
terms of the Maass operators studied in section \ref{ss:maass}. We 
assume that $N>3$.

Let $\hbox{\bf Sp}_D$ the full subcategory of the category of rings 
consisting of the rings $B$ such that $\calR_{1}\otimes B\simeq M_{2}(B)$. 
Note that if $(\tau,v,e)$ is a $p$-ordinary test triple, then $\calO_v\in{\rm Ob}\hbox{\bf Sp}_D$.
For a $\KS$-normalized section $s(z)$ in \eqref{eq:kanormsec} the assignment
\begin{equation}
	f(z)\mapsto f^\ast(z)=f(z)s(z)^{\otimes k}
	\label{eq:assignsplit}
\end{equation}
sets up an identification
$$
M_{k,1}(\De,N)\simeq H^0(X_1(\De,N),\calL^{\otimes k})
$$
defined up to a sign (the ambiguity obviously disappears for $k$ even).
The identification extends naturally to an identification of the bigger space 
$M_{k,\vep}^{\infty}(\De,N)$ of $C^\infty$-modular forms with the global sections of the associated $C^\infty$-bundle $\calL^{\otimes k}_\infty$.
This ``geometric'' interpretation of modular forms can be used to endow 
the space $M_{k,1}(\De,N)$ with a canonical $B$-structure for any subring 
$B\subset\C$ of definition for $\calL$ in $\hbox{\bf Sp}_D$.
In fact for \emph{any} ring $B$ in $\hbox{\bf Sp}_D$ such that
$\calL$ is defined over $B$ 
the space of modular forms defined over $B$ may be defined as
$$
M_{k,1}(\De,N;B)=H^0(\calX_1(\De,N)\otimes B,\calL^{\otimes k}).
$$
Remark \ref{re:LoverC}.2 shows that this $B$-structure does not depend 
on the choice of $\calL$, i.e. on the choice of idempotent $e$.
If $\pr{B}$ is a flat $B$-algebra, the identification
$M_{k,1}(\De,N;B)\otimes\pr{B}=M_{k,1}(\De,N;\pr{B})$
follows from the usual properties of flat base change.
By smoothness, if $1/N\De\in B\subset\C$  
then $M_{k,1}(\De,N;B)\otimes\C=M_{k,1}(\De,N;\C)=M_{k,1}(\De,N)$.
In fact the assignment \eqref{eq:assignsplit} is normalized so that 
$f\in M_{k,1}(N;B)$ if and only if its Fourier coefficients belong to $B$ 
($q$-expansion principle, e.÷g. \cite[Ch.~1]{Katz73} 
\cite[theorem 4.8]{Harris81}).

Let $\calA/\calX$ be either universal family \eqref{eq:univEC}
or \eqref{eq:univQMAV}. Let $x\in\calX(\calO_{(v)})$ be represented by a 
point $\tau\in\CM_{\De,K}$ in a split $p$-ordinary test triple 
$(\tau,v,e)$. Denote $A_x$ the fiber of $\calA/\calX$ over $x$ and 
$A_{\tau}$ the corresponding complex torus.  
We will implicitely identify the ring $W_{\overline{k(v)}}$ of Witt 
vectors for the algebraic closure of $k(v)$ with $\nr\calO_{v}$.
For each $n\geq0$, let
$J_{x,n}=\calO_{\calX,x}/\gotm_x^{n+1}$ and
$J_{x,\infty}=\limproj{n}J_{x,n}=\wh{\calO}_{\calX,x}$.
By smoothness, there is a non-canonical isomorphism 
$J_{x,\infty}\simeq\calO_{(v)}[[u]]$.
For $n\in\N\cup\{\infty\}$ the family $\calA/\calX$ restricts to
abelian schemes $\calA_{x,n}{}_{/J_{x,n}}$.
Tautologically $A_x=\calA_{x,0}$, and
$\calA_{x,n}=\calA_{x,\infty}\otimes J_{x,n}$ with respect to
the canonical quotient map $J_{x,\infty}\map J_{x,n}$.
Also, let $\nr{J}_{x,n}=J_{x,n}\wh{\otimes}\nr\calO_{v}$ and 
$\nr{\calA}_{x,n}=\calA_{x,n}\otimes\nr{J}_{x,n}$.

Let $\calM=\Spf(\calR)$ be either the full local moduli functor (in
the split case) or its subtorus described in proposition \ref{th:NR}
(in the non-split case) associated with the reduction
$\wtA_x=A_x\otimes{\overline{k}}_v$ with universal formal deformation
$\calA_x/\calM$. In either case $\calM\simeq\Hom(T,\wh{\G}_m)$ 
where $T$ is a free $\Z_p$-module of rank 1.
Since the rings $\nr{J}_{x,n}$ are pro-$p$-Artinian,
there are classifying maps
$$
\phi_{x,n}\colon\calR\mmap\nr{J}_{x,n},\qquad
\hbox{for all $n\in\N\cup\{\infty\}$}
$$
such that $\nr{\calA}_{x,n}=\calA_x\otimes_{\phi_{x,n}}\nr{J}_{x,n}$. Since
the abelian schemes $\calA_{x,n}$ are the restriction of the
universal (global) family, the map $\phi_{x,\infty}$ is an isomorphism.
We will use it to transport the Serre-Tate parameter
$q_S-1\in\calR$ and the formal sections $\om_u$ constructed in section
\ref{se:STtheory} out of a choice of a $Z_p$-generator $S$ of $T$
to the $p$-adic disc of points in $\calX$ that reduce
modulo $\gotp_v$ to the same geometric point in
$\calX\otimes{\overline{k}}_v$.  Also, we can pull back the parameter along 
the translation by $\inv x$ in $\calM$ to obtain a local parameter $u_x$ at $x$ (depending on $S$), namely
$$
	\nr{J}_{x,\infty}=\nr{\calO}_v[[u_x]],\qquad
	\hbox{with $u_x=\inv{q_S(x)}q_S-1$}.
$$
The complex uniformization of $A_x$ associated with the choice of
$\tau$ can be used to define transcendental periods. For any 
$\om_o\in H^0(A_{x}(\C),\calL(x))$ write
$$
\om_o=p(\om_o,\tau)s(\tau)
\qquad p(\om_o,\tau)\in\C,
$$
under the isomorphism $A_x(\C)\simeq A_\tau$.
For $f\in M_{k,1}(\De,N)$ define complex numbers
\begin{equation}
    c^{(r)}(f,x,\om_o)=
    \frac{\de_k^{(r)}(f)(\tau)}{p(\om_o,\tau)^{k+2r}}
    \qquad r=0,1,2,\ldots.
    \label{eq:defcrfx}
\end{equation}
The use of $x$ in the definition \eqref{eq:defcrfx} is justified by 
the following fact.

\begin{pro}
	Suppose that $f\in M_{k}(\Ga)$ for some Fuchsian group
	of the first kind $\calR_1^1\geq\Ga\geq\GaU(\De,N)$. 
	Then the numbers $c^{(r)}(f,x,\om_o)\in\C$ do not depend on the 
	choice of $\tau$ in its $\Ga$-orbit.
\end{pro}

\pf For any $\ga\in\Ga$, multiplication by $j(\ga,\tau)^{-1}$ induces 
an isomorphism of complex tori $A_\tau\isom A_{\ga\tau}$. Since $s$ is a global 
constant section of $\calL$ over $\gotH$, under the standard 
identifications of invariant forms, $s(\ga\tau)=s(\tau)j(\ga,\tau)^{-1}$. 
The assertion follows at once.\qed

The periods $p(\om_o,\tau)$ (and conseguently the numbers 
$c^{(r)}(f,x,\om_o)$) can be normalized by choosing $\om_o$ as in 
proposition \ref{teo:Thkalgebraic}. For such a choice,
defined up to a $v$-unit, set
$$
\Om_\infty=\Om_\infty(\tau)=p(\om_o,\tau),\qquad
c^{(r)}_v(f,x)=c^{(r)}(f,x,\om_o).
$$
Also, define the \emph{$p$-adic period} 
$\Om_p=\Om_p(x)\in\calO_v^{{\rm nr},\times}$ 
(again defined up to a $v$-unit) as
$$
\om_o=\Om_p\om_u(x).
$$

Let $f\in M_{k,1}(\De,N;\nr{\calO}_v)$. Over $\Spf(\nr{J}_{x,\infty})$
write $f^\ast=f_x\om_u^{\otimes k}$ and
$$
    \jet_x(f^\ast)=x^*\jet(f_x)\otimes\om_u(x)^{\otimes k}
    =\left(\sum_{n=0}^\infty\frac{b_n(f,x)}{n!}U_x^n\right)\,
    \om_u(x)^{\otimes k}
$$
where $f_x$ is expanded at $x$ in terms of the formal local parameter $U_x=\log(1+u_x)$.

\begin{thm}\label{thm:equality}
	Let $x\in\calX_{\De,N}(\calO_{(v)})$ be represented by a split
	$p$-ordinary test triple $(\tau,v,e)$ 
	and $f\in M_{k}(\De,N;\calO_{(v)})$. Then, for all $r\geq0$,
	$$
	\frac{b_r(f,x)}{\Om_p^{k+2r}}=c^{(r)}_v(f,x)\in\calO_{(v)}.
	$$
\end{thm}

\pf The case $r=0$ is clear, so let us assume that $r\geq1$.

We have $\nabla(f^\ast)=\nabla(f_x\om_u^{\otimes k})=
df_x\otimes\om_u^{\otimes k}+kf_x\om_u^{\otimes k-1}\nabla(\om_u)$.
Since $\nabla(\om^u(P))\in H^0(\calM,\calU)$ for each 
$P\in T_p(\wt{A})$, \cite[theorem 4.3.1]{Katz81}, the term containing 
$\nabla(\om_u)$ is killed by the projection $\PFr$.
Also, $dU_x=d\log(u_x+1)=d\log(q+1)$ doesn't depend on $x$ and we obtain
$\Th_{k,p}(f^\ast)=(df_x/dU_x)\om_u^{\otimes k+2}$. Iterating the 
latter computation $r$ times and evaluating the result at $x$ yields
$$
\Th_{k,p}^{(r)}(f^\ast)(x)=
\frac{d^rf_x}{d{U_x}^r}(x)\om_u^{k+2r}(x)=
\frac{b_r(f,x)}{\Om_p^{k+2r}}\om_o^{k+2r}.
$$
On the other hand, applying $r$ times the proposition
\ref{th:maassincoord} and evaluating at $\tau$ yields
$$
\Th_{k,\infty}^{(r)}(f^\ast)(x)=
\de_k^{(r)}(f)(\tau)s(\tau)^{\otimes k+2r}=
c^{(r)}(f,x,\om_o)\om_o^{k+2r}.
$$
The result follows from proposition \ref{teo:Thkalgebraic}. \qed

This result has a converse. For, we need the following preliminary discussion.
Let $\calD$ be any domain of characteristic $0$ and field of quotients $\calK$. The formal substitution $u=e^U-1=U+\frac1{2!}U^2+\frac1{3!}U^3+\cdots$ defines a bijection between the rings of formal power series $\calK[[u]]$ and $\calK[[U]]$. Under this bijection the ring $\calD[[u]]$ is identified with a subring of the ring of \textit{Hurwitz series}, namely power series of the form
$$
\sum_{n=0}^\infty\frac{\be_n}{n!}\,U^n,\qquad
\hbox{with $\be_n\in\calD$ for all $n=0,1,2,...$}.
$$
We say that a power series $\Phi(U)\in\calK[[U]]$ is $u$-integral if $\Phi(U)=F(e^U-1)$ for some 
$F(u)\in\calD[[u]]$.
Denote $c_{n,r}$ the coefficients defined by the polynomial identity
$$
n!\vvec Xn=X(X-1)\cdots(X-n+1)=\sum_{r=0}^nc_{n,r}X^r.
$$
The following possibly well known result is closely related to \cite[Th\'eor\`eme 13]{Serre73}.

\begin{thm}
    A Hurwitz series $\Phi(U)=\sum_{n=0}^\infty\frac{\be_n}{n!}\,Z^n$ is $u$-integral if and only if 
    $\frac1{d!}(c_{d,1}\be_1+c_{d,2}\be_2+\cdots+c_{d,d}\be_d)\in\calD$ for all $d=1,2,...$.
\end{thm}

\pf Let $F(u)=\Phi(\log(1+u))\in\calK[[u]]$. For any polynomial 
$P(X)=p_0+p_1X+\cdots+p_dX^d\in\calK[X]$ of degree $d$, an immediate chain rule computation yields
$$
\left.P\left((u+1)\frac d{du}\right)F(u)\right|_{u=0}=p_0\be_0+p_1\be_1+...+p_d\be_d.
$$
Also 
$\left.P\left((u+1)\frac d{du}\right)F(u)\right|_{u=0}=\left.P\left((u+1)\frac d{du}\right)F_d(u)\right|_{u=0}$
where $F_d(u)$ is the degree $d$ truncation of $F(u)$, i.÷e. $F(u)=F_d(u)+U^{d+1}H(u)$. On the space of polynomials of degree $\leq d$ the substituition $u=v-1$ is defined over $\calD$ and
$\left.P\left((u+1)\frac d{du}\right)F_d(u)\right|_{u=0}=
\left.P\left(v\frac d{dv}\right)F_d(v-1)\right|_{v=1}$.
Since $\left.P\left(v\frac d{dv}\right)v^k\right|_{v=1}=P(k)v^k$, the argument shows that if $\Phi(U)$ is 
$u$-integral, then the expression $p_0\be_0+p_1\be_1+...+p_d\be_d$ is a $\calD$-linear combination of the values $P(0)$, $P(1), ..., P(d)$.

On the other hand, the argument also shows that
$$
\frac1{d!}(c_{d,1}\be_1+c_{d,2}\be_2+\cdots+c_{d,d}\be_d)=
\left.\vvec{v\,d/dv}{d}F_d(u-1)\right|_{v=1}
$$
is the coefficient of $U^d$ in $F$. 

Therefore, we obtain that $\Phi(U)$ is $u$-integral if and only 
$p_0\be_0+p_1\be_1+...+p_d\be_d\in\calD$ for every polynomial 
$P(X)=p_0+p_1X+\cdots+p_dX^d\in\calK[X]$ such that $P(0)$, 
$P(1), ..., P(d)\in\calD$. We conclude observing that a degree $d$ polynomial $P(X)\in\calK[X]$ such that $P(0)$, $P(1), ..., P(d)\in\calD$ is necessarily \textit{numeric}, i.÷e. $P(\N)\subset\calD$ and that the 
$\calD$-module of numeric polynomials is free, generated by the binomial coefficients. \qed

Note that when $\calD$ is a ring of algebraic integers, or one of its non-archimedean completions, the conditions of the theorem can be readily rephrased in terms of congruences, known as \textit{Kummer-Serre congruences}.

Denote $L^{v, {\rm sc}}$ the compositum of all finite extensions 
$L\subseteq F$ such that $v$ splits completely in $F$ and let $\calO^{\rm sc}_{(v)}$ 
be the integral closure of $\calO_{(v)}$ in $L^{v, {\rm sc}}$.

\begin{thm}[Expansion principle]\label{thm:expanprinc}
    Let $f\in M_{k,1}(\De,N)$ and $x\in\calX(\De,N)(\calO_{(v)})$ be represented
    by a split $p$-ordinary test triple $(\tau,v,e)$ such that the numbers 
    $c^{(r)}_v(f,x)\in\calO_{(v)}$ for all $r\geq0$ and 
    the $p$-adic numbers $\Om_p^{2r}c^{(r)}_v(f,x)$ satisfy the Kummer-Serre
    congruences. Then $f$ is defined over $\calO^{\rm sc}_{(v)}$.
\end{thm}

\pf Choose a field embedding $\imath\colon\C\map\C_{p}$ to 
view $f\in M_{k,1}(\De,N;\C_{p})$. For all $r\geq0$ set $c_r=c^{(r)}_v(f,x)$ and 
$\be_{r}=c_r\Om_{p}^{\ep_Dk+2r}\in\C_{p}$.
Unwinding the computations that led to the equality in theorem 
\ref{thm:equality} shows that 
$\jet_x(f^\ast)=\left(\sum_{r\geq0}\frac{\be_r}{r!}U_x^r\right)
\om_u(x)^{\otimes k}\in(\nr{J}_{x,\infty}\otimes\C_p)\om_u(x)^{\otimes k}$.
We claim that $\jet_x(f^\ast)$ is defined over $\calO_v$. 
Write
$$
\jet_x(f^\ast)=
\left(\sum_{r\geq0}\frac{\be_r}{r!}U_x^r\right)\om_u(x)^{\otimes k}=
\left(\sum_{r\geq0}\frac{c_r}{r!}(\Om_p^2U_x)^r\right)\om_o^{\otimes k}.
$$
Since the $\be_r$ are $v$-integral and satisfy the Kummer-Serre 
congruences, the first equality shows that $\jet_x(f^\ast)$ is $v$-integral. 
Since the formal substitution $u=e^U-1$ preserves the field of 
definition, the claim follows from the second equality if we check 
that the formal local parameter $\Om_p^2U_x$ is defined over $L_v$.
The group $\Aut(\C_p/L_v)$ acts on the section $\om_u$ via the action 
of its quotient $\Gal(\nr{L}_v/L_v)\simeq\Gal(\overline{k}_v/k_v)$ 
on $T_p(\wt{A}_x)$ which is scalar because $\wt{A}_x$ is either an
elliptic curve or isogenous to a product of elliptic curves. 
Thus, the section $\Om_p\om_u$, whose restriction at $x$ is 
defined over $L_v$, is itself defined over $L_v$. Therefore
$\Om_p^2\,dU_x=\KS(\Om_p^2\om_u^{\otimes 2})$ is defined over 
$L_v$, and so is $\Om_p^2U_x$ because it is a priori defined over 
$\nr{L}_v$ and its value at the point $x$, defined over 
$L_v$, is $0$.

We can now use the very same arguments of Katz's proof \cite{Katz73} 
of the $q$-expansion principle 
to conclude that the section $f^\ast$ is defined over $\calO_v$. 
For, observe that the $q$-expansion of a modular form $f$ at the cusp $s$ 
multiplied by the right power of the canonical Tate form is 
$\jet_s(f^\ast)$. The specific nature of a cusp in the modular curve 
plays no role in Katz's proof, which works as well when the former 
are replaced by any point in a smooth curve.

Since $f^\ast$ is defined over $\C$ and over $\calO_v$, the modular 
form $f$ is defined over the integral closure of $\calO_{(v)}$ in the largest 
subfield $F\subset\C$ such that $\imath(F)\subseteq L_v$. The 
assertion follows from the arbitrariness of the choice of $\imath$, 
since $L^{\rm sc}_{(v)}$ can be characterized as the largest subfield 
of $\C$ whose image under all the embeddings $\C\map\C_p$ is contained 
in $L_v$. \qed


\section{$p$-adic interpolation}

\subsection{$p$-adic $K^{\times}$-modular forms}\label{se:padicforms}

A \textit{weight} for the quadratic imaginary field $K$ is a formal linear combination 
$\uw=w_1\si_1+w_2\si_2\in\Z[I_K]$, which
will be also written $\uw=(w_1,w_2)$. Following our conventions, write 
$z^\uw=z^{w_1}\bar{z}^{w_2}$ for all $z\in\C$. Also, let 
$\bar{\uw}=(w_{2},w_{1})$ and $\vass\uw=w_1+w_2$, so that 
$\uw+\bar{\uw}=\vass\uw\underline{1}$ with $\underline{1}=(1,1)$.

\begin{dfn}[Hida \cite{Hida86}]
	Let $E\supseteq K$ be a subfield of $\C$. The space
	$\wtS_{\uw}(\gotn;E)$ of $K^{\times}$-modular form of weight $\uw$ and level $\gotn$ with values 
	in $E$ is the space of functions $\tif\colon\calI_{\gotn}\map E$ such that
	$$
	\tif((\la)I)=\la^{\uw}\tif(I)
	$$
	for all $\la\in K^{\times}_{\gotn}$.
\end{dfn}

\noindent A remarkable subset of $\wtS_{\uw}(\gotn)=\wtS_{\uw}(\gotn;\C)$ is the set of algebraic 
Hecke characters of type {$\mathrm{A}_{0}$}, 
$$
\wtXi_{\uw}(\gotn)=\wtS_{\uw}(\gotn)\cap\Hom(\calI_{\gotn},\C^{\times}).
$$
A well-known property noted by Weil \cite{Weil55} 
is that for every $\tixi\in\wtXi_{\uw}(\gotn)$ there exists a 
number field $E_{\tixi}$ such that $\tixi\in\wtS_{\uw}(\gotn;E_{\tixi})$.
A classical construction identifies the space $\wtS_{\uw}(\gotn)$ 
with the space $S_{\uw}(\gotn)$ of functions $f\colon\ideles\map\C^{\times}$ such that
\begin{equation}
f(s\la zu)=z^{-\uw}f(s)
\quad
\hbox{for all $\la\in K^{\times}$, $z\in\C^{\times}$ and $u\in U_{\gotn}$.}
\label{eq:ftilde}
\end{equation}
If $\tif\leftrightarrow{f}$ under this identification, then 
\begin{equation}
	\tif(I)={f}(s)\quad
	\hbox{whenever $I=[s]$ and $s_{v}=1$ for $v=\infty$ and $v|\gotn$.}
	\label{eq:weilrel}
\end{equation}
This relation can be used to recognize $\wtS_{\uw}(\gotn;E)$ in ${S}_{\uw}(\gotn)$.
Since $U_{c\calO_K}<\whO_{K,c}^\times$, the functions 
${f}\colon\ideles\map\C^{\times}$ satisfying the relation in \eqref{eq:ftilde}
for all $\la\in K^{\times}$,  $z\in\C^{\times}$ and 
$u\in\whO_{K,c}^{\times}$ form a linear subspace 
${S}_{\uw}(\calO_{K,c})\subset{S}_{\uw}(c\calO_{K})$. 
The subspace ${S}_{\uw}(\calO_{K,c})$ 
includes the Hecke characters trivial on $\whO_{K,c}^\times$, namely 
${\Xi}_\uw(\calO_{K,c})={S}_\uw(\calO_{K,c})\cap
\Hom(\ideles/K^\times\whO_{K,c}^\times,\C^\times)$.
Denote 
$$
\CC_{\gotn}=\ideles/K^{\times}\C^{\times}U_{\gotn}\simeq 
\calI_{\gotn}/P_{\gotn},\qquad
\CC^{\sharp}_{c}=\ideles/K^\times\C^\times\whO_{K,c}^\times.
$$ 
and let $h_{\gotn}=\vass{\CC_{\gotn}}$ and 
$h^{\sharp}_{c}=\vass{\CC^{\sharp}_{c}}$.  Clearly $h^\sharp_c|h_c$.

\begin{lem}\label{th:charchar}
	\begin{enumerate}
	  \item  ${\Xi}_\uw(\calO_{K})\neq\emptyset$ if and only if $(\calO_{K}^{\times})^{\uw}=1$.
	  \item ${\Xi}_\uw(\calO_{K,2})={\Xi}_\uw(2\calO_K)$ and they are non-trivial if and only if 
	            $\vass\uw$ is even.
	  \item If $c>2$ then ${\Xi}_\uw(c\calO_{K})\neq\emptyset$  and 
	            ${\Xi}_\uw(\calO_{K,c})\neq\emptyset$ if and only if $\vass\uw$ is even.
	  \item ${\Xi}_\uw(\calO_{K,c})$ and ${\Xi}_\uw(c\calO_{K})$ are 
	bases for ${S}_\uw(\calO_{K,c})$ and ${S}_\uw(c\calO_{K})$ respectively.
	\end{enumerate}
\end{lem}

\pf Let $U<U_1$ and $\CC_U=\ideles/K^{\times}\C^{\times}U$. Then there is a short exact sequence 
$$
1\mmap\frac{\C^{\times}}{H_U}\mmap
\frac{\ideles}{K^\times U}\mmap\CC_U\mmap1
$$
where $H_U=\C^\times\cap K^{\times}U$. The first three points follow from the observation that
$$
H_{U_c}=
\begin{cases}
   \calO_K^\times   & \text{if $c=1$}, \\
   \{\pm1\}  & \text{if $c=2$}, \\
   \{1\}    & \text{if $c>2$},
\end{cases}
\qquad
H_{\whO_{K,c}^\times}=
\begin{cases}
   \calO_K^\times   & \text{if $c=1$}, \\
   \{\pm1\}    & \text{if $c\geq2$}.
\end{cases}
$$
For the last part, observe that multiplication by any 
$\xi\in{\Xi}_\uw(c\calO_K)$ defines, for every weight $\uw^\prime$, 
an isomorphism ${S}_{\uw^\prime}(c\calO_K)\isom{S}_{\uw+\uw^\prime}(c\calO_K)$ 
which identifies the respective sets of Hecke characters.
When $\xi\in\Xi_\uw(\calO_{K,c})\neq\emptyset$, the isomorphism restricts to an isomorphism of the subspaces  ${S}_{\uw^\prime}(\calO_{K,c})\isom{S}_{\uw+\uw^\prime}(\calO_{K,c})$.
So, we are reduced to check the assertion in the case of the null weight 
$\underline{0}=(0,0)$, which is clear because 
${S}_{\underline{0}}(c\calO_{K})$ and ${\Xi}_{\underline{0}}(c\calO_{K})$
(respectively ${S}_{\underline{0}}(\calO_{K,c})$ and ${\Xi}_{\underline{0}}(\calO_{K,c})$)
are  the set of functions on the finite abelian group
$\CC_c$ (respectively $\CC^{\sharp}_{c}$) and its Pontryagin dual.\qed

If $\gotm|\gotn$ the inclusion $\calI_{\gotn}<\calI_{\gotm}$ defines a 
natural restriction map 
\begin{equation}
	\wtS_{\uw}(\gotm)\map\wtS_{\uw}(\gotn).
	\label{eq:restr}
\end{equation}

\begin{lem}\label{th:injec}
	The restriction maps \eqref{eq:restr} are injective.
\end{lem}

\pf We can assume that $\gotn=\gotm\gotp$ with $\gotp$ prime and 
$(\gotp,\gotm)=1$. Let $\tif\in\wtS_{\uw}(\gotm)$ and suppose that $\tif(I)=0$ 
for all ideals $I\in\calI_{\gotn}$. Let $\la\in K^{\times}_{\gotm}$ such 
that $\la\calO_{\gotp}=\gotp\calO_{\gotp}$. 
Then $\gotp[\la^{-1}]\in\calI_{\gotn}$ and 
$0=\tif(\gotp[\la^{-1}])=\la^{-\uw}\tif(\gotp)$, i.e. $\tif(\gotp)=0$ proving 
that $\tif=0$ identically. \qed

For ${f}\in{S}_{\uw}(\gotn)$ and ${g}\in{S}_{\uw^{\prime}}(\gotn)$ let
\begin{equation}
	\scal{f}{g}=
	\left\{
	\begin{array}{ll}
		h_{\gotn}^{-1}\sum_{\si\in\CC_{\gotn}}{f}(s_{\si}){g}(s_{\si})
		=h_{\gotn}^{-1}\sum_{\si\in\CC_{\gotn}}\tif(I_{\si})\tilde{g}(I_{\si})
		  & \hbox{if $\uw^{\prime}=-\uw$}  \\
		  & \\
		0 & \hbox{if $\uw^{\prime}\neq-\uw$}
	\end{array}
	\right.,
	\label{eq:padicpair}
\end{equation}
where $\{s_{\si}\}$ and $\{I_{\si}\}$ are full set of representatives of 
$\CC_{\gotn}$ in $\ideles$ and in $\calI_{\gotn}$ respectively.
The bilinear form $\scal{\cdot}{\cdot}$ extends 
by linearity to a pairing on
${S}(\gotn)=\bigoplus_{\uw\in\Z[I_K]}{S}_{\uw}(\gotn)$,
or on the corresponding space
$\wtS(\gotn)=\bigoplus_{\uw\in\Z[I_K]}\wtS_{\uw}(\gotn)$,
compatible with the restriction maps \eqref{eq:restr}. Note that for 
Hecke characters ${\xi}\in{\Xi}_{\uw}(\gotn)$ and 
${\xi}^{\prime}\in{\Xi}_{-\uw}(\gotn)$ one has the orthogonality relation
$$
\scal{\xi}{\xi^\prime}=
\left\{
\begin{array}{ll}
	1 & \hbox{if ${\xi}^{\prime}={\xi}^{-1}$} \\
	0 & \hbox{otherwise}
\end{array}
\right..
$$
\begin{rem}
\rm It follows at once from the definition \eqref{eq:padicpair}
that the pairing $\scal{\cdot}{\cdot}$ takes values in $E$ on $E$-valued forms. 
\end{rem}

Let $p$ be a prime number and let $F$ be a $p$-adic local field with 
ring of integers $\calO_{F}$. 
Following \cite{Hida86, Tilo96}, the space of $p$-adic 
$K^{\times}$-modular forms of level $\gotn$ with coefficients in $F$ is the 
space $\gotS(\gotn;F)=\calC^{0}(\gotC_{\gotn},F)$ of $F$-valued continuous 
functions on 
$\gotC_{\gotn}=\limproj{}{}_{r\geq0}\CC_{\gotn p^{r}}$.
It is a $p$-adic Banach space under the sup norm 
$\vvass{\phi}=\sup_{x\in\gotC_{\gotn}}\vass{\phi(x)}$
and we denote 
$\gotS(\gotn;\calO_{F})$ its unit ball. 
Assume that $E$ is a subfield of $F$ (e.g. $F$ is the completion of 
$E$ at a prime dividing $p$) and write 
$\wtS_{\uw}(\gotn;F)=\wtS_{\uw}(\gotn;E)\otimes F$ and 
$\wtS_{\uw}(\calO_{K,c};F)=\wtS_{\uw}(\calO_{K,c};E)\otimes F$.

\begin{pro}[\cite{Tilo96}]\label{th:padicembed}
	  For every ideal $\gotm|\gotn$ and for every ideal $\gotq$ with support included in the set of primes dividing $p$ there is a natural embedding
			$$
			\wtS(\gotm\gotq;F)=
			\bigoplus_{\uw\in\Z[I_K]}\wtS_{\uw}(\gotm\gotq;F)
			\hookrightarrow\gotS(\gotn;F).
			$$
\end{pro}

\pf We may use Lemma \ref{th:injec} to assume 
that $\gotm\gotq=\gotn p^{a}$ for some $a\geq1$. 
Since $\bigcap_{r\geq0}P_{\gotn p^{r}}=\{1\}$ the 
group $I_{\gotn p}$ embeds as a dense subset in $\gotC_{\gotn}$. 
The restriction of $\tif\in\wtS_{\uw}(\gotn p^{a};F)$
to a coset $I\cdot P_{\gotn p^{a}}$ is the function 
$I(\la)\mapsto\tif(I)\la^{\uw}$. Since $\uw\in\Z[I_{K}]$ the 
character $\la^{\uw}$ is continuous for the $p$-adic topology on 
$K^{\times}$ and so extends to a character $\chi_{\uw}$ of
$(K\otimes\Q_{p})^{\times}$. Therefore
$\tif$ extends locally to cosets of 
$1+p^{a}(R_{K}\otimes\Z_{p})$ and globally to the whole of 
$\gotC_{\gotn}$. The injectivity of the direct sum space 
$S(\gotn p^{a};F)$ follows 
from the linear independence of characters. \qed

We shall denote $\wh{f}$ the $p$-adic modular form associated to the 
$K^{\times}$-modular form $\tif$. If $\tif=\tilde{\xi}$ is an Hecke character, the 
$p$-adic form $\wh{\xi}$ is again a character which is sometimes 
called the \textit{$p$-adic avatar} of $\tilde{\xi}$ (or of ${\xi}$).
The density of $I_{\gotn p}$ in $\gotC_{\gotn}$ implies also that the image of 
$\wtS_{\uw}(\gotn p^{a};F)$ in $\gotS(\gotn;F)$ is characterized by the 
functional relations $\tif(\la s)=\la^{\uw}\tif(s)$ for all
$\la\equiv1\bmod\gotn p^{a}$. Thus, the association $\tif\mapsto\wh{f}$ 
identifies $\wtS_{\uw}(\gotn p^{a};F)$ with the closed linear subspace
\begin{equation}
\gotS_{\uw,a}(\gotn;F)=\left\{\sopra
{\mbox{$\phi\in\gotS(\gotn;F)$ such that $\phi(sx)=\phi(s)\chi_{\uw}(x)$}}
{\mbox{ for all $x\in1+p^{a}(\calO_{K}\otimes\Z_{p})$}}
\right\}
\label{eq.clospace}
\end{equation}
(when $a=0$ the domain for $x$ is $(\calO_{K}\otimes\Z_{p})^{\times}$).
Let $\overline{S}(\gotn p^{a};F)=\wh{\bigoplus}_{\uw}
\gotS_{\uw,a}(\gotn;F)$ be the closure of
$\wtS(\gotn p^a;F)$ in $\gotS(\gotn;F)$.
Since $\gotS_{\uw,a}(\gotn;F)$ is closed the projection onto the $w$-th summand extends to a projection
$\pi_{\uw,a}\colon\overline{S}(\gotn p^{a};F)\map\gotS_{\uw,a}(\gotn;F)$.
Define a pairing
$$
	\scalq{\cdot}{\cdot}\colon
	\overline{S}(\gotn p^{a};F)\times\overline{S}(\gotn p^{a};F)\mmap F
$$
as the composition
$$
\overline{S}(\gotn p^{a};F)\times\overline{S}(\gotn p^{a};F)
\stackrel{m}{\mmap}\overline{S}(\gotn p^{a};F)
\stackrel{\pi_{\underline{0},a}}{\mmap}
\gotS_{\underline{0},a}(\gotn;F)
\stackrel{\mu_{\mathrm{H}}}{\mmap}F
$$
where $m$ is multiplication and $\mu_{\mathrm{H}}$ is the Haar 
distribution which is bounded on the space $\gotS_{\underline{0},a}(\gotn;F)$. 

\begin{pro}\label{th:spiden}
     The pairing $\scal{\cdot}{\cdot}$ extends to a continuous pairing on
     $\overline{S}(\gotm\gotq;F)$ which coincides with
     $\scalq{\cdot}{\cdot}$. 	
\end{pro}

\pf The pairing $\scalq{\cdot}{\cdot}$ is continuous as composition 
of continuous mappings. Thus it is enough to check the identity 
$\scal{\tif}{\tilde{g}}=[{\wh{f}},{\wh{g}}]$ for 
$\tif\in\wtS_{\uw}(\gotn p^{a};F)$ and $\tilde{g}\in\wtS_{\uw^{\prime}}(\gotn p^{a};F)$. 
It follows from the definition \eqref{eq:padicpair} and the 
density of $\calI_{\gotn p}$ in $\gotC_{\gotn}$, since on $\calI_{\gotn p}$
the restrictions of $\tif$ and $\wh{f}$ and of $\tilde{g}$ and $\wh{g}$ coincide.
\qed

Recall that a $p$-adic distribution on $\Z_p$ with values in the $p$-adic Banach space $W$ over $F$ is a linear operator $\calC^0(\Z_{p},F)\map W$. Given two $p$-adic distributions on $\Z_{p}$ with values in 
$\overline{S}(\gotn p^{a};F)$
we construct a new distribution $\mu_{\scalq{\mu_1}{\mu_2}}$ with values in $F$ as the composition
$$
\calC^0(\Z_{p},F)\stackrel{\mu_1\ast\mu_2}{\mmap}
\overline{S}(\gotn p^{a};F)
\stackrel{\pi_{\underline{0},a}}{\mmap}
\gotS_{\underline{0},a}(\gotn;F)
\stackrel{\mu_{\mathrm{H}}}{\mmap}F,
$$
where $\mu_1\ast\mu_2$ is the convolution product of $\mu_1$ and $\mu_2$. If $\mu_1$ and $\mu_2$ are measures (bounded distributions) $\mu_{\scalq{\mu_1}{\mu_2}}$ is not a measure in general since the map 
$\pi_{\underline{0},a}$ is not bounded.
Denote 
$m_k(\mu)=\int_{\Z_p}x^k\,d\mu(x)$, $k\geq0$ the $k$-th moment of the distribution $\mu$.

\begin{lem}\label{th:measuremix}
	Let $M\in\N\cup\{\infty\}$ and suppose that there exist
	pairwise distinct weights $\{\uw_{k}\}$ for $0\leq k<M$ 
	such that 
	$m_{k}(\mu_{1})\in\wtS_{\uw_{k}}(\gotn p^{a};F)$ and
	$m_{k}(\mu_{2})\in\wt S_{-\uw_{k}}(\gotn p^{a};F)$ for all 
	$0\leq k<M$. Then
	$$
	m_{k}(\mu_{\scalq{\mu_1}{\mu_2}})=
	\left\{
	\begin{array}{ll}
		0 & \mbox{if $0\leq k<M$ is odd,}  \\
		\binom{2l}{l}\scalq{m_{l}(\mu_{1})}{m_{l}(\mu_{2})} & 
		\mbox{if $0\leq k=2l<M$ is even.}
	\end{array}
	\right.
	$$
	If $M=\infty$ the latter formulae characterize the distribution
	$\mu$ completely.
\end{lem}

\pf By direct computation
$m_k(\mu_{\scalq{\mu_1}{\mu_2}})=\mu_H\circ\pi_{\underline{0},a}
\left(\iint_{\Z_{p}^2}(x+y)^{k}\,d\mu_1(x) d\mu_2(y)\right)=
\sum_{i=0}^{k}\vvec{k}{i}\mu_H\circ\pi_{\underline{0},a}\left(m_i(\mu_1)m_{k-i}(\mu_2)\right)=
\sum_{i=0}^{k}\vvec{k}{i}\scalq{m_i(\mu_1)}{m_{k-i}(\mu_2)}.
$
The formula follows at once from the orthogonality relations in 
\eqref{eq:padicpair} since $\uw_{i}=\uw_{k-i}$ only if $k=2l$ is 
even and $i=l$. The final assertion is also clear.\qed

Let $\mu$ be a $p$-adic distribution on $\Z_p$ with values in a $p$-adic space $\calS$ of continuous $F$-valued functions on a profinite space $T$.  For   every $t\in T$, evaluation at $t$ defines an $F$-valued distribution $\mu(t)$ on $\Z_p$, $\mu(t)(\phi)=\mu(\phi)(t)$.
Conversely, a family $\{\mu_t\}_{t\in T}$ of $F$-valued distributions such that the function
$\mu(\phi)(t)=\mu_t(\phi)$ is in $\calS$ for all $\phi\in\calC^0(\Z_{p},F)$ defines a $p$-adic distribution 
$\mu$ on $\Z_p$ with values in $\calS$ and $\mu(t)=\mu_t$ for all $t\in T$,  which is obviously unique for this property.

\begin{lem}\label{le:UnBoundPr}
    Let $T$ be a profinite space, $\cal S$ a $p$-adic space of continuous $F$-valued functions and 
    $\mu$ a $p$-adic distribution on $\Z_p$ with values in $\cal S$. Then $\mu$ is a $p$-adic measure 
    if and only if $\mu(t)$ is a $p$-adic measure for all $t\in T$.
\end{lem}

\pf If $\mu$ is bounded, the distributions $\mu(t)$ are obviously bounded. 

Suppose that $\mu(t)$ is bounded for all $t\in T$. Let $\{\phi_k\}$ $k=0,1,2,...$ be functions in 
$\calC^0(\Z_{p},F)$ with $\vvass{\phi_k}=1$ and let $\varphi_k=\mu(\phi_k)$.
If $\vvass{\varphi_k(t)}=p^{r_k}$ choose $t_k\in T$ such that 
$\vass{\varphi_k(t_k)}_p=p^{r_k}$.
If the set of values $\left\{p^{r_k}\right\}$ is not bounded we may assume without loss of generality that $r_1<r_2<r_3<\cdots$ and since each $\mu(t_k)$ is bounded also that $\left\{t_k\right\}$ is an infinite set.
By compactness of $T$, there exists $\bar{t}\in T$, ${\bar t}\neq t_k$ for all $k$, such that every
neighborhood of $\bar t$ meets $\left\{t_k\right\}$. This contradicts the boundedness of 
$\mu({\bar t})$ since $\vass{\varphi(t)}_p$ is locally constant for all $\varphi\in\cal S$.

In particular, the sequence $\mu\left(\binom xk\right)$ is bounded and $\mu$ is a measure. \qed

As an application, let $\tilde{\chi}\in\wtXi_{\uw}(\gotn p^a)$ and 
$\tixi\in\wtXi_{\uw^{\prime}}(\gotn p^a)$ be Hecke characters 
taking values respectively in $\calO_{F}^{\times}$ and  
$\calO_{F^\prime}^{\times}$ where $\Q_p\subseteq F^\prime$ is a totally ramified subextension of $F$. 
For every $x\in\gotC_{\gotn}$ the series 
$$
\sum_{k=0}^{\infty}\frac{1}{k!}\wh{\chi}{\wh{\xi}}^{k}(x)Z^{k}
=\wh{\chi}(x)(1+T)^{\wh{\xi}(x)},\qquad
T=e^{Z}-1=Z+\frac12Z^{2}+\frac1{3!}Z^{3}+\cdots,
$$ 
has integral coefficients in the variable $T$ (the corresponding measure is 
$\wh{\chi}(x)\partial_{\wh{\xi}(x)}$, 
where $\partial_t$ denotes the Dirac measure concentrated at $t$). 
Thus, there exists a unique measure $\mu_{\chi,\xi}$ on $\Z_{p}$ with values in 
$\overline{S}(\gotn;\calO_{F})$ such that 
$m_{k}(\mu_{\chi,\xi})=\wh{\chi}\wh{\xi}^{k}$. When 
$\uw^{\prime}\neq\underline{0}$ the moments' weights are pairwise distinct.

\subsection{Expansions as distributions}

Let $\jmath\colon K\hmap D$ a normalized embedding of conductor 
$c=c_{\tau,N}$ with corresponding $\tau\in\CM_{\De,K}$ and
$x\in\CM(\De,N;\calO_{K,c})$. Let $y=\mathrm{Im}(\tau)$.
The embedding $\jmath$ defines by scalar extension a diagram
\begin{equation}
	\begin{array}{ccccccc}
		 & & K_\A^\times/K^\times\R^\times & \mmap & 
		D^\times\bs D_\A^\times/Z_\infty  \\
		& & \downarrow &  & \downarrow  \\
		& & K_\A^\times/K^\times\R^\times\wh{\calO}_{K,c}^\times & \mmap & 
		D^\times\bs D_\A^\times/Z_\infty\wh{\calR}_N^\times  \\
		& & \downarrow &  & \downarrow  \\
		\CC^{\sharp}_{c} & \simeq & K_\A^\times/K^\times\C^\times\wh{\calO}_{K,c}^\times & \mmap  & 
		D^\times\bs D_\A^\times/\jmath(\C^\times)\wh{\calR}_N^\times & \simeq & 
		\GaZ(\De,N)\bs\gotH
	\end{array}
	\label{eq:maindiag}
\end{equation}
where the vertical maps are the natural quotient maps and $Z_\infty$ is the center of $D_\infty^\times$.
Under the decomposition 
\begin{equation}
	D_\A^\times=D_\Q^\times\GL_2^+(\R)\wh{\calR}_N^\times
	\label{eq:decomp}
\end{equation}
the idele $d=d_og_\infty u$ corresponds to the point represented by $g_\infty\tau$.
Classfield theory provides an identification $\CC^{\sharp}_{c}\simeq\Gal(H_c/K)$
where $H_c$ is the ray classfield of conductor $c$. It is also 
well-known that the points in the image of the bottom map in 
\eqref{eq:maindiag} are defined over $H_c$, so that
$\Gal(H_c/K)$ acts naturally on them, and that the two actions are 
compatible (Shimura reciprocity law, \cite{ShiRed}). In particular, 
if $s_{\si}\in\ideles$ represents $\si\in\Gal(H_c/K)$, then $s_{\si}$ 
maps to $x^{\si}$ and $A_{x^\si}=A_x^{(\inv\si)}$. Write 
$A_x(\C)=A_{\tau}=\C^{\ep}/\La_\tau$ 
with $\La_\tau=\La\vvec{\tau}{1}$ where $\ep=1$ and
$\La=\Z^2\subset\C$ if $D$ is split and $\ep=2$ and
$\La=\Phi_{\infty}(\calR_{1})\subset\C^2$ if $D$ is non-split.
The theory of complex multiplication implies that 
$A_{x^\si}(\C)\simeq\C^{\ep}/s_{\si}\La_\tau$ where
$s_{\si}\La_\tau=\La d_{\si}^{-1}\vvec{\tau}{1}$ if 
$s_{\si}^{-1}=d_{\si}g_{\si}u_{\si}$ under \eqref{eq:decomp}.
For a fixed  prime $p$ one can choose representatives
$\left\{s_{\si}\right\}\subset\ideles$ normalized as follows:
$$
\left\{\begin{array}{l}
    s_{\si,\infty}=1,  \\
    s_{\si,v}\mbox{ is $v$-integral at all finite places $v$ and a 
     $v$-unit at the places $v|pc$.}
\end{array}\right.
$$
For each such representative $s$ there is a diagram of complex tori
$$
\begin{CD}
	A_{g\tau}=\C^{\ep}/\La_{g\tau} @>{j(g,\tau)}>> 
	\C^{\ep}/s\La_{\tau} @>{\pi_{s}}>> \C^{\ep}/\La_{\tau}\\
	@.  @|  @| \\
	{} @.  A_{x^\si}(\C) @.  A_x(\C)\\
\end{CD}
$$
where $s=dgu$ under \eqref{eq:decomp} and $\pi_{s}$ is the natural 
quotient map arising from the inclusion 
$s\La_{\tau}\subset\La_\tau$. The element 
$g\in\GL_{2}^{+}(\R)$ is defined by $g\tau\in\gotH$ only up to an
element in $\calO_{K,c}^{\times}$.
Choose $p$ and  a place $v$ over $p$ in  a number field $L$ 
large enough so that for each $s$ the triple 
$(g\tau,v,e)$ is a $p$-ordinary test triple and that the isogenies 
$\pi_{s}$ are defined over $L$.

\begin{lem}\label{lem:compperiods}
         With the above notations, it is possible to choose for every 
         $\si\in\Gal(H_c/K)$ an invariant 1-form on $A_{x^\si}$ that generates 
         ${\calL}(x^\si)\otimes\calO_{(v)}$ 
         and for which
	\begin{enumerate}
 	\item  $\Om_{\infty}(g\tau)\sim_{\calO_{(v)}^\times}j(g,\tau)\Om_{\infty}(\tau)$;
	\item  $\Om_{p}(x^{\si})\sim_{\calO_{(v)}^\times}\Om_{p}(x)$.
	\end{enumerate}
\end{lem}

\pf Take $\om_o\in H^0(A_{x}(\C),\calL(x)\otimes\C)$. The quotient map $\pi_{s}$ 
is the identity on (co)tangent spaces and commutes with the action of 
the endomorphisms. Thus 
$\om_s=\pi_s^*(\om_o)\in H^0(A_{x^\si}(\C),\calL(x^\si))$ and
$p(\om_s,g\tau)=j(g,\tau)p(\om_o,\tau)$. Furthermore, $p$ doesn't 
divide the degree of $\pi_s$ and so $\pi_s^*$ is an isomorphism between 
the natural $p$-adic structures on the spaces of invariant forms. 
This proves part 1.

For part 2 observe that the reduction mod $p$ of the dual map 
$\pi_s^t$ gives an isomorphism of the rank 1 tate module quotient $T$ of \S 3.2. 
Thus $\pi_s^*(\om_{u}(P))$ 
is a universal form on the deformations of $\wtA_{x^\si}$ by formula 
\eqref{eq:omcomp} and the equality follows. \qed

If $s$ and $s^{\prime}=s\la zu$ with $\la\in K^{\times}$, 
$z\in\C^{\times}$ and $u\in\wh{\calO}_{K,c}^\times$ are two normalized 
representants of the same $\si\in\CC^\sharp_{c}$ a comparison of the 
relations in lemma \ref{lem:compperiods} for the decompositions 
$s=dgr$ and $s^{\prime}=(\la d)(gz)(ru)$ shows that 
$\om_{s^{\prime}}\sim_{\calO_{K,c}^{\times}}z\om_{s}$. Therefore the 
construction of $\om_{s}$ can be extended modulo 
$\calO_{K,c}^{\times}$-equivalence to all $s\in\ideles$ by setting
\begin{equation}
    \om_{s\la zu}\sim_{\calO_{K,c}^{\times}}z\om_{s}
    \quad\mbox{for all $\la\in K^{\times}, z\in\C^{\times}, 
    u\in\wh{\calO}_{K,c}^\times$ and $s$ normalized.}
    \label{eq:definoms}
\end{equation}
Let $f\in\M_{2\ka,0}(\De,N)$ and normalize the invariant form as in 
proposition \ref{teo:Thkalgebraic}. For all integers $r\geq0$ such 
that $(\calO_{K,c}^{\times})^{2(\ka+r)}=1$ define a function
${c}_{(r)}(f,x)\colon\ideles\map\C$ as
$$
   {c}_{(r)}(f,x)(s)=\frac{\de_{2\ka}^{(r)}f(g\tau)}{p(\om_{s},g\tau)^{2(\ka+r)}}
$$
where $s=dgu$ as above.

\begin{pro}\label{th:meascrfx}
Suppose that $f$ is defined over $\calO_{(v)}$ and assume that 
$(\calO_{K,c}^{\times})^{2(\ka+r)}=1$. Then  
$c_{(r)}(f,x)\in{S}_{(2(\ka+r),0)}(\calO_{K,c})\cap\calS(cR_{K},\calO_{v})$.
\end{pro}

\pf The modular relation for $c_{(r)}(f,x)$ follows at once from 
\eqref{eq:definoms} and the definition
since $gz\tau=g\tau$.
For an idele $s$ satisfying the conditions \eqref{eq:weilrel} for 
$\gotn=(pc)$ the 
invariant form $\om_{s}$ satisfies proposition \ref{teo:Thkalgebraic} 
and then theorem \ref{thm:equality} together with lemma \ref{lem:compperiods}
shows that as a $p$-adic $K^{\times}$-modular form 
${c}_{(r)}(f,x)$ has coefficients in 
$L_{v}$ and in fact belongs to the unit ball. \qed

Assume that $\calO_{K,c}^{\times}=\{\pm1\}$. 
Let $\mu_{f,x}$ be the $p$-adic distribution on $\Z_p$ with values in 
$\overline{S}(c\calO_K;L_v)$ such that $m_r(\mu_{f,x})=\wh{c}_{(r)}(f,x)$
and let $\mu_{\chi,\xi}$ be the $p$-adic measure associated to a choice of Gr\"ossencharakters
$\chi\in\Xi_{(-2\ka,0)}(\calO_{K,c})$, $\xi\in\Xi_{(-2,0)}(\calO_{K,c})$ as in the discussion after
lemma \ref{le:UnBoundPr}.

\begin{thm}\label{th:measexists}
     There exist a $p$-adic field $F$ and a $p$-adic measure $\mu(f,x;\chi,\xi)$ on 
     $\Z_p$ with values in $\calO_F$ such that
     $$
     m_r\left(\mu_{[\mu_{f,x},\mu_{\chi,\xi}]}\right)
     =\left\{
	\begin{array}{ll}
		0 & \mbox{if $0\leq r$ is odd,}  \\
		(h^\sharp_c)^{-1}\Om_p^{-2(\ka+l)}\binom{2l}{l}m_l(\mu(f,x;\chi,\xi)) & 
		\mbox{if $0\leq r=2l$ is even,}
	\end{array}
	\right.
     $$
\end{thm}

\pf Let $F$ be large enough to contain $L_v$, the field of values of $\chi$ and $\xi$
and the $p$-adic period $\Om_p$.
The expression follows from Lemma \ref{th:measuremix} and the fact that for a suitable
choice of  representants  for $\CC^{\sharp}_{c}$ we have, combining the definition 
\eqref{eq:padicpair} with theorem \ref{thm:equality}, proposition \ref{th:spiden} and lemma \ref{lem:compperiods},
$\scalq{\wh{\chi}{\wh{\xi}}^{l}}{\wh{c}_{(l)}(f,x)}=
(h^\sharp_c)^{-1}\Om_p^{-2(\ka+l)}\sum_\si\wh{\chi}{\wh{\xi}}^{l}(s)b_l(x^\si)$. 
Finally, each term ${\wh{\xi}}^{l}(s)b_l(x^\si)$ is the $l$-th moment of a 
suitable $p$-adic measure on $\Z_p$ because the identification
$\sum_{n=0}^\infty(b_n(x^\si)/n!)T^n=\sum_{n=0}^\infty a_nU^n$ with $a_n\in\calO_F$ through the substituition $U=e^T-1$ yields an identification
$\sum_{n=0}^\infty(b_n(x^\si)/n!)z^nT^n=\sum_{n=0}^\infty a_nV^n$ where 
$V=(U+1)^z-1$ and this substitution preserves $\calO_F$-integrality when $z$ is a unit in a field with residue field $\F_p$. Conclude using the linearity of measures.
\qed

\subsection{Special $L$-values}

For $f\in M_{2\ka,0}^\infty(\De,N)$, let 
$\phi_f\in L^2(D_\Q^\times\bs D_\A^\times)$ be the usual
$\wh{\calR}_N^\times$-invariant $C^\infty$ lift of $f$ to $D_\A^\times$.
Namely,
$\phi_f(d)=f(g_\infty\cdot i)j(g_\infty,i)^{-2\ka}\det(g_\infty)^\ka$ if
$d=d_{\Q}g_\infty u$ under \eqref{eq:decomp}.
The Lie algebra $\gotg=\gotg\gotl_2\simeq\Lie(D_\infty^\times)$ 
acts on the $\C$-valued $C^\infty$ functions on $D_\A^\times$ by 
$(A\cdot\varphi)(d)=\left.\frac{d}{dt}\varphi(de^{tA})\right|_{t=0}$. 
By linearity and composition the action extends to the complexified 
universal enveloping algebra $\gotA(\gotg)_{\C}$. Let
$$
I=\left(
\begin{array}{cc}
	1 & 0  \\
	0 & 1
\end{array}
\right),
\qquad
H=\left(
\begin{array}{cc}
	0 & -i  \\
	i & 0
\end{array}
\right),
\qquad
X^\pm=\frac12\left(
\begin{array}{cc}
	1 & \pm i  \\
	\pm i & -1
\end{array}
\right)
$$
be the usual eigenbasis of $\gotg_{\C}$ for the adjoint action
of the maximal compact subgroup
$$
\SO(2)=\left\{\hbox{$r(\th)=\left(
\begin{array}{cc}
	\cos\th & -\sin\th  \\
	\sin\th & \cos\th
\end{array}
\right)$ such that $\th\in\R$}\right\}.
$$
Since ${\rm Ad}(r(\th))X^\pm=e^{\mp 2i\th}X^\pm$, we have
$X^\pm\cdot\varphi_{f}\in M_{2\ka\pm2,0}^\infty(\De,N)$.
A standard computation (e.g. \cite[\S\S2.1--2]{Bu96})
links the Lie action to the Maass operators of 
\S\ref{ss:maass}, namely
$$
X^+\cdot\phi_f=-4\pi\phi_{\de_{2\ka}f}.
$$
For $r\geq0$ let
\begin{equation}
	\phi_r=\left(-\frac{1}{4\pi}X^{+}\right)^{r}\cdot\phi_f=
	\phi_{\de_{2\ka}^{r}f}.
	\label{eq:defphir}
\end{equation}

\begin{dfn}\label{th:Jintegral}
	\rm 
	Let $f\in M_{2\ka,0}(\De,N)$, $\xi\in{\Xi}_\uw(c\calO_{K})$ for a 
	weight $\uw$ such that $\vass{\uw}=0$
	and $\tau=t+iy\in\CM_{\De,K}$ with $c_{\tau,N}=c$ and associated normalized
	embedding $\jmath$. For each $r\geq0$, let
	$$
	J_r(f,\xi,\tau)=\int_{\ideles/K^\times\R^\times}\phi_r(\jmath(t)d_\infty)\xi(t)\,dt
	$$
	where $d_\infty=\smallmat{y^{1/2}}{ty^{1/2}}0{y^{-1/2}}$ and $dt$ is the Haar measure on 
	$\ideles$ whose archimedean component is normalized so that 
         $\mathrm{vol}(\C^\times/\R^\times)=\pi$ and such that the local grups of units have volume 1
         (hence $m_c=\mathrm{vol}(\wh{\calO}_{K,c}^\times)=
         [(\calO_K/c\calO_K)^\times\colon(\Z/c\Z)^\times]^{-1}$).
\end{dfn}

We show that $J_r(f,\xi,\tau)$ can be expressed in terms 
of the pairing introduced in \S\ref{se:padicforms}.
Write $w_{K,c}=\vass{\calO_{K,c}^{\times}}$.

\begin{thm}\label{th:compJ}
   Let $f\in M_{2\ka,0}(\De,N)$ and $\xi\in\Xi_{(w,-w)}(\calO_{K,c})$.
   Assume that $(\calO_{K,c}^{\times})^{2w}=1$.
   Then
   $$
   J_r(f,\xi,\tau)=
   	\frac{\pi m_c}{w_{K,c}}h_{c}^{\sharp}y^{-w}\Om_\infty(\tau)^{-2w}
   	{\scal{c_{(r)}(f,x)}{\xi\vvass{N_{K/\Q}}^{-w}}}.
   $$
\end{thm}

\pf Since the integrand function is right $\wh{\calO}^\times_{K,c}$-invariant, we have
$J_r(f,\xi,\tau)=m_c\int_{\ideles/K^\times\R^\times\wh{\calO}^\times_{K,c}}
\phi_r(\jmath(t)d_\infty)\xi(t)\,dt$. For a chosen set of 
representatives $\{s_\si\}$ of $\CC^{\sharp}_{c}$ there is a decomposition
$$
\ideles/K^\times\R^\times\wh{\calO}^\times_{K,c}=
\bigcup_{\si\in\CC^{0}_{c}}\C^\times s_\si/\R^\times\calO_{K,c}^{\times}
\qquad
\mbox{(disjoint union).}
$$
Therefore,
$J_r(f,\xi,\tau)=m_c\sum_{\si\in\CC^{\sharp}_{c}}\xi(s_\si)\int_{\C^\times/\R^\times\calO_{K,c}^{\times}}
\phi_r(\jmath(s_\si z)d_\infty)\xi_\infty(z)\,d^\times z$.

Since the standard normalized embedding of $\C$ in $\M_2(\R)$ is 
$\rho e^{i\th}\mapsto\rho r(\th)$, we can write 
$D_\infty^\times\ni\jmath(z)=\rho d_\infty r(\th)d_\infty^{-1}$. 
Therefore
$\int_{\C^\times/\R^\times\calO_{K,c}^{\times}}\phi_r(\jmath(s_\si z)d_\infty)
\xi_\infty(z)\,d^\times z=w_{K,c}^{-1}\int_0^{\pi}\phi_r(\jmath(s_\si d_\infty)r(\th))
\xi_\infty(e^{i\th})\,d\th=w_{K,c}^{-1}\phi_r(\jmath(s_\si d_\infty))
\int_0^{\pi}e^{-2i(\ka+r)\th}e^{-2iw\th}\,d\th$ 
and
\begin{equation}
    J_r(f,\xi,\tau)=\left\{
    \begin{array}{ll}
        \pi m_cw_{K,c}^{-1}\sum_{\si\in\CC^{\sharp}_{c}}
        \xi(s_\si)\phi_r(\jmath(s_\si d_\infty)) & \hbox{if $w=-\ka-r$}  \\
        0 & \hbox{otherwise}
    \end{array}\right..
    \label{eq:formforJ}
\end{equation}
Note that this proves the claimed formula when $w\neq-\ka-r$ since the inner product in its right hand side vanishes in this case. Thus, we may now assume that $w=-\ka-r$. Put 
$I_\si=\xi(s_\si)\phi_r(s_\si d_\infty)$ and write 
$s=s_\si=d_sg_su_s$ under 
\eqref{eq:decomp} and $\tau_s=g_s\tau$. Note that 
$\vvass{N_{K/\Q}(s)}=\det(g_s)$. Under the hypothesis 
$(\calO_{K,c}^{\times})^{2(\ka+r)}=1$ we have
\begin{eqnarray*}
	I_\si & = & \xi(s)\de_{2\ka}^{(r)}f(g_sd_\infty\cdot i)j(g_s d_\infty,i)^{-2(\ka+r)}\det(g_s)^{\ka+r}  \\
	 & = & y^{\ka+r}\xi(s)\de_{2\ka}^{(r)}f(\tau_{s})j(g_s,\tau )^{-2(\ka+r)}\vvass{N_{K/\Q}(s)}^{\ka+r}  \\
	 & = & y^{\ka+r}\xi(s)c_{(r)}(f,x)(s)
	p(\om_s,\tau_{s})^{2(\ka+r)}j(g_s,\tau )^{-2(\ka+r)}\vvass{N_{K/\Q}(s)}^{\ka+r}  \\
	 & = & y^{\ka+r}\xi(s)c_{(r)}(f,x)(s)
	p(\om_{s_o},\tau_{s_o})^{2(\ka+r)}j(g_{s_o},\tau )^{-2(\ka+r)}\vvass{N_{K/\Q}(s)}^{\ka+r}  \\
	 & = & y^{\ka+r}\Om_\infty(\tau)^{2(\ka+r)}\xi(s)c_{(r)}(f,x)(s)\vvass{N_{K/\Q}(s)}^{\ka+r}  \\
\end{eqnarray*}
where $s_o$ is a normalized representant. 
It is now clear that the formula follows. \qed

\begin{dfn}
     Let $M$ be a proper divisor of $N$,  $x\in\CM(\De,N;\calO_{K,c})$ and 
     $\pr x\in\CM(\De,M;\calO_{K,\pr c})$ the image of $x$ under the natural 
    quotient map. 
    A character $\xi\in\Xi_{\uw}(\calO_{K,c})$ 
    is called $(x,M)$-primitive if it is 
    not trivial on $\wh{\calO}_{K,\pr c}^{\times}$.
\end{dfn}

For a divisor $d$ of $N/M$ there is an embedding
$\iota_{\De,d}:M_{2\ka,0}(\De,M)\mmap M_{2\ka,0}(\De,N)$.
When $\De=1$ the embedding is simply $f(z)\mapsto f(dz)$. 
When $\De>1$ the explicit description of 
$\iota_{\De,d}$ is less immediate, e.÷g. \cite[\S3]{MoTe99}. We 
denote $M_{2\ka,0}(\De,N)^{M-\mathrm{old}}$ the span of the 
images of the embeddings $\iota_{\De,d}$ for all $d$.
After theorem \ref{th:compJ} the following result can be read  
as an orthogonality statement between primitive characters and 
$K^{\times}$-modular forms arising from oldforms. 

\begin{pro}\label{th:oldforms}
    Let  $\tau\in\CM_{\De,K}$ and $x\in\CM(\De,N;\calO_{K,c})$
    be the point represented by $\tau$. Let $f\in M_{2\ka,0}(\De,N)^{M-\mathrm{old}}$
    and suppose that $\xi\in\Xi_{(-\ka-r,\ka+r)}(\calO_{K,c})$ is $(x,M)$-primitive.
    Then $J_r(f,\xi,\tau)=0$.
\end{pro}

\pf Consider again the first expression in \eqref{eq:formforJ}. 
Let $\pr x\in\CM(\De,M;\calO_{K,\pr c})$ the point image of $x$ and choose a 
system of representants $\{s_{\pr\si}\}$ of $\CC^{\sharp}_{\pr c}$
and a system of representants $\{r_{i}\}$ of 
$\wh{\calO}^\times_{K,\pr c}/\wh{\calO}^\times_{K,c}$. Then the set of products 
$\{s_{\si}r_{i}\}$ is a system of representatives of $\CC^{\sharp}_{c}$ 
and since $d_{\infty}$ commutes with each $r_{i}$ and $f$ is $M$-old 
we obtain the expression
$$
J_{r}(f,\xi,\tau)=\frac{\pi m}{w_{K,c}}
\left(\sum_{\wh{\calO}^\times_{K,\pr c}/\wh{\calO}^\times_{K,c}}
\xi(r_{i})\right)\left(\sum_{\pr\si\in\CC^{\sharp}_{\pr c}}
\xi(s_{\pr\si})\phi_r(s_{\pr\si}d_\infty)\right)
$$
which vanishes because $\xi$ is non trivial on 
$\wh{\calO}^\times_{K,\pr c}/\wh{\calO}^\times_{K,c}$. \qed

We shall assume from now on that the modular form $f$ is a 
holomorphic newform with associated automorphic representation 
$\pi^{D}=\pi_{f}$. Let $\pi$ be the automorphic representation of 
$\GL_{2}(\A)$ corresponding to $\pi^{D}$ under the Jacquet-Langlands 
correspondence.

Other than the Weil representation $r_{\psi}$ of $\SL_{2}(\A)$, the adelic 
Schwartz-Bruhat space $\calS_{\A}(D)=\bigotimes_{p\leq\infty}\calS_{p}$
supports the unitary representation of $\GO(D)(\A)$ given by
$$
L(h)\varphi(x)=\vvass{\nu_{0}(h)}_{\A}^{-1}\varphi(h^{-1}x),
\qquad x\in D_{\A}.
$$
We assume that the archimedean space $\calS_{\infty}$ consists
only of the Schwartz functions on $D_{\infty}$ which are 
$K^{1}_{\infty}\times K^{1}_{\infty}$-finite under the action 
of $D^{\times}\times D^{\times}$ via the group $\GO(D)$ 
(\S\ref{se:CMpts}). Here $K^{1}_{\infty}$ is the maximal compact subgroup of 
$\jmath(K^{\times}\otimes\R)\subset D_{\infty}^{\times}\simeq\GL_{2}(\R)$. 
As explained in \cite[\S5]{HaKu92}, the two representations mingle into
one single representation, still denoted $r_{\psi}$, of the group
$R(D)=
\{\hbox{$(g,h)\in\GL_{2}\times\GO(D)$ such that $\det(g)=\nu_{0}(h)$}\}$
given by $r_{\psi}(g,h)\varphi=r_{\psi}(g_{1})L(h)\phi$ where 
$g_{1}=g\smallmat{1}{}{}{\nu_{0}(h)}^{-1}$. Note that 
\begin{itemize}
	\item the assignment $(g,h)\mapsto(g_{1},h)$ sets up an isomorphism 
          $R(D)\stackrel{\sim}{\map}\SL_{2}\ltimes\GO(D)$; 
	\item the group $R(D)$ is naturally a subgroup of the symplectic 
	      group $\Sp(W)$, where $W=P\otimes D$ with $P$ the standard 
	      hyperbolic plane, via $(g,h)x\otimes y=gx\otimes h^{-1}y$.
\end{itemize}
The groups $(\SL_{2},{\rm O}(D))$ form a dual reductive pair in 
$\Sp(W)$ and the extended Weil representation $r_{\psi}$ allows to realize the 
theta correspondence between the similitude groups. The theta kernel 
associated to a choice of $(g,h)\in R(D)$ and $\varphi\in\calS_{\A}(D)$ is
$\vartheta(g,h;\varphi)=\sum_{d\in D}r_{\psi}(g,h)\varphi(d)$.
The theta lift to $\GO(D)$ of a cuspidal automorphic form $F$ on 
$\GL_{2}(\A)$ is the automorphic form on $\GO(D)(\A)$ given by
\begin{equation}
     \theta_{\varphi}(F)(h)=\int_{\SL_{2}(\Q)\backslash\SL_{2}(\A)}
     \vartheta(gg^{\prime},h;\varphi)F(gg^{\prime})\,dg^{\prime}
     \label{eq:thetalift}
\end{equation}
where $\det(g)=\nu_{0}(h)$ and $dg^{\prime}$ is induced by a choice
of a  Haar measure $dg=\prod dg_{p}$ on $\GL_{2}(\A)$. 
A straightforward substitution yields
\begin{equation}
	\theta_{r_{\psi}(g_1,h_1)\varphi}(F)(h)=
	\theta_{\varphi}(\pi(g_1^{-1})F)(hh_1),
	\qquad
	\forall (g_1,h_1)\in R(D).
	\label{eq:tlautom}
\end{equation}
An automorphic form $\Phi$ on $\GO(D)(\A)$ pulls back via the map 
$\varrho$  of \eqref{eq:similitudes} to an automorphic form 
$\widetilde{\Phi}$ on the 
product group $D^{\times}\times D^{\times}$. Let $\widetilde{\Theta}(\pi)$ 
be the space of automorphic forms on $D^{\times}\times D^{\times}$ which are 
pull-backs of theta lifts \eqref{eq:thetalift} with $F\in\pi$. 
If $\check{\pi}^{D}$ denotes the contragredient 
representation of $\pi^{D}$ the crucial result is, 
with a slight abuse of notation, the following, \cite{Shimi72}.

\begin{thm}[Shimizu]\label{th:Shimizu}
	$\widetilde{\Theta}(\pi)=\pi^{D}\otimes\check{\pi}^{D}$.
\end{thm}

\begin{rems}\label{rm:onchoice}
   \rm 
   \begin{enumerate}
     \item In our case of interest $\pi^D=\check{\pi}^{D}$.
     \item  The Schwartz functions, hence the theta lifts 
             $\widetilde{\theta}_{\varphi}(F)$, are 
	    $K^{1}_{\infty}\times K^{1}_{\infty}$-finite. 
	    Thus, in Shimizu's theorem the representation space $\pi^{D}$ consists of 
	    $K^{1}_{\infty}$-finite automorphic forms. Note that the functions 
	    $\pi(d_{\infty})\phi_{r}$ are $K^{1}_{\infty}$-finite.
   \item An explicit version of Shimizu's theorem has been worked out by Watson \cite{Wat03}, 
            see also \cite[\S3.2]{Pra06} and \cite[\S12]{HaKu92}. Namely, if 
            $\varphi=\otimes_{p\leq\infty}\varphi_p$ is chosen as
            \begin{equation}
            \varphi_\infty(z_1,z_2)=
            \frac{(-1)^\ka}\pi z_2^{2\ka}e^{-2\pi(z_{1}\bar{z}_{1}+z_{2}\bar{z}_{2})},
            \quad\varphi_p=
            \frac{ \mathrm{ch}_{\calR_{N}\otimes\Z_p}}{\mathrm{vol}((\calR_{N}\otimes\Z_p)^{\times})}
            \label{eq:choiceofphi}
           \end{equation}
           where $z_1$ and $z_2$ are the complex coordinates in $D_\infty$ of \S1.3, then
           $$
           \pi(d_{\infty})\phi_{f}\otimes\pi(d_{\infty})\phi_{f}=\widetilde{\theta}_\varphi(F)
           $$
           where $F\in\pi$ is the adelic lift of an eigenform normalized so to have an equality of 
           Petersson norms 
           $\langle\pi(d_{\infty})\phi_{f},\pi(d_{\infty})\phi_{f}\rangle=\langle F,F\rangle$.
  \end{enumerate}   
\end{rems}

Let $\underline{\xi}=(\xi,\xi^{\prime})\in\Xi_\uw(c\calO_{K})\times\Xi_{\uw^{\prime}}(c\calO_{K})$ 
thought of as a character of the torus $\ideles\times\ideles$. 
Let $\tilde{H}(t)$ be any function on $\ideles\times\ideles$ such that 
$\tilde{H}(t)\underline{\xi}(t)$ is $(K^{\times}\R^{\times})^{2}$-invariant. 
Following \cite[\S14]{HaKu91} \cite[\S1.4]{Harris93} we let
$$
    L_{\underline{\xi}}(\tilde{H})=
    \int_{(K^{\times}\R^{\times}\backslash\ideles)^{2}}
    \tilde{H}(t)\underline{\xi}(t)\,dt.
$$
In particular, for $\xi$ as in definition \ref{th:Jintegral},
$$
L_{(\xi,\xi)}( \pi(d_{\infty})\phi_r\otimes\pi(d_{\infty})\phi_r)=
J_r(f,\xi,\tau)^2.  
$$
When $\xi=\xi^\prime$ is unitary, $\vass{\uw}=\vass{\uw^{\prime}}=0$,
the integral $L_{\underline{\xi}}(\widetilde{\theta}_{\varphi}(F))$
can also be read, via the map $\al$ of \eqref{eq:similitudes}, 
as the Petersson scalar product of two automorphic forms
on the similitude group $T=G({\rm O}(K)\times{\rm O}(K^{\perp}))$ 
associated with the decomposition $D=K\oplus K^\perp$, namely
$L_{(\xi,\xi)}(\widetilde{\theta}_{\varphi}(F))=\int_{T(\Q)T(\R)\backslash T(\A)}
\widetilde{\theta}_{\varphi}(F)((a,b))\xi(b)\,d^\times ad^\times b$,
where $\alpha(t)=(a,b)$.
Thus the seesaw identity 
\cite{Kudla84} associated with the seesaw dual pair
$$
\begin{array}{ccc}
	\GL_{2}\times\GL_{2} &  & \GO(D) \\
	\uparrow & \mbox{{\Huge $\times$}} & \uparrow  \\
	\GL_{2} &  & G({\rm O}(K)\times{\rm O}(K^{\perp}))
\end{array}
$$
identifies, up to a renormalization of the Haar measures, 
the value $L_{(\xi,\xi)}(\widetilde{\theta}_{\varphi}(F))$ with a scalar 
product on $\GL_{2}$,
\begin{equation}
	L_{\underline{\xi}}(\widetilde{\theta}_{\varphi}(F))=
	\int_{\GL_{2}(\Q)\A^{\times}\backslash\GL_{2}(\A)}
	F(g)\theta_{\varphi}^t(1,\xi)(g,g)\,dg,
	\label{eq:RankinSelberg}
\end{equation}
where $\theta_{\varphi}^t$ denotes the theta lift to 
$\GL_{2}\times\GL_{2}$. If $\varphi$ is 
split and primitive, i.e. admits a decomposition 
$\varphi=\varphi_{1}\otimes\varphi_{2}$ under 
$D_{\infty}=(K\oplus K^\perp)\otimes\R$ and 
each component decomposes in a product of local factors, 
$\varphi_{i}=\bigotimes_{p\leq\infty}\varphi_{i,p}$ for $i=1,2$,
then $\theta_{\varphi}^t(1,\xi)$ splits as a product of two separate lifts. In fact
$$
\theta_{\varphi}^t(1,\xi)(g_{1},g_{2})=
E(0,\Phi,g_{1})\theta_{\varphi_{2}}(\tilde{\xi})(g_{2})
$$
where: 
\begin{itemize}
    \item  $E(0,\Phi,g)$ is the value at $s=0$ of the holomorphic 
           Eisenstein series attached to the unique flat section (\cite[\S3.7]{Bu96})
           extending the function $\Phi(g)=r_\psi(g,k)\varphi_1(0)$ 
	   where $k\in\ideles$ is such that $N(k)=\det(g)$ and 
	   $r_\psi$ denotes here the extended adelic Weil representation 
           attached to $K$ as a normed space (Siegel-Weil formula),
    \item  $\theta_{\varphi_{2}}(\xi)(g)$ is a binary form in the automorphic representation 
           $\pi(\xi)$ of $\GL_{2}$ attached to $\xi$.
\end{itemize}
This expression yields a relation between the right hand side 
of \eqref{eq:RankinSelberg} and the value at the centre of symmetry of a
Rankin-Selberg convolution integral. If the Whittaker function $W_F$
of $F$ decomposes as a product of local Whittaker functions, 
the Rankin-Selberg integral admits an Euler 
decomposition \cite{Ja72} and $L_{(\xi,\xi)}(\widetilde{\theta}_{\varphi}(F))$
is equal to the value at $s=\frac12$ of the analytic continuation of
$$
\frac 1{h_K}\prod_{q\leq\infty}L_q(\varphi_q,\xi_q,s),
$$
where
\begin{multline}
L_q(\varphi_q,\xi_q,s)=
\int_{K_q}\int_{\Q_q^\times}
W^{\psi_q}_{F,q}\left(\left(\begin{array}{cc}a & 0 \\0 & 1\end{array}\right)k\right)
W^{\psi_q}_{\theta_{\varphi_{2},q}}\left(\left(\begin{array}{cc}-a & 0 \\0 & 1\end{array}\right)k\right)\cdot\\
\Phi^s_{\varphi_1,q}\left(\left(\begin{array}{cc}a & 0 \\0 & 1\end{array}\right)k\right)\inv{\vass{a}}\,
d^\times a\,dk_q.
\label{eq:localterm}
\end{multline}
The local measures are normalized so that $K_\infty=\SO_2(\R)$ has volume $2\pi$ and
$K_q=\GL_2(\Z_q)$ has volume $1$ for finite $q$.
Also $W_{\theta_{\varphi_{2}}}$ is the Whittaker function
and $\Phi^s(g)=\vvass{a}^{s-\frac12}\Phi(g)$ if $g=nak$ under the $NAK$-decomposition 
where $\vvass{\smallmat a{}{}b}=\vass{a/b}$.
Since the local term \eqref{eq:localterm} does not vanish and for almost all $q$ is the local Euler factor of some automorphic $L$-function, one obtains, as in \cite{Harris93, HaKu91}, a version of Waldspurger's result \cite{Waldsp85}. Namely,
$$
   L_{\underline{\xi}}(\widetilde{\theta}_{\varphi}(F))=\left.
   \Lambda(\varphi,\xi,s)L(\pi_K\otimes\xi,\frac s2)L(\eta_K,2s)^{-1}\right|_{s=1/2},
$$
where $\Lambda(\varphi,\xi,s)$ is a finite product of local integrals, $\pi_K$ is the base change to $K$ of the automorphic representation $\pi$ and $L(\eta_K,2s)$ is the Dirichlet $L$-function attached to 
$\eta_K$, the quadratic character associated to $K$

When $\varphi^f=\bigotimes_{p<\infty}\varphi_p$ and $F$ are chosen as in Remark \ref{rm:onchoice}.3 
the local non-archimede\-an terms in the Rankin-Selberg integral have been explicitely computed by Prasanna \cite[\S3]{Pra06} under the simplifying assumptions that $N$ is squarefree, $c=1$ and $\xi$ is unramified. The effect of these assumptions is that
\begin{enumerate}
  \item the local component of $\xi$ can be written either as  
           $\xi_q=(\xi_q^{\rm sp},(\xi_q^{\rm sp})^{-1})$ for some unramified character 
           $\xi_q^{\rm sp}$ of $\Q_q^\times$ at a prime $q$ split in $K$ under the isomorphism 
           $(K\otimes\Q_q)^\times\simeq\Q_q^\times\times\Q_q^\times$, or as
           $\xi_q=\xi^{\rm in}_q\circ{\rm N}_{K_q/\Q_q}$ for an unramified character 
           $\xi_q^{\rm in}$ of $\Q_q^\times$ at a prime $q$ inert in $K$, or as
           $\xi_q=\xi^{\rm rm}_q\circ{\rm N}_{K_q/\Q_q}$ at a ramified prime $q$ where 
           $\xi_q^{\rm rm}$ is the unramified character of $\Q_q^\times$ obtained by a trivial extension; 
  \item at a prime $q|N\De$ the local component $\pi_q$ is equivalent to the special representation
           $\si(\vass{\cdot}^{\frac12+it_q},\vass{\cdot}^{-\frac12+it_q})$ with $q^{2it_q}=1$.
\end{enumerate}
Note that the former condition remains true for a split prime $q$ that does not divide $c$ and 
$\xi\in\Xi_\uw(R_{K,c})$ with $\vass\uw=0$. Thus we can apply Prasanna's computations to this more general case to write down a formula in which only the local factors at primes in $\Si=\left\{q|cN^\prime\right\}$ are left implicit, where
$N=N_{\rm sf}(N^\prime)^2$ and $N_{\rm sf}$ is square-free.
Namely,
\begin{equation}
L_{\underline{\xi}}(\widetilde{\theta}_{\varphi_\infty\otimes\varphi^f}(F))
=\left.\frac{V_N}{h_K}\la_\infty(\varphi_\infty,\xi_\infty,s)\left(\prod_{q\leq\infty}\nu_q(\xi_q)\right)
L(\pi_K\otimes\xi,\frac s2)L(\eta_K,2s)^{-1}\right|_{s=1/2}
\label{eq:finalexp}     
\end{equation}        
where $V_N=\prod_q{\mathrm{vol}((\calR_{N}\otimes\Z_q)^{\times})}$,
 $\la_\infty(\varphi_\infty,\xi_\infty,s)=\vass{{\rm N}u}^{\frac12}_\infty\xi_\infty(z_u)^{-1}
L_\infty(\varphi_\infty,\xi_\infty,\frac s2)$ ($z_u$ denotes the complex coordinate of $u$ in the chosen identification $(Ku)\otimes\R\simeq\C$)
and 
$$
\begin{cases}
    \nu_q(\xi_q)=\xi_q^{\rm sp}(q)^{n_1-n_2}  & \text{if $q$ splits, $q\notin\Si$, $(q,N_{\rm sf})=1$ }, \\
    \nu_q(\xi_q)=-\frac{1}{q+1}q^{-\frac12+t_q+s}\xi_q^{\rm sp}(q)^{n_1-n_2}  & \text{if $q|N_{\rm sf}$}, \\
    \nu_q(\xi_q)=\xi_q^{\rm in}(q)^{-2n} & \text{if $q$ is inert, $q\notin\Si$}, \\
    \nu_q(\xi_q)=\xi_q^{\rm rm}(-{\rm N}\inv u) & \text{if $q$ ramifies}, \\
    \nu_\infty(\xi_\infty)=\xi_\infty(z_u) &  \\
    \end{cases}
$$
where the ideal $\calJ$ of proposition \ref{prop:decomporder} in $K\otimes\Q_q$ is generated by $q^{-n}$ when $q$ is inert and decomposes as $q^{-n_{q,1}}\Z_q\times  q^{-n_{q,2}}\Z_q$ under 
$K\otimes\Q_q\simeq\Q_q\times\Q_q$ when $q$ is split.

\begin{rem}
     \rm It is clear that $\nu_q(\xi_q)=1$ for almost all $q$.
     The local terms $\la_q(\nu_q)$ do depend 
     on the choice of $u$ in \eqref{eq:Dsplit} (replacing $u$ with $xu$
     the local Whittaker function $W^{\psi_q}_{\theta_{\varphi_{2},q}}$ gets modified by the factor 
     $\vass{{\rm N}x}^{-\frac12}_q\xi_q(x)^{-1}$), but the quantity
     $$
     \nu(\xi,\tau,s)=\prod_{q\leq\infty}\nu_q(\xi_q)
     $$
     depends only on $\xi$ and the chosen embedding $\jmath:K\rightarrow D$.
\end{rem}

For a pair of non-negative integers $(m,q)$ consider the function of two complex variables
$\varphi^{(l,q)}(z_1,z_2)=(z_1\bar{z}_1)^{l}{z}_2^{q}e^{-2\pi(z_1\bar{z}_1+z_2\bar{z}_2)}$.

\begin{lem}
     Let $\varphi(z)=(z\bar z)^le^{-2\pi z\bar z}$. Then the Fourier transform of $\varphi$ is
     $$
     \hat\varphi(w_1+w_2i)=e^{-2\pi(w_1^2+w_2^2)}\sum_{0\leq\al+\be\leq l}\ga(\al,\be;l)w_1^{2\al}w_2^{2\be},
     $$
     where
     $$
     \ga(\al,\be;l)=\sum_{\substack{j+k=l\\ \al\leq j, \be\leq k}}(-4\pi)^{\al+\be-l}\binom lj\binom{2j}{2\al}\binom{2k}{2\be}     
     (2j-2\al-1)!!(2k-2\be-1)!!
     $$
\end{lem}

\pf
One has
$\hat\varphi(w_1+w_2i)=2\int_{\R^2}e^{4\pi i(w_1x_1+w_2x_2)}(x_1^2+x_2^2)^le^{-2\pi(x_1^2+x_2^2)}\,dx_1dx_2=
2\sum_{j+k=l}\binom lj\left(\int_{\R}e^{2\pi iw_1x_1}x_1^{2j}e^{-2\pi x_1^2}\,dx_1\right)
\left(\int_{\R}e^{2\pi iw_2x_2}x_2^{2k}e^{-2\pi x_2^2}\,dx_2\right)$
and the result follows from
$ \int_{\R}e^{4\pi itx-2\pi x^2}x^{2v}\,dx=\frac1{\sqrt{2}}e^{-2\pi t^2}\sum_{i=0}^v(-4\pi)^{-i}\binom{2v}{2i}(2i-1)!!\,t^{2(v-i)}$.
\qed

\begin{lem}\label{le:localarch}
    Let $r\geq l\geq0$ be integers, $F$ the lift of a weight $2\ka$ eigenform
    and $\xi_\infty$ the character 
    $\xi_\infty(z)=(z/\bar{z})^{\ka+r}$ of $\C^\times$. Then
    $$
    L_\infty(\varphi^{(l,2(\ka+r))},\xi_\infty,s)=
    \begin{cases}
    0  & \text{if $l<r$}, \\
    \frac{(-1)^{r}2\pi(-{\rm N}\inv u)^\frac12\xi_\infty(z_u)r!}{(4\pi)^{s+2(\ka+r)-\frac12}}\Ga(s+2\ka+r-\frac12) & 
    \text{if $l=r$}.
\end{cases}
    $$
\end{lem}

\pf
It is well known that 
$W_F^{\psi_\infty}\left(\smallmat a001r(\th)\right)=a^\ka{\rm ch}_{\R^+ }(a)e^{-2\pi a-2\ka i\th}$.
We compute the other two terms in the integrand of \eqref{eq:localterm} separately with
$\varphi_1(z_1)=(z_1\bar{z}_1)^le^{-2\pi z_1\bar{z}_1}$ and 
$\varphi_2(z_2)=\bar{z}_2^{2(\ka+r)}e^{-2\pi z_2\bar{z}_2}$.
\begin{enumerate}
  \item To compute $\Phi^s_{\varphi_1}\left(\smallmat a001 r(\th)\right)=
            \vass{a}^sr_{\psi_\infty}(r(-\th))\varphi_1(0)$
            we use the definitions \eqref{eq:Weilrep} together with the decomposition
            \begin{equation}
            r(\th)=\smallmat{1}{-\tan\th}{0}{1}\smallmat{0}{-1}{1}{0}
            \smallmat{1}{-\sin\th\cos\th}{0}{1}\smallmat{0}{1}{-1}{0}\smallmat{1/\cos\th}{0}{0}{\cos\th}.
            \label{eq:decrth}
            \end{equation}
            Some straightforward passages yield 
            $r_{\psi_\infty}(r(\th))\varphi_1(0)=(-\cos\th)\varphi_1^\sharp(0)$ where
            $\varphi_1^\sharp(z)$ is the Fourier transform of 
            $e^{-2\pi i\cos\th\sin\th\vass{z}^2}\hat\varphi_1((\cos\th)z)$.     
            Since
            \begin{align*}
            \varphi_1^\sharp(0) &= \int_{\C}e^{-2\pi i\sin\th\cos\th\vass{z}^2}{\hat\varphi}_1((\cos\th)z)\,dz\\
            \intertext{(for $z=x+yi$ and from the previous lemma)}
                                           &= 2\sum_{0\leq\al+\be\leq l}\ga(\al,\be;l)(cos\th)^{2(\al+\be)}
                                                 \int_{\R^2}e^{-2\pi(i\sin\th\cos+(\cos\th)^2)(x^2+y^2)}x^{2\al}y^{2\be}\,dxdy\\
                                           &=-\sum_{0\leq\al+\be\leq l}\ga(\al,\be;l)\frac{(2\al-1)!!\, (2\be-1)!!}{(4\pi)^{\al+\be}}
                                                 \frac{(\cos\th)^{2(\al+\be)}}{(-\sin\th\cos\th-(\cos\th)^2)^{\al+\be+1}},
            \end{align*}
            eventually
            \begin{multline*}
            \Phi^s_{\varphi_1}\left(\smallmat a001 r(\th)\right)=\\
            -\vass{a}^s\sum_{0\leq\al+\be\leq l}\frac{\ga(\al,\be;l)(2\al-1)!!(2\be-1)!!}{(-4\pi)^{\al+\be}}
            (\cos\th)^{\al+\be}e^{-(\al+\be+1)i\th}.
            \end{multline*}
  \item To compute $W^{\psi_\infty}_{\theta_{\varphi_{2},\infty}}\left(\smallmat{-a}001r(\th)\right)$ we need 
            to use again
            \eqref{eq:Weilrep} together with the decomposition \eqref{eq:decrth}. 
            For, it should be noted that this time the 
            norm in $(Ku)\otimes\R\simeq\C$ is $-{\rm N}_{\C/\R}$ (in particular, definite negative)
            and the main involution is $z\mapsto -z$. Thus, we get
            $W^{\psi_\infty}_{\theta_{\varphi_{2},\infty}}\left(\smallmat{-a}001r(\th)\right)=
            e^{i(2\ka+2r+1)\th}W^{\psi_\infty}_{\theta_{\varphi_{2},\infty}}\left(\smallmat{-a}001\right)$.
            On the other hand, for a choice of $h\in\C^\times$ such that ${\rm N}h=-a{\rm N}\inv u>0$,
            \begin{align*}
            W^{\psi_\infty}_{\theta_{\varphi_{2},\infty}}\left(\smallmat{-a}001\right)
                &= \frac1{2\pi}\int_{S^1}r_{\psi_\infty}\left(\smallmat{-a{\rm N}\inv u}001,h\vth\right)
                      \varphi_2(u)\xi_\infty(h\vth)\,d\vth\\
                &= \frac{(-a{\rm N}\inv u)^\frac12}{2\pi}\int_{S^1}\varphi_2(-a{\rm N}\inv u\inv{(h\vth)}u))
                      \xi_\infty(h\vth)\,d\vth\\
                &= \frac{(-a{\rm N}\inv u)^\frac12}{2\pi}\int_{S^1}\varphi_2(\bar{h}{\inv\vth}u)\xi_\infty(h\vth)\,d\vth\\
                &= \frac{(-a{\rm N}\inv u)^\frac12}{2\pi}\int_{S^1}(\bar{h}\inv\vth{z}_u)^{2(\ka+r)}e^{-2\pi a}
                      (h\vth)^{\ka+r}({\bar h}\inv\th)^{-\ka-r}\,d\vth\\
                &= (-{\rm N}\inv u)^\frac12\xi_\infty(z_u)a^{\ka+r+\frac12}e^{-2\pi a}
            \end{align*}    
\end{enumerate}
Putting all the ingredients together
\begin{align*}
L_\infty&(\varphi^{(l,2(\ka+r))},\xi_\infty,s)=\\ 
             &=- (-{\rm N}\inv u)^\frac12\xi_\infty(z_u)\int_{\R^{>0}}\int_{S^1}
                                                                    a^{s+2\ka+r-\frac12}e^{-4\pi a}e^{i(2\ka+2r+1)\th}\\
             & \qquad\qquad\times\left(\sum_{0\leq\al+\be\leq l}\frac{\ga(\al,\be;l)(2\al-1)!!(2\be-1)!!}{(-4\pi)^{\al+\be}}
                                                                   (\cos\th)^{\al+\be}e^{-(\al+\be+1)i\th}\right)\,d^\times ad\th\\
             & =-(-{\rm N}\inv u)^\frac12\xi_\infty(z_u)\int_{\R^{>0}}a^{s+2\ka+r-\frac12}e^{-4\pi a}\,d^\times a\\
             & \qquad\qquad\times\sum_{0\leq\al+\be\leq l}\frac{\ga(\al,\be;l)(2\al-1)!!(2\be-1)!!}{(-4\pi)^{\al+\be}}
                                                    \int_{S^1}(\cos\th)^{\al+\be}e^{-(\al+\be+1)i\th}\,d\th
\end{align*}
Since $\al+\be\leq l\leq r$ we have
\begin{multline*}
\int_{S^1}(\cos\th)^{\al+\be}e^{i(2r-\al-\be+1)\th}\,d\th=\\
\frac1{2^{\al+\be}}\sum_{j=0}^{\al+\be}\binom{\al+\be}{j}\int_{S^1}e^{2i(r-j)\th}\,d\th=
\begin{cases}
   2^{1-r}\pi   & \text{if $\al+\be=l=r$ }, \\
    0  & \text{otherwise}.
\end{cases}
\end{multline*}
hence $L_\infty(\varphi^{(l,2(\ka+r))},\xi_\infty,s)=0$ if $l<r$. When $l=r$, since $\ga(j,r-j;r)=\binom rj$ and
$\sum_{j=0}^r\binom rj(2j-1)!!(2r-2j-1)!!=2^rr!$ as readily proved by induction, we have
\begin{align*}
L_\infty(\varphi^{(r,2(\ka+r))},\xi_\infty,s)&=\frac{2\pi(-{\rm N}\inv u)^\frac12\xi_\infty(z_u)r!}{(-4\pi)^r}
                                                                   \int_{\R^{>0}}a^{s+2\ka+r-\frac12}e^{-4\pi a}\,d^\times a\\
 &=\frac{(-1)^r2\pi(-{\rm N}\inv u)^\frac12\xi_\infty(z_u)r!}{(4\pi)^{s+2(\ka+r)-\frac12}}\Ga(s+2\ka+r-\frac12).
\end{align*}
\qed

We shall now state and prove the main result of this section.

\begin{thm}\label{thm:maininterpolation}
    Let $N$ be a positive integer and fix a decomposition 
    $N=\De N_{o}$ with $\De$ a product of an even number of distinct 
    primes and $(\De,N_{o})=1$.
    Let $\pi$ be an automorphic cuspidal representation for $\GL_{2}$ of 
    conductor $N$ such that
    \begin{enumerate}
        \item  $\pi_{\infty}\simeq\si(\mu_{1},\mu_{2})$, the discrete series 
               representation with $\mu_{1}\mu_{2}^{-1}(t)=t^{2\ka-1}{\rm sgn}(t)$.
        \item  $\pi_{\ell}$ is special for each $\ell|\De$. 
    \end{enumerate}
    Let $K$ be a quadratic imaginary field such that all $\ell|\De$ 
    are inert in $K$ and all $\ell|N_{o}$ are split in $K$. Let $c$ be a positive integer
     with $(c,N)=1$ and $p$ an odd prime number not dividing
    $N$ that splits in  $K$. Assume that $\calO_{K,c}^{\times}=\{\pm1\}$.
    Suppose that there exist Gr\"ossencharakters $\chi\in\Xi_{(-2\ka,0)}(\calO_{K,c})$
    and $\xi\in\Xi_{(-2,0)}(\calO_{K,c})$ such that the $p$-adic avatar $\wh\xi$ takes values 
    in a totally ramified extension of $\Q_p$.
    
    Then, there exists $x\in\CM(\De,N;\calO_{K,c})$ represented by
    $\tau=t+yi\in\CM_{\De,K}$ with associated periods $\Om_\infty$ and $\Om_p$
    such that for all $r\geq0$
    \begin{multline*}
    \Om_{p}^{-4(\ka+r)}\int_{\Z_p}z^r\,d\mu(f,x;\chi,\xi)= \\
    \frac{2\varpi V_Nw_{K,c}^2}{m_{c}^2h_K}
    \frac{(-1)^{\ka+r}r!(2\ka+r)!}
    {4^{2\ka+3r}\pi^{2(\ka+r+1)}y^{2(\ka+r)}\Om_\infty^{4(\ka+r)}}
    \nu(\xi_r,\tau,\frac12)L(\pi_K\otimes\xi_r,\frac12)L(\eta_K,1)^{-1}
    \end{multline*}
    where $\xi_r=\xi\chi^r\vvass{N_{K/\Q}}^{-\ka-r}$ and $\varpi$ is a (fixed) ratio of Petersson norms.
\end{thm}

\pf
Let $D$ be the quaternion algebra with $\De_{D}=\De$. By hypothesis 
the representation $\pi$ is the 
image of an automorphic representation $\pi^{D}$ of $D^{\times}$ under 
the Jacquet-Langlands correspondence and let 
$f\in S_{2k,0}(\De,N_{o})$ be a holomorphic newform in $\pi^{D}$.
For all integers $r\geq0$ let $\phi_{r}$ be as in \eqref{eq:defphir}. 

By proposition \ref{teo:existCM} there exists $x\in\CM(\De,N;\calO_{K,c})$ 
and choose a split $p$-ordinary test triple $(\tau,v,e)$, $\tau=t+iy$, 
representing $x$ with corresponding $d_\infty\in\SL_2(\R)$.
By taking $\calO_{v}$ and $F$ large enough, we can 
assume that $f$ is defined over $\calO_{(v)}\subset\calO_F$, and that the measure 
$\mu(f,x;\chi,\xi)$ has values in $\calO_{F}$

By remark \ref{rm:onchoice}.3 we can write 
$\pi(d_{\infty})\phi_{0}\otimes\pi(d_{\infty})\phi_{0}=\varpi\widetilde\theta_{\varphi}(F)$
with $\varphi=\varphi_\infty\otimes\varphi^f$ as in \eqref{eq:choiceofphi},
$F$ the adelization of the normalized eigenform in $\pi$
and $\varpi\in\C^\times$ a Petersson normalization constant.
We claim that for all $r\geq0$
\begin{equation}
    \pi(d_{\infty})\phi_{r}\otimes\pi(d_{\infty})\phi_{r}=
    \varpi\frac{(-1)^\ka}{4^r\pi}
    \widetilde\theta_{\phi^{r,2(\ka+r)}\otimes\varphi^{f}}(F)
    +\sum_{l=0}^{r-1}a_{r,l},\widetilde\theta_{\phi^{l,2(\ka+r)}\otimes\varphi^{f}}(F)
    \label{eq:thetaclaim}
\end{equation}
where $a_{r,l}\in\varpi\Z[\pi]$. 
For, the short exact sequence \eqref{eq:sesGOD} 
gives a Lie algebras identification
$\gotg\gotoo(D)\simeq(D_{\infty}\times D_{\infty})/\R$ 
and in particular
$\gotoo(D)=\{(A,B)\in D_{\infty}\times D_{\infty}\,|\,\tra A=\tra B\}/\R
\simeq\gots\gotl_{2}\times\gots\gotl_{2}$.
Under this identification, differentiating \eqref{eq:tlautom} yields
$$
\widetilde\theta_{H\varphi}(F)=\left.\frac{d}{dt}
\widetilde\theta_{\varphi}(F)(h\exp(tH))\right|_{t=0} 
\hbox{ with }
H\varphi(x)=\left.\frac{d}{dt}\varphi(e^{-tH_{1}}xe^{tH_{2}})\right|_{t=0}
$$
for all $H=(H_{1},H_{2})\in\Lie(\mathrm{O}(D))$. If 
$A\in\gots\gotl_{2}$ a repeated application of the last formula with 
$\pr{A}=(A,0)$ and $\pr{A}{}^{\prime}=(0,A)$ shows that the diagonal action 
of $A$ on $\pi^{D}\otimes\pi^{D}$ corresponds to the action 
of the second order operator 
$A_{2}=\pr{A}\pr{A}{}^{\prime}=\pr{A}{}^{\prime}\pr{A}
\in\gotA(\Lie(\mathrm{O}(D)))$ on Schwartz functions, i.e.
$$
    A_{2}\varphi(x)=\left.\frac{\partial^{2}}{\partial u\partial 
    v}\varphi(e^{-uA}xe^{vA})\right|_{u=v=0}.
$$
We are interested in the expression of the operator $A_{2}$ in the 
normalized coordinates for 
$A=d_{\infty}X^{+}d_{\infty}^{-1}$. Up to conjugation, 
this is the same as to compute the second order operator associated 
to $A=X^{+}$ under the standard coordinates \eqref{eq:standardcoord}.
A straightforward computation using the obvious real coordinates associated to the underlying
real decomposition
$D_{\infty}= \R\smallmat 1{}{}1\oplus\R\smallmat{}{-1}1{}\oplus
\R\smallmat {}11{}\oplus\R\smallmat{-1}{}{}1$
shows that
$\pr{A}=-i\left(z_2\frac{\partial}{\partial{\bar z}_1}+z_1\frac{\partial}{\partial{\bar z}_2}\right)$
and
$\pr{A}{}^\prime=i\left(z_2\frac{\partial}{\partial z_1}+{\bar z}_1\frac{\partial}{\partial{\bar z}_2}\right)$, so that
$$
A_2=z_2^2\frac{\partial^2}{\partial z_1\partial{\bar z}_1}+
{\bar z}_1{z}_2\frac{\partial^2}{\partial{\bar z}_1\partial{\bar z}_2}+
z_1z_2\frac{\partial^2}{\partial z_1\partial{\bar z}_2}+
z_1{\bar z}_1\frac{\partial^2}{{\partial{\bar z}_2}^2}+
z_2\frac{\partial}{\partial{\bar z}_2}.
$$
Since
$$
A_2\phi^{m,q}=
\begin{cases}
 -2\pi\phi^{0,q+2}+4\pi^2\phi^{1,q+2}     & \text{if $m=0$}, \\
 m^2\phi^{m-1,q+2}-(4m-2)\pi\phi^{m,q+2}+4\pi^2\phi^{m+1,q+2}     & \text{if $m\geq1$},
\end{cases}
$$
formula \eqref{eq:thetaclaim} follows from an $r$-fold iteration using the linearity of the theta lift and the definitions \eqref{eq:defphir} and \eqref{eq:choiceofphi} of $\phi_r$ and $\varphi_\infty$ respectively.

Let $\chi_r$ be a Gr\"ossencharakter of weight $(-2(\ka+r),0)$ 
and trivial on $\wh{R}_{c}^{\times}$ such that 
$\xi_r=\chi_r\vvass{N_{K/\Q}}^{-\ka-r}$ is unitary. 
Combining \eqref{eq:thetaclaim} and \eqref{eq:finalexp} with lemma \ref{le:localarch} we get
\begin{equation}
\label{eq:end}
J_{r}(f,\xi_r,\tau)^{2} = 
\frac{(-1)^{\ka+r}2\varpi V_Nr!(2\ka+r)!}{4^{2\ka+3r}h_K\pi^{2\ka+2r}}
\nu(\xi_r,\tau,\frac12)L(\pi_K\otimes\xi_r,\frac12)L(\eta_K,1)^{-1}.
\end{equation}
On the other hand, from theorem \ref{th:compJ},
$$
J_{r}(f,\xi_{r},\tau)^{2}=
\frac{m_c^2\pi^2}{w_{K,c}^2}(h_c^\sharp)^2y^{2(\ka+r)}\Om_\infty^{4(\ka+r)}
\scal{c_{(r)}(f,x)}{\chi_r}^{2}.
$$
When $\chi_r=\chi\xi^r$ we use proposition \ref{th:spiden} to rewrite the last formula as
$$
\Om_p^{-4(\ka+r)}m_r(\mu(f,x;\chi,\xi))^2=
\frac{w_{K,c}^2}{m_c^2\pi^2}\Om_\infty^{-4(\ka+r)}y^{-2(\ka+r)}J_{r}(f,\xi_r,\tau)^{2}.
$$
Substituting \eqref{eq:end} into the latter formula proves the theorem. \qed



\begin{thebibliography}{10}
\bibliographystyle{plain}
    
\bibitem{BerDar96}
M.~Bertolini and H.~Darmon.
\newblock Heegner points on {M}umford-{T}ate curves.
\newblock {\em Inventiones Math.}, 126:413--456, 1996.

\bibitem{BoCa91}
J.~F. Boutot and H.~Carayol.
\newblock Uniformisation $p$-adique des courbes de {S}himura: les theoremes de
  {C}erednik et de {D}rinfeld.
\newblock In {\em Courbes modulaires at courbes de {S}himura}, volume 196 of
  {\em Ast{\'e}risque}, pages 45--158, 1991.
  
\bibitem{Bu96}
D.~Bump.
\newblock {\em Automorphic {F}orms and {R}epresentations}, volume~55 of {\em
  Cambridge Studies in Adv. Math..}
\newblock Cambridge University Press, 1996.

\bibitem{Dar04}
H.~Darmon.
\newblock {\em Rational {P}ooints on {M}odular {E}lliptic {C}urves}, volume~101 of {\em
  Regional Conf. Ser. Math..}
\newblock American Mathematical Society, 2004.

\bibitem{Del82}
P.~Deligne.
\newblock Hodge cycles on abelian varieties.
\newblock In {\em Hodge cycles, {M}otives, and {S}himura varieties}, volume 900
  of {\em Lecture Notes in Math.}, pages 9--100, 1982.

\bibitem{DelRap73}
P.~Deligne and M.~Rapoport.
\newblock Les sch{\'e}mas de modules de courbes elliptiques.
\newblock In {\em Modular Functions of One Variable II}, volume 349 of {\em
  Lecture Notes in Math.}, pages 143--316, 1973.

\bibitem{DiaIm95}
F.~Diamond and J.~Im.
\newblock Modular forms and modular curves.
\newblock In {\em Seminar on Fermat's Last Theorem (Toronto, ON, 1993--1994)},
  volume~17 of {\em CMS Conf. Proc.}, pages 39--133. American Math. Soc., 1995.

\bibitem{DiaTay94}
F.~Diamond and R.~Taylor.
\newblock Non-optimal levels of mod $l$ representations.
\newblock {\em Inventiones Math.}, 115:435--462, 1994.

\bibitem{FalCha90}
G.~Faltings and C.-L. Chai.
\newblock {\em Degenerations of Abelian Varieties}, volume~22 of {\em
  Ergebnisse der Math.}
\newblock Springer, 1990.

\bibitem{Harris79}
M.~Harris.
\newblock A note on three lemmas of {S}himura.
\newblock {\em Duke Math. J.}, 46:871--879, 1979.

\bibitem{Harris81}
M.~Harris.
\newblock Special values of zeta functions attached to {S}iegel modular forms.
\newblock {\em Ann. Sc. {\'E}c. Norm. Sup.}, 14:77--120, 1981.

\bibitem{Harris93}
M.~Harris.
\newblock Non-vanishing of ${L}$-functions on $2\times2$ unitary groups.
\newblock {\em Forum Math.}, 5:405--419, 1993.

\bibitem{HaKu91}
M.~Harris and S.~Kudla.
\newblock The central critical value of a triple product ${L}$-function.
\newblock {\em Annals Math.}, 133:605--672, 1991.

\bibitem{HaKu92}
M.~Harris and S.~Kudla.
\newblock Arithmetic automorphic forms for the nonholomorphic discrete series
  of ${\GSp}(2)$.
\newblock {\em Duke Math. J.}, 66:59--121, 1992.

\bibitem{HaTi01}
M.~Harris and J.~Tilouine.
\newblock $p$-adic measures and square roots of special values of triple
  product ${L}$-functions.
\newblock {\em Math. Ann.}, 320:127--147, 2001.

\bibitem{Hashim95}
K.~Hashimoto.
\newblock Explicit form of quaternion modular embeddings.
\newblock {\em Osaka J. Math.}, 32:533--546, 1995.

\bibitem{Hida86}
H.~Hida.
\newblock Hecke algebras for ${\GL}_1$ and ${\GL}_2$.
\newblock In {\em Seminaire de theorie des nombres, Paris 1983--84}, volume~63
  of {\em Progr. Math.}, pages 133--146, 1986.

\bibitem{Hida93}
H.~Hida.
\newblock {\em Elementary Theory of ${L}$-functions and Eisenstein Series},
  volume~26 of {\em Student Texts}.
\newblock London Mathematical Society, 1993.

\bibitem{Howe89}
R.~Howe.
\newblock Transcending classical invariant theory.
\newblock {\em J. Amer. Math. Soc}, 2:536--552, 1989.

\bibitem{Ja72}
H.~Jacquet.
\newblock {\em Automorphic forms on $\GL(2)$, II}, volume 278 of {\em Lecture Notes
  in Math}.
\newblock Springer-Verlag, 1972.

\bibitem{JaLa70}
H.~Jacquet and R.~P. Langlands.
\newblock {\em Automorphic forms on $\GL(2)$}, volume 114 of {\em Lecture Notes
  in Math}.
\newblock Springer-Verlag, 1970.

\bibitem{Jord86}
B.~Jordan.
\newblock Points on {S}himura curves rational over number fields.
\newblock {\em J. reine u. angew. Math.}, 371:92--114, 1986.

\bibitem{Katz70}
N.~Katz.
\newblock Nilpotent connections and the monodromy theorem: applications of a
  result of {T}urrittin.
\newblock {\em Publ. Math. IHES}, 39:355--412, 1970.

\bibitem{Katz73}
N.~Katz.
\newblock $p$-adic properties of modular schemes and modular forms.
\newblock In {\em Modular Functions of One Variable III}, volume 350 of {\em
  Lecture Notes in Math.}, pages 70--189, 1973.

\bibitem{Katz76}
N.~Katz.
\newblock $p$-adic interpolation of real analytic {E}isenstein series.
\newblock {\em Annals of Math.}, 104:459--571, 1976.

\bibitem{Katz78}
N.~Katz.
\newblock $p$-adic ${L}$-function for {C}{M} fields.
\newblock {\em Inventiones Math.}, 49:199--297, 1978.

\bibitem{Katz81}
N.~Katz.
\newblock Serre-{T}ate local moduli.
\newblock In {\em Surfaces Algebraiques}, volume 868 of {\em Lecture Notes in
  Math.}, pages 138--202, 1981.

\bibitem{KatMaz85}
N.~Katz and B.~Mazur.
\newblock {\em Arithmetic moduli of elliptic curves}, volume 108 of {\em Annals
  of Math. Studies}.
\newblock Princeton Univ. Press, 1985.

\bibitem{KatOda68}
N.~Katz and T.~Oda.
\newblock On the differentiation of de {R}ham cohomology classes with respect
  to parameters.
\newblock {\em J. Math. Kyoto Univ.}, 8(2):199--213, 1968.

\bibitem{Kudla84}
S.~Kudla.
\newblock Seesaw dual reductive pairs.
\newblock In {\em Automorphic forms of several variables (Katata 1983)},
  volume~46 of {\em Progr. Math}, pages 244--268. Birkh\"{a}user Boston, 1984.

\bibitem{Maass53}
H.~Maass.
\newblock Die {D}ifferentialgleichungen in der {T}heorie der {S}iegelschen
  {M}odulfunktionen.
\newblock {\em Math. Ann.}, 126:44--68, 1953.

\bibitem{Milne79}
J.~S. Milne.
\newblock Points on {S}himura varieties $\bmod p$.
\newblock In {\em Automorphic {F}orms, {R}epresentations and ${L}$-functions},
  volume~33 of {\em Proc. Symp. Pure Math.}, pages 165--183, 1979.

\bibitem{Mori94}
A.~Mori.
\newblock A characterization of integral elliptic automorphic forms.
\newblock {\em Ann. Sc. Norm. Sup. Pisa}, IV(21):45--62, 1994.

\bibitem{MoTe99}
A.~Mori and L.~Terracini.
\newblock A canonical map between {H}ecke algebras.
\newblock {\em Boll. Un. Mat. It. Sez. B}, 8(2):429--452, 1999.

\bibitem{Mum65}
D.~Mumford.
\newblock {\em Geometric Invariant Theory}.
\newblock Springer-Verlag, Heidelberg, 1965.

\bibitem{Noot92}
R.~Noot.
\newblock {\em Hodge classes, {T}ate classes, and local moduli of abelian
  varieties}.
\newblock PhD thesis, Rijksuniversiteit Utrecht, 1992.

\bibitem{Oort71}
F.~Oort.
\newblock Finite group schemes, local moduli for abelian varieties and lifting
  problems.
\newblock {\em Compositio Math.}, 23:265--296, 1971.

\bibitem{Pra06}
K.~Prasanna.
\newblock {\em Integrality of  a ratio of {P}etersson norms and level-lowering congruences}.
\newblock {\em Ann. Math.}, 163:901--967, 2006.

\bibitem{Robert89}
D.~Roberts.
\newblock {\em Shimura curves analogous to ${X}_0({N})$}.
\newblock PhD thesis, Harvard University, 1989.

\bibitem{Serre62}
J.-P. Serre.
\newblock Endomorphismes completement continus des espaces de {B}anach
  $p$-adiques.
\newblock {\em IHES Publ. Math.}, 12:69--85, 1962.

\bibitem{Serre73}
J.-P. Serre.
\newblock Formes modulaires et fonctions z\^eta $p$-adiques.
\newblock In {\em Modular Functions of One Variable, III}, volume 350 of {\em Lecture Notes in
  Math.}, pages 191--268, 1973.

\bibitem{Shimi72}
H.~Shimizu.
\newblock Theta series and automorphic forms on ${\GL}_2$.
\newblock {\em Math. Soc. Japan}, 24:638--683, 1972.

\bibitem{ShiRed}
G.~Shimura.
\newblock {\em Introduction to the Arithmetic Theory of Automorphic Functions}.
\newblock Iwanami Shoten and Princeton University Press, 1971.

\bibitem{Tilo96}
J.~Tilouine.
\newblock {\em Deformations of {G}alois representations and {H}ecke algebras}.
\newblock Narosa Publ. House, 1996.

\bibitem{Vigner80}
M.-F. Vigneras.
\newblock {\em Arithm{\'e}tique des Alg{\`e}bres de Quaternions}, volume 800 of
  {\em Lecture Notes in Math.}
\newblock Springer-Verlag, 1980.

\bibitem{Waldsp85}
J.-L. Waldspurger.
\newblock Sur les valeurs des certaines fonctions ${L}$ automorphes en leur
  centre de sym{\'e}trie.
\newblock {\em Compositio Math.}, 54:173--242, 1985.

\bibitem{Wat03}
T.~Watson.
\newblock {\em Rankin triple products and quantum chaos }.
\newblock PhD thesis, Princeton University, 2003.

\bibitem{Weil55}
A.~Weil.
\newblock On a certain type of characters of the idele-class group of an
  algebraic number-field.
\newblock In {\em Proceedings of the international symposium on algebraic
  number theory, Tokio \& Nikko 1955}, pages 1--7, 1956.

\end{thebibliography}
\end{document}